\documentclass[reqno]{amsart}

\usepackage[utf8]{inputenc}
\usepackage[T1]{fontenc}
\usepackage{ae,aecompl}

\usepackage{amsmath,amsfonts,amssymb}%,amsthm
\usepackage{color}
\usepackage{dsfont} % for characteristic function
\usepackage{bm}

\usepackage{graphicx}
\usepackage{tikz}
\usepackage[european]{circuitikz}
\usepackage{pgfplots}
\usepgfplotslibrary{patchplots}
\usepgfplotslibrary{groupplots}
\pgfplotsset{compat=newest}
\usetikzlibrary{calc}

\usepackage{fancyhdr}
\usepackage{footnote}
\usepackage{nicefrac} % for other frac
\usepackage{setspace}
\usepackage{ifthen}

\newtheorem{theorem}{Theorem}[section]
\newtheorem{proposition}[theorem]{Proposition}
\newtheorem{lemma}[theorem]{Lemma}
\newtheorem{corollary}[theorem]{Corollary}
\newtheorem{definition}[theorem]{Definition}
\newtheorem{remark}[theorem]{Remark}
\newtheorem{assumption}[theorem]{Assumption}
\newtheorem{example}{Example}

\numberwithin{equation}{section}
\newcommand{\myempty}{\hspace*{2.9ex}\mbox{}}
\newcommand{\nullsp}{\mathcal{N}}

\begin{document}

\title{Splitting Techniques for DAEs with port-Hamiltonian Applications}

\author[Bartel et al.]{Andreas Bartel$^{1,\ddag}$ \and \mbox{Malak Diab}$^{1}$ \and \mbox{Andreas Frommer}$^{1}$ \and 
\mbox{Michael Günther}$^{1}$ \and \mbox{Nicole Marheineke}$^{2}$}

\date{\today\\
$^1$ IMACM, University of Wuppertal, Gaußstraße 20, D-42097 Wuppertal, Germany, email: \{bartel, mdiab, frommer, guenther\}@uni-wuppertal.de \\
$^2$ Trier University, Department IV, Universitätsring 15, D-54296 Trier, Germany; email: marheineke@uni-trier.de\\
$^\ddag$ corresponding author, orcid: 0000-0003-1979-179X}

\begin{abstract}
In the simulation of differential-algebraic equations (DAEs), it is essential to employ numerical schemes that take into account the inherent structure and maintain explicit or hidden algebraic constraints without altering them. This paper focuses on operator-splitting techniques for coupled systems and aims at preserving the structure in the port-Hamiltonian framework. The study explores two decomposition strategies: one considering the underlying coupled subsystem structure and the other addressing energy-associated properties such as conservation and dissipation. We show that for coupled index-$1$ DAEs with and without private index-2 variables, the splitting schemes on top of a dimension-reducing decomposition achieve the same convergence rate as in the case of ordinary differential equations. Additionally, we discuss an energy-associated decomposition for index-1 pH-DAEs and introduce generalized Cayley transforms to uphold energy conservation. The effectiveness of both strategies is evaluated using port-Hamiltonian benchmark examples from electric circuits. 
\end{abstract}

\maketitle

\smallskip


%-----------------------------------------------------
\section{Introduction}
\label{intro}
%-----------------------------------------------------
 
Port-Hamiltonian differential-algebraic systems (pH-DAEs) arise from port-based network modeling of multi-physics problems. For this, a physical system is decomposed into smaller subsystems that are interconnected through energy exchange. The subsystems may belong to various different physical domains, e.g. electrical, mechanical, or hydraulic ones. The energy-based formulation is advantageous since the physical properties, such as energy-conservation and -dissipation, are directly encoded in the structure of the port-Hamiltonian model equations, the port-Hamiltonian character is inherited by the coupling, and different scales are brought on a single level. Algebraic equations naturally come from the interconnections in the form of network conditions, such as Kirchhoff’s laws in electrical circuits, or from constraints that are directly modeled, like, e.g., position or velocity constraints in mechanical systems, or mass balances in chemical engineering problems, see, e.g., \cite{IVP_DAE95,kunkel2006,riaza2008}.

Port-Hamiltonian systems can be derived in two different ways, via a formulation as descriptor systems with special structured coefficient matrices \cite{beattie2018} or via an energy-based formulation on top of a Dirac structure \cite{schaft2013}.
The properties of pH-DAEs have been studied in, e.g., \cite{beattie2018,schaft2013,schaft2018}.
The systems typically consist of explicit as well as implicit (hidden) constraints. For simulation and optimization, numerical schemes are required that are aware of the structure and keep the algebraic constraints unchanged so as not to destroy crucial properties \cite{kotyczka2018}. If the differentiation-index is larger than one, an index reduction is possible, e.g., via derivative arrays or minimal extension, cf.\ \cite{kunkel2006} for general DAEs. For pH-DAEs the index reduction has to be performed in a structure-preserving way, \cite{beattie2018}. The differentiation-index has been shown to be at most two, see \cite{mehl2018} for the linear constant coefficient case. 

Operator splitting is a powerful numerical tool to deal with dynamical systems
\cite{Blanes2010,mclachlan2002}. The numerical procedure consists of three ingredients/steps: (1) decomposition of the right hand into subproblems of profoundly different behavior regarding, e.g., dynamics, stiffness, or into subproblems of smaller size to enhance the computational efficiency,
(2) splitting scheme for the exact subproblems' fluxes, and 
(3) numerical flux approximation. While the decomposition (first step) is problem-dependent, the actual splitting (second step) is general. Splitting schemes differ in the sequence of and the step size of the respective exact subproblems' fluxes. Convergence theory is available in case of ordinary differential equations (ODE), e.g., via the Baker-Campbell-Hausdorff formula, \cite{Blanes2008}. Certainly, the most prominent schemes are the Lie-Trotter splitting of first order and the symmetric Strang splitting of second order \cite{Strang1968}. The numerical flux approxi\-mation/integrator (third step) is finally chosen or tailored with respect to the properties of the respective subproblem and the convergence order of the underlying splitting scheme.

Operator splitting for port-Hamiltonian applications is a topic of \cite{ECMI_2023,Frommer_2023aa}.
This article aims at a generalization of the latter mentioned works. It deals with operator splitting techniques that are structure-aware or even structure-preserving. A crucial point is the handling of the algebraic constraints. In view of the port-Hamiltonian framework, we propose and investigate two different decomposition strategies: A) the \textit{dimension-reducing decomposition} accounting for an underlying coupled subsystem structure and B) the \textit{$J$-$R$ decomposition} accounting for the energy-associated properties, such as energy-conservation and dissipation. 
Presupposing coupled index-1 DAEs---without and with private index-2 variables---as they typically arise from network mo\-deling, we show that a clever doubling of the constraints in the dimension-reducing decomposition allows to treat the single subsystems as inherent ODEs. As consequence, the well-established convergence theory on splitting schemes for ODEs can be transferred to the DAE case. In particular, we analytically establish the second-order convergence rate for the Strang scheme. The numerical integration of the subsystems deserves special attention, since the constraints must not be violated here, as we will comment on. The $J$-$R$ decomposition is generally applicable to implicit pH-ODEs, but the interaction with constraints has a crucial influence on the applicability for 
index-1 pH-DAEs. We show that for certain assignments of the energy parts in the constraint, the $J$-$R$ decomposition leads to two subproblems of an inherent ODE and a lower dimensional ODE. In this case, convergence results and properties of the splitting schemes can be deduced from the splitting theory for ODEs. 
Schemes of Lie-Trotter-type and Strang-type satisfy the characteristic port-Hamiltonian dissipation equality, whereas the negative step sizes occurring in splitting schemes of order $p\geq 3$ destroy the dissipativity, cf.\ also \cite{hairer2006,suzuki1991}. For the Strang-type splitting, we discuss the numerical flux discretization and introduce a new concept of generalized Cayley transforms targeted at energy-conservation.
The performance of the two decomposition strategies in combination with the second-order Strang-type splittings is numerically investigated for port-Hamiltonian benchmark examples originating from electric circuits.

The paper is structured as follows:  
In Section~\ref{sec:phs-framework}, we give a short review on the essentials of the port-Hamiltonian DAE framework and discuss some details which are relevant in our setting. Section~\ref{sec:coupled} and Section~\ref{sec:energy-based-splitting} then deal with the dimension-reducing decomposition and the $J$-$R$ decomposition, respectively, and establish splitting convergence properties for DAEs. 
We demonstrate our findings with numerical examples, mainly from electric circuit models, in  Section~\ref{sec:numerics} and finally present our conclusions in Section~\ref{sec:conclusion}. Details on the generalized Cayley transforms, on splitting schemes and flux approximations as well as on the modeling of electrical networks are reported in the appendices.

%------------------------------------------------------------
\section{Model Framework of pH-DAEs}
\label{sec:phs-framework}
%------------------------------------------------------------

In this section, we summarize structural properties and simplified representations of pH-DAEs. For details we refer to \cite{beattie2018,schaft2018}.

\begin{definition}[pH-DAE] \label{def:pHS-DAE} A linear constant coefficient DAE system of the form
\begin{equation}\label{eq:linear_PHS}
\begin{aligned}
E\dot{x}&=(J-R)\,Qx+(B-P)\,u,\\
y&=(B+P)^\top\,Qx+(S+N)\,u,
\end{aligned}
\end{equation}
with $E$, $Q$, $J$, $R\in\mathbb{R}^{n\times n}$, $B$, $P\in\mathbb{R}^{n\times m}$, $S=S^\top$, $N=-N^\top\in\mathbb{R}^{m\times m}$, and $x = x(t) \in \mathbb{R}^n, u= u(t), y = y(t) \in \mathbb{R}^m$ with $t \in \mathbb I\subset \mathbb{R}$, $\mathbb I=[0,T]$ a compact interval, is called a port-Hamiltonian differential-algebraic system (pH-DAE) if the following properties are satisfied:
\begin{enumerate}
\item the differential-algebraic operator
\begin{align*}
Q^\top E\frac{d}{dt}-Q^\top JQ : \mathcal{X}\subset \mathcal{C}^1(\mathbb{I},\mathbb{R}^n)\rightarrow \mathcal{C}^0(\mathbb{I},\mathbb{R}^n)
\end{align*}
is skew-adjoint, i.e.,
$Q^\top J^\top Q=-Q^\top JQ$ and $Q^\top E=E^\top Q$ hold,
\item the product $Q^\top E$ is positive semidefinite, denoted as $Q^\top E=E^\top Q\geq 0$, 
\item the passivity matrix
\begin{align*}
W=\begin{pmatrix}Q^\top RQ \;& \;Q^\top P\\ P^\top Q & S\end{pmatrix}\in\mathbb{R}^{(n+m)\times(n+m)}
\end{align*}
is symmetric positive semi-definite, $W=W^\top \geq 0$.
\end{enumerate}
The underlying quadratic Hamiltonian ${\mathcal H}:\mathbb{R}^n\rightarrow\mathbb{R}$ of the system is given by
\begin{align*}
{\mathcal H}(x)=\frac{1}{2}x^\top Q^\top Ex.
\end{align*}
\end{definition}

\begin{proposition}[\cite{beattie2018,schaft2018}, e.g.]\label{th:pHDAE} 
If the pH-DAE system \eqref{eq:linear_PHS} has a (classical) solution $x\in \mathcal{C}^1(\mathbb{I},\mathbb{R}^n)$ for given input function $u$, then
\begin{align} \label{eq:dissipation_inequality}
\frac{d}{dt}{\mathcal H}(x)=u^\top y-\begin{pmatrix}x\\u\end{pmatrix}^\top W\begin{pmatrix}x\\u\end{pmatrix}.
\end{align}
Furthermore, if $W=0$, then $\frac{d}{dt}{\mathcal H}(x)=u^\top y$.
\end{proposition}

\begin{remark}
The classical solution $x\in \mathcal{C}^1(\mathbb{I},\mathbb{R}^n)$ of \eqref{eq:linear_PHS} exists if the matrix pencil $\{E,(J-R)Q\}$ is regular, i.e., the null spaces of $E$ and $(J-R)Q$ have only trivial intersection, and the input satisfies $u\in \mathcal{C}^2(\mathbb{I},\mathbb{R}^m)$.
%\hfill $\Box$
\end{remark}

Proposition~\ref{th:pHDAE} implies some important properties of a pH-DAE. First of all, its Hamiltonian is an energy storage function, and the system is passive with the decrease in energy described by the dissipation inequality \eqref{eq:dissipation_inequality}. Furthermore, a pH-DAE is implicitly Lyapunov stable as $\mathcal H$ defines a Lyapunov function, cf.\ \cite{beattie2018}. The physical properties are encoded in the algebraic structure of the coefficient matrices.
In this sense, $E^\top Q$ is the energy matrix, $Q^\top RQ$ is the dissipation matrix, $Q^\top JQ$ is the structure matrix describing the energy flux among the energy storage elements, $B  \pm P$ are port matrices for energy in-  and output, and  $S$, $N$ are the matrices associated with a direct feed-through from input $u$ to output $y$. In the case that $E$ is the identity matrix, $E = I$, the pH-DAE reduces to a standard pH-ODE as studied, e.g., in \cite{schaft2014}.

From now on, we will always assume that $Q$ is non-singular. 

\begin{remark}\label{rem:QE-ID} Using linear transformations, we can reduce \eqref{def:pHS-DAE} to the case $Q=I$ and $E \geq0$  diagonal, i.e., to the semi-explicit representation   \eqref{eq:phs.ode.transformed} below.
Indeed, left-multiplying \eqref{def:pHS-DAE} with $Q^\top$ and then using the 
eigen-decomposition $Q^\top E = V^\top DV$ 
with $D=\mathrm{diag}(d_1,\ldots,d_r,0,\ldots,0)$, $d_i>0$, $V^\top V=I$, together with the transformation $\bar x= \tilde{D}^{1/2}V x$ where $\tilde{D}^{1/2}$ is the diagonal matrix  $\tilde{D}^{1/2}:=$diag$(\sqrt{d_1},\ldots,\sqrt{d_r},1,
\ldots,1)$, we obtain 
\begin{equation}\label{eq:phs.ode.transformed}
\begin{aligned}
\bar{E} \dot{\bar x}&=(\bar J- \bar R)\, \bar x+(\bar B- \bar P)\,u,\\
y&=(\bar B+\bar P)^\top\, \bar x+(S+N)\,u,
\end{aligned}
\end{equation}
where
\begin{align*}
    \bar E &= (\tilde{D}^{-1/2}V) Q^TE \,(\tilde{D}^{-1/2}V)^\top \,\,\,= \left( 
\begin{smallmatrix}
    I_r & 0 \\ 0 & 0
\end{smallmatrix} 
\right),  
    \\
    \bar J &= (\tilde{D}^{-1/2}V) Q^TJQ (\tilde{D}^{-1/2}V)^\top= - \bar J^\top, 
    \\ 
    \bar R &= (\tilde{D}^{-1/2} V) Q^TRQ (\tilde{D}^{-1/2} V)^\top \ge 0, 
     \\
    \bar  B &=  (\tilde{D}^{-1/2}V) Q^TB, \enspace \qquad \bar  P=(\tilde{D}^{-1/2}V) Q^TP .
   % \tag*{\mbox{$\Box$}}
\end{align*}
\end{remark}

If $x$ is a solution of \eqref{eq:linear_PHS}, then $\bar{x} = \tilde{D}^{1/2}Vx$ is a solution of \eqref{eq:phs.ode.transformed}, and since we assume $Q$ to be non-singular, the reverse also holds, i.e., if $\bar{x}$ solves \eqref{eq:phs.ode.transformed}, then $x = V^T\tilde{D}^{-1/2}\bar x$ solves \eqref{eq:linear_PHS}. Moreover, the pencil $\{E,(J-R)Q\}$ is regular if and only if the pencil $\{\bar E, (\bar J - \bar R)\}$ is regular.

\begin{remark}
The index-1 condition for the transformed system \eqref{eq:phs.ode.transformed} reads 
$$
(0 \;\; I_{n-r}) (\bar J- \bar R) 
(0 \;\; I_{n-r})^\top
\qquad \text{ is non-singular}
$$
(where $(0 \;\; I_{n-r}) \in \mathbb{R}^{(n-r)\times n}$).
%\hfill $\Box$ 
\end{remark}

Runge-Kutta schemes are invariant under linear transformation. When analyzing such schemes, it is therefore sometimes convenient to consider the semi-explicit pH-DAE formulation \eqref{eq:phs.ode.transformed} from the very beginning. Also note that in the  applications considered later in this paper, we have no feed-through, i.e., $S=0$ and $N=0$, and a vanishing port matrix $P=0$.

%-------------------------------------------
\section{Dimension-Reducing Decomposition for DAE-Systems}
\label{sec:coupled} 
%-------------------------------------------   

High-dimensional DAEs typically arise in models that couple several lower-dimensional DAEs. In this section, we reverse this approach for a given general DAE $E\dot{x} = f$ by investigating decompositions $f= \mathrm{f}_1 + \mathrm{f}_2$, so that each DAE $E_i \dot{x}= \mathrm{f}_i$, $i=1,2$, effectively has a reduced dimension since components in the equations and in the solution can be eliminated a priori. This amounts to requiring zeros in some components of the $\mathrm{f}_i$. While this is a common technique for ODEs, see, e.g., \cite{Rueth2018}, it is new for DAEs and has been recently introduced in \cite{ECMI_2023} for semi-implicit index-1 DAEs. The challenge is how to deal with the algebraic constraints. It was shown that the algebraic constraints need to be doubled to achieve convergence when applying an operator splitting scheme.  We refer to this approach as dimension-reducing decomposition. In the following, we first discuss the index-1 case and then generalize the concept to private index-2 variables.

%%%%%%%%%%%
\subsection{Coupled index-1 DAEs} \label{subsec:coupled_index1}

Coupled index-1 DAEs often arise in port-Hamiltonian modeling.
Port-Hamiltonian systems are compatible with a coupling process. Two port-Hamil\-tonian systems that are linked via a linear coupling condition on the in- and outputs form a larger new system that is again port-Hamiltonian. As a motivation, consider two linear pH-DAEs without feed-through of the form \eqref{eq:linear_PHS} that are coupled through a condition on the energy balance, i.e.,
\begin{align*}
E_i\dot{x_i}&=(J_i-R_i)Q_i x_i +B_i\,u_i,\qquad
y_i=B_i^\top\,Q_ix_i,\qquad i=1,2,\\
&\quad \begin{pmatrix}
 u_1 \\ u_2 
\end{pmatrix}
+ C
\begin{pmatrix}
    y_1 \\ y_2
\end{pmatrix}=0, \,\,\, \qquad C=-C^\top.
\end{align*}
In this case the coupled closed system can be expressed in a compressed generic pH-form \eqref{eq:linear_PHS} as
\begin{align*}
\begin{pmatrix}
    E_1 & 0 \\ 0 &  E_2 
\end{pmatrix} 
\begin{pmatrix}
    \dot x_1 \\ \dot x_2
\end{pmatrix}
&=\left\{ 
\left[ \begin{pmatrix}
    J_1 & 0\\ 0 & J_2
\end{pmatrix} -
\begin{pmatrix}
B_1 & 0 \\ 0 & B_2
\end{pmatrix}
C
\begin{pmatrix}
B_1 & 0 \\ 0 & B_2
\end{pmatrix}^{\!\!\top}
\right] 
- \begin{pmatrix}
    R_1 & 0 \\ 0 & R_2
\end{pmatrix}
\right\}
\begin{pmatrix}
    Q_1 & 0 \\ 0 & Q_2
\end{pmatrix}
\begin{pmatrix}
    x_1 \\ x_2
\end{pmatrix},
\end{align*}
since the skew-symmetry of the coupling matrix $C$ allows the input terms to be shifted into the off-diagonal blocks of the skew-symmetric system matrix.
If the constraints in the individual DAEs are formulated explicitly---which is the case in many applications---, $E_i$ is a block-diagonal matrix of type 
\begin{align*}
E_i = \begin{pmatrix}
    D_i & 0 \\ 0 & 0
\end{pmatrix}, \qquad D_i \text{ non-singular}.
\end{align*}
With the matching partition $x_i^\top = (x_{i1}^\top,z_{i1}^\top)$ into differential and algebraic variables, the coupled system has then an additional structure, i.e.,
\begin{equation*}
  \begin{aligned}
\begin{pmatrix}
\dot x_{11} \\
0 \\
\dot x_{21} \\
0
\end{pmatrix}
& =
\begin{pmatrix}
    \begin{pmatrix}
        D_1^{-1} & 0 \\ 0 & I
    \end{pmatrix} & 0 \\ 0 &  \begin{pmatrix}
        D_2^{-1} & 0 \\ 0 & I
    \end{pmatrix}  
\end{pmatrix} \\
& \quad \cdot 
\left\{ 
\left[ \begin{pmatrix}
    J_1 & 0\\ 0 & J_2
\end{pmatrix} -
\begin{pmatrix}
B_1 & 0 \\ 0 & B_2
\end{pmatrix}
C
\begin{pmatrix}
B_1 & 0 \\ 0 & B_2
\end{pmatrix}^{\!\!\top}
\right] 
- \begin{pmatrix}
    R_1 & 0 \\ 0 & R_2
\end{pmatrix}
\right\}\\
& \quad \cdot 
 \begin{pmatrix}
    \begin{pmatrix}
        D_1^{-1} & 0 \\ 0 & I
    \end{pmatrix} & 0 \\ 0 &  \begin{pmatrix}
        D_2^{-1} & 0 \\ 0 & I
    \end{pmatrix}  
\end{pmatrix}^{\!\!\top} \cdot
\begin{pmatrix}
    \begin{pmatrix}
        D_1 & 0 \\ 0 & I
    \end{pmatrix} & 0 \\ 0 &  \begin{pmatrix}
        D_2 & 0 \\ 0 & I
    \end{pmatrix}  
\end{pmatrix}^{\!\!\top}
\begin{pmatrix}
 Q_1 \begin{pmatrix} x_{11} \\ z_{11} \end{pmatrix} \\
 Q_2 \begin{pmatrix} x_{21} \\ z_{21} \end{pmatrix}
\end{pmatrix} 
\\
& =: \left\{ 
\begin{pmatrix}
    \tilde J_1 & \tilde J \\ - \tilde J^\top & \tilde J_2
\end{pmatrix}
- \begin{pmatrix}
    \tilde R_1 & 0 \\ 0 & \tilde R_2
\end{pmatrix}
\right\}
\begin{pmatrix}
    \tilde Q_1 & 0 \\ 0 & \tilde Q_2
\end{pmatrix}
\begin{pmatrix}
x_{11} \\ z_{11} \\
x_{21} \\ z_{21} 
\end{pmatrix}.
\end{aligned}
\end{equation*}

This pH-system is a special case of the more general index-1 DAE
\begin{equation} \label{eq:coupled_index1-DAE}
  \begin{aligned}
\begin{pmatrix}
\dot x_{11} \\
0 \\
\dot x_{21} \\
0
\end{pmatrix}
& = 
\begin{pmatrix}
f_1(x_{11}, x_{21}, z_{11}, z_{21})\\ 
g_1(x_{11}, x_{21}, z_{11}, z_{21})\\ 
f_2(x_{11}, x_{21}, z_{11}, z_{21})\\
g_2(x_{11}, x_{21}, z_{11}, z_{21})
\end{pmatrix},
\end{aligned} 
\end{equation}
that was considered in \cite{ECMI_2023} and for which we now summarize the results on dimension-reducing decompositions.\footnote{In the notation of the variables, the first index stands for the associated subsystem. The second index $_1$ is not needed here, but we use it for compatibility with the model equations \eqref{eq:private_index2_system} which are considered later in Section~\ref{subsec:coupled_index2}.}  

We consider the system \eqref{eq:coupled_index1-DAE} on the compact interval $\mathbb{I} = [0,T] \subset \mathbb{R}$. Equipped with initial values $x(0)^\top=x^{(0),\top}=(x_{11}^{(0),\top},x_{21}^{(0),\top})\in \mathbb{R}^{n_x}$ and $z(0)^\top=z^{(0),\top}=(z_{11}^{(0),\top},z_{21}^{(0),\top})\in \mathbb{R}^{n_z}$, we assume that \eqref{eq:coupled_index1-DAE} possesses a unique solution $x^\top = (x_{11}^\top,x_{21}^\top) : \mathbb{I} \rightarrow  \mathbb{R}^{n_x}$ and $z^\top = (z_{11}^\top,z_{21}^\top) : \mathbb{I} \rightarrow  \mathbb{R}^{n_z}$. 
The functions $f_i$, $g_i$, $i\in\{1,2\}$ are supposed to be sufficiently differentiable in the neighborhood of the solution (note that this is trivially fulfilled in the linear case), and $\partial (g_1,\, g_2)/\partial (z_{11},z_{21}), \partial g_1/\partial z_{11}$ and $\partial g_2/\partial z_{21}$ are supposed to be non-singular in a neighborhood of the solution. 

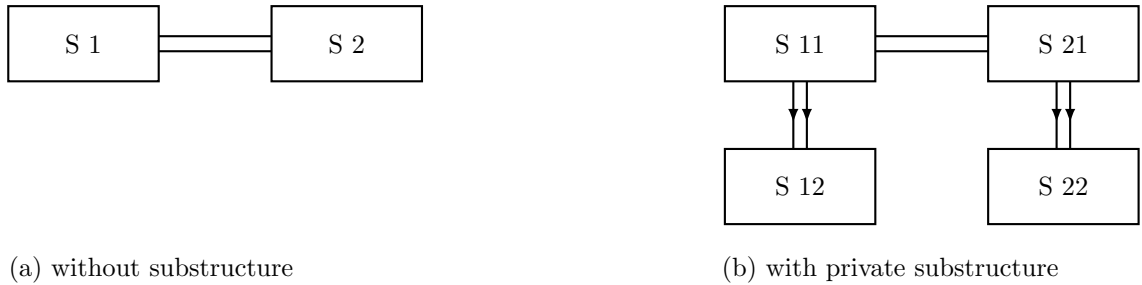
\begin{figure}[tb]
\begin{center}
\begin{tikzpicture}[scale=1.]
\draw[thick] (0,0)--(2,0)--(2,1) -- (0,1) -- cycle;
\draw[xshift=1cm,yshift=0.5cm] node {S 1}; 
\draw[thick,xshift=3.5cm] (0,0)--(2,0)--(2,1) -- (0,1) -- cycle;
\draw[xshift=4.5cm,yshift=0.5cm] node {S 2}; 
\draw[thick] (2,0.4)--(3.5,0.4);
\draw[thick] (2,0.6)--(3.5,0.6);
\draw[xshift=1.9cm,yshift=-2.5cm] node{(a) without substructure};
\end{tikzpicture}\hfill
\begin{tikzpicture}[scale=1.]
\draw[thick] (0,0)--(2,0)--(2,1) -- (0,1) -- cycle;
\draw[xshift=1cm,yshift=0.5cm] node {S 11}; 
\draw[thick,-latex] (0.91,0)--(0.91,-0.55);
\draw[thick] (0.91,-0.45)--(0.91,-0.9);
\draw[thick,-latex] (1.09,0)--(1.09,-0.55);
\draw[thick] (1.09,-0.45)--(1.09,-0.9);
\draw[thick,yshift=-1.9cm] (0,0)--(2,0)--(2,1) -- (0,1) -- cycle;
\draw[xshift=1cm,yshift=-1.4cm] node {S 12}; 

\draw[thick,xshift=3.5cm] (0,0)--(2,0)--(2,1) -- (0,1) -- cycle;
\draw[xshift=4.5cm,yshift=0.5cm] node {S 21}; 
\draw[thick] (2,0.4)--(3.5,0.4);
\draw[thick] (2,0.6)--(3.5,0.6);
\draw[xshift=2.2cm,yshift=-2.5cm] node{(b) with private substructure};
\draw[thick,xshift=3.5cm,-latex] (0.91,0)--(0.91,-0.55);
\draw[thick,xshift=3.5cm,-latex] (1.09,0)--(1.09,-0.55);
\draw[thick,xshift=3.5cm] (0.91,-0.45)--(0.91,-0.9);
\draw[thick,xshift=3.5cm] (1.09,-0.45)--(1.09,-0.9);
\draw[thick,,xshift=3.5cm,yshift=-1.9cm] (0,0)--(2,0)--(2,1) -- (0,1) -- cycle;
\draw[xshift=4.5cm,yshift=-1.4cm] node {S 22}; 
\end{tikzpicture}
\end{center}
 \caption{Coupled system consisting of two subsystems (S).}\label{fig:coupled_index_1_networks}
\end{figure}

According to \cite{ECMI_2023}, doubling the algebraic conditions is a way to obtain a decomposition for which we can use operator splitting schemes for DAEs, i.e., we take the decomposition 
\begin{equation}
    \label{eq:DAE_dimred_decomposition}
    \begin{pmatrix}
        \dot{x}_{11} \\ 0 \\ \dot{x}_{21} \\ 0 
    \end{pmatrix} 
    = \begin{pmatrix}
        f_1 \\ 2g_1 \\ f_2 \\ 2g_2 
    \end{pmatrix} 
    = \underbrace{%
       \begin{pmatrix}
            f_1 \\ g_1  \\ 0 \\ g_2 
        \end{pmatrix}%
                  }_{=:\mathrm{f}_1
                  }
        + \underbrace{%
         \begin{pmatrix}
            0 \\ g_1 \\ f_2 \\ g_2 
         \end{pmatrix}%
                }_{=:\mathrm{f}_2
                }.
\end{equation}
In the two resulting DAEs $E\dot x = \mathrm{f}_i, i=1,2$, the third and first block then become trivial, reading $\dot x_{21} = 0$ and $\dot x_{11} = 0$, respectively, in this manner reducing the effective dimensions.  
We thus have the coupled subsystems
\begin{alignat}{2}  
   &\text{(subsystem 1)} & & \text{(subsystem 2)} \nonumber \\
    \dot{x}_{11} & = f_1 (x_{11}, x_{21}, z_{11}, z_{21}),  
            \qquad\qquad & \dot{x}_{21} & = f_2 (x_{11}, x_{21}, z_{11}, z_{21}), \label{eq:DAE_dimred_subsystems}\\
      0 & = g_1  (x_{11}, x_{21}, z_{11}, z_{21}), 
            \qquad &   0 & =g_1 (x_{11}, x_{21}, z_{11}, z_{21}), \nonumber \\
      0 & = g_2 (x_{11}, x_{21}, z_{11}, z_{21}),   \;   
            \qquad & 0 & =g_2 (x_{11}, x_{21}, z_{11}, z_{21}), \nonumber
\end{alignat}
where $x_{21}$ is the coupling variable for subsystem 1, and $x_{11}$ is the coupling variable for subsystem 2. This is schematically depicted in Fig.~\ref{fig:coupled_index_1_networks}~(a).

Due to the index-1 and regularity assumptions, the implicit function theorem implies that we can reduce the two DAEs in \eqref{eq:DAE_dimred_subsystems} to the two ODEs 
\begin{equation}\label{eq:ODE_dimred_subsystems}
\begin{array}{rcl}
\dot{x}_{11} &=\,& f_1(x_{11}, x_{21},
		       \varphi(x_{11},x_{21})) \\
         &=:& h_1(x_{11},x_{21}),
\end{array}\;   \qquad
\begin{array}{rcl}
 \dot{x}_{21} &=\,& f_2(x_{11}, x_{21},
 			   \varphi(x_{11},x_{21})) \\
       &=:& h_2(x_{11},x_{21}),
\end{array}
\end{equation}
where $\varphi: \mathbb{R}^{ n_{x}} \rightarrow \mathbb{R}^{ n_{z}}$ is the continuously differentiable function with $(z_{11},z_{21}) = \varphi(x_{11},x_{21})$. The two systems in \eqref{eq:ODE_dimred_subsystems} arise from the dimension-reducing decomposition of the ODE
\begin{equation} \label{eq:ODE_dimred_decomposition}
\begin{pmatrix}
\dot x_{11} \\
\dot x_{21}
\end{pmatrix} 
= \begin{pmatrix}
 h_1(x_{11},x_{21}) \\
 h_2(x_{11},x_{21}) 
\end{pmatrix}
=
\underbrace{\begin{pmatrix}
h_1(x_{11},x_{21}) \\
0 
\end{pmatrix}
}_{=: \mathrm{h}_1}
+
\underbrace{
\begin{pmatrix}
0 \\
h_2(x_{11},x_{21}) 
\end{pmatrix}
}_{=: \mathrm{h}_2}.
\end{equation}

From an analytical viewpoint, any DAE operator splitting scheme applied to the decomposition \eqref{eq:DAE_dimred_decomposition} and thus working with the two DAEs \eqref{eq:DAE_dimred_subsystems} is equivalent to an ODE splitting scheme applied to the decomposition \eqref{eq:ODE_dimred_decomposition}, working with the two ODEs from \eqref{eq:ODE_dimred_subsystems}. This gives the following result. 

\begin{proposition} \label{prop:dimred_order}
Let all the assumptions formulated for system \eqref{eq:coupled_index1-DAE} be fulfilled. Assume that a given splitting scheme has convergence order $p$ for the dimension-reducing ODE decomposition \eqref{eq:ODE_dimred_decomposition}. Then, the same splitting scheme applied to the dimen\-sion-reducing DAE decomposition \eqref{eq:DAE_dimred_decomposition} also has order $p$. 
\end{proposition}
Consequently, the Lie-Trotter scheme has convergence order $p=1$, and the Strang splitting scheme has order $p=2$.

In practice, we usually cannot solve the subsystems in \eqref{eq:DAE_dimred_subsystems} exactly when working with an operator splitting scheme. We rather resort to a numerical, time-discrete integrator. And appropriate numerical schemes for the DAEs \eqref{eq:DAE_dimred_subsystems} may not necessarily rely on using a numerical scheme for the ODEs \eqref{eq:ODE_dimred_subsystems} with explicitly given function $\varphi(x_{11},x_{21})$. This is why from a numerical perspective the following result is worth mentioning even though it might be regarded a straight-forward consequence of Proposition~\ref{prop:dimred_order}.

\begin{corollary}\label{cor:dimred_order_discrete}
Under the assumptions of Proposition~\ref{prop:dimred_order}, if within the given splitting scheme of convergence order $p$ we apply numerical DAE time integrators which have minimal order $p$ to the subsystems \eqref{eq:DAE_dimred_subsystems}, the discrete numerical approximations are also of order $p$. 
\end{corollary}

\begin{remark}
The numerical integrator for the DAEs must preserve the convergence order of the splitting scheme in the approximation of the differential and algebraic variables. Combining for example a fourth-order Triple Jump splitting scheme with a 2-stage Gauss-Runge-Kutta method might lead to a loss of convergence in the algebraic variables, whereas the use of a 3-stage Lobatto-IIIC method ensures the fourth-order convergence also for the algebraic variables in the index-1 case. For details on the convergence orders of established splitting schemes and numerical integrators, we refer to Appendix~\ref{appendix:splitting-approximation}. 
%\hfill $\Box$
\end{remark}

\begin{remark}
 Proposition~\ref{prop:dimred_order} and Corollary~\ref{cor:dimred_order_discrete} can be generalized to systems which have multiple subsystems similar to \eqref{eq:coupled_index1-DAE}, where all subsystems and the overall system are of index~1. 
 %\hfill $\Box$
\end{remark}

%%%%%%%%%%%
\subsection{Coupled DAEs with private index-2 variables}\label{subsec:coupled_index2}

Proceeding from \eqref{eq:coupled_index1-DAE}, we extend the model system by incorporating substructures with private index-2 variables, see Fig.~\ref{fig:coupled_index_1_networks}~(b). We thus allow for index-2 variables which are not seen from the other subsystem. This can be modeled, for instance, by,
\begin{subequations}\label{eq:private_index2_system}
\begin{align}
\dot{x}_{11} &= f_1(x_{11}, \myempty, x_{21}, \myempty,
		       z_{11}, \myempty, z_{21}, \myempty),  \label{eq:sub1a}\\ 
0 		&= g_1(x_{11}, \myempty, {x_{21}}, \myempty,
			   z_{11}, \myempty, {z_{21}}, \myempty), \label{eq:sub1b}\\ 
\dot{x}_{12} &= h_1(x_{11}, x_{12}, \myempty, \myempty, 					
		 	   z_{11}, z_{12}, \myempty, \myempty) , \label{eq:sub1c}\\
0 		&= k_1(x_{11}, x_{12}, \myempty, \myempty, 		
			    z_{11},z_{12}, \myempty, \myempty), \label{eq:sub1d} 
\\[1ex]
\dot{x}_{21} &= f_2(x_{11}, \myempty, x_{21}, \myempty, 
 			   z_{11}, \myempty, z_{21}, \myempty), \label{eq:sub2a}\\
0 		&= g_2(x_{11}, \myempty, x_{21}, \myempty,
			   z_{11}, \myempty, z_{21}, \myempty), \label{eq:sub2b} \\
\dot{x}_{22} &= h_2(\myempty, \myempty, x_{21}, x_{22}, 
               \myempty, \myempty, z_{21}, z_{22}), \label{eq:sub2c}\\
0 		&= k_2(\myempty, \myempty, x_{21}, x_{22}, 
			   \myempty, \myempty, z_{21}, z_{22}). \label{eq:sub2d} 
\end{align}
\end{subequations}
The first subsystem \eqref{eq:sub1a}--\eqref{eq:sub1d} is composed of two parts: The part \eqref{eq:sub1a}--\eqref{eq:sub1b} manages the coupling, and we assume that it is an index-1 system for $x_{11}$, $z_{11}$ when all other variables are given. The other part \eqref{eq:sub1c}--\eqref{eq:sub1d} describes a private substructure, which we assume to be of index~2 for $x_{12}$, $z_{12}$, where all other variables are given. Analogous index conditions are assumed for the second subsystem \eqref{eq:sub2a}--\eqref{eq:sub2d}. 
To fix the dimensions, let
 \begin{alignat*}{1}
   &f_i: \mathbb{R}l^{m}\to\mathbb{R}^{n_{xi1}}, \quad g_i: \mathbb{R}^{m}\to\mathbb{R}^{n_{zi1}}, \\
   &h_i: \mathbb{R}^{n_i}\to\mathbb{R}^{n_{xi2}}, \quad
  %   \\
   k_i: \mathbb{R}^{n_i}\to\mathbb{R}^{n_{zi2}},  
\end{alignat*}
using $m = n_{x11} + n_{x21} + n_{z11} + n_{z21}$ and $n_i=n_{xi1}+n_{xi2}+n_{zi1}+n_{zi2}$, $i\in\{1,2\}$.
We equip the model equations \eqref{eq:private_index2_system} with consistent initial values: 
\begin{equation*}
\label{cdt_initial}
    x(0) = x^{(0)} =\begin{pmatrix} x_{11}^{(0)}
    \\x_{12}^{(0)} \\ x_{21}^{(0)}\\ x_{22}^{(0)}  \end{pmatrix}
    \in \mathbb{R}^{n_x}, 
    \qquad z(0) = z^{(0)} = \begin{pmatrix} 
    z_{11}^{(0)} \\ z_{12}^{(0)} \\ z_{21}^{(0)} \\ z_{22}^{(0)}  \end{pmatrix} \in \mathbb{R}^{n_z}.
\end{equation*}
The resulting initial value problem is assumed to have a unique solution $x^\top=(x_{11}^\top, x_{12}^\top, x_{21}^\top, x_{22}^\top) : \mathbb{I} \rightarrow \mathbb{R}^{n_x}$ and $z^\top=(z_{11}^\top,z_{12}^\top,z_{21}^\top,z_{22}^\top) : \mathbb{I} \rightarrow \mathbb{R}^{n_z}$ on $\mathbb{I} = [0,T]$, where $n_\alpha = n_{\alpha 11}+n_{\alpha 12}+n_{\alpha 21}+n_{\alpha 22}$ for $\alpha\in \{x,z\}$. Moreover, all functions $f_i$, $g_i$, $h_i$, $k_i$, $i \in \{1,2\}$ are supposed to be sufficiently differentiable in a neighborhood of the solution. 

The private index-$2$ substructures (\ref{eq:sub1c}--\ref{eq:sub1d})  and (\ref{eq:sub2c}--\ref{eq:sub2d}) are incorporated via the index-1 part, i.e., (\ref{eq:sub1a}--\ref{eq:sub1b},\ref{eq:sub2a}--\ref{eq:sub2b}). This is a one-way coupling (one-sided dependence), because the index-1 part (describing the overall coupling) does not depend on the solutions $x_{12}, z_{12}, x_{22}, z_{22}$ of the index-2 substructures. 

We now decompose the system \eqref{eq:private_index2_system} into two subproblems according to
\begin{equation}
    \label{DAE_decomposition-private}
    \begin{pmatrix}
        \dot{x}_{11} \\ 
        0 \\
        \dot{x}_{12} \\
        0 \\
        \dot{x}_{21} \\
        0 \\ 
        \dot{x}_{22} \\ 0 
    \end{pmatrix} \;
    = \;
    \begin{pmatrix}
        f_1 \\ 
        2g_1 \\ 
        h_1 \\
        k_1 \\
        f_2 \\ 
        2g_2 \\
        h_2 \\
        k_2
    \end{pmatrix}
    \; = \underbrace{%
       \begin{pmatrix}
        f_1 \\ 
        g_1 \\ 
        h_1 \\
        k_1 \\
        0 \\ 
        g_2 \\
        0 \\
        0
    \end{pmatrix}
                  }_{=:\mathrm{f}_1, \text{ subproblem 1}} 
        + \underbrace{%
         \begin{pmatrix}
        0 \\ 
        g_1 \\ 
        0 \\
        0 \\
        f_2 \\ 
        g_2 \\
        h_2 \\
        k_2
    \end{pmatrix}%
                }_{=:\mathrm{f}_2, \text{ subproblem 2}},
\end{equation} 
where we double the algebraic constraints $g_i$ corresponding to the coupling part, but not the algebraic constraints $k_i$ on the private substructures. To get a standard splitting, all variables need to be considered in each subproblem. For this purpose, the occurring terms $0=0$ need to be interpreted, such that the splitting procedure is well-defined. In fact, for subsystem 1 the information of $x_{22}, z_{22}$ is not required and there is no update of these variables. Therefore, analogously to equation $\dot{x}_{22}=0$ (given by \eqref{DAE_decomposition-private}), we also use  $\dot{z}_{22}=0$ in the operator splitting to describe that these variables remain unchanged when the first subsystem is solved. The same applies to subsystem 2. Thus our decomposition deals with the subproblems $E_i \dot x=\mathrm{f}_i$, $i\in\{1,2\}$ given by
\begin{alignat}{2}E_1\dot x=
     \begin{pmatrix}
        \dot{x}_{11} \\ 
        0 \\
        \dot{x}_{12} \\
        0 \\
        \dot{x}_{21} \\
        0 \\ 
        \dot{x}_{22} \\ 
        \dot{z}_{22}
    \end{pmatrix} \;
    = \;
       \begin{pmatrix}
        f_1 \\ 
        g_1 \\ 
        h_1 \\
        k_1 \\
        0 \\ 
        g_2 \\
        0 \\
        0
    \end{pmatrix}=\mathrm{f}_1,
   \quad\qquad E_2 \dot x=
     \begin{pmatrix}
        \dot{x}_{11} \\ 
        0 \\
        \dot{x}_{12} \\
        \dot{z}_{12} \\
        \dot{x}_{21} \\
        0 \\ 
        \dot{x}_{22} \\ 0 
    \end{pmatrix} \;
    = \;
         \begin{pmatrix}
        0 \\ 
        g_1 \\ 
        0 \\
        0 \\
        f_2 \\ 
        g_2 \\
        h_2 \\
        k_2
    \end{pmatrix}=\mathrm{f}_2
    . \label{DAE_decomposition-private2}
\end{alignat}
The effective dimension of the subsystems can be obviously reduced by solving only the following equations, while the remaining variables are kept fixed in the respective splitting step: 
\begin{subequations}\label{eq:DAE_dimred}
\begin{align}\nonumber
&\text{\hspace*{-1.cm}(subsystem~1)} &&\text{\hspace*{-1.cm}(subsystem~2)}\\ \label{DAE_dimred_index1}
\text{coupling:\hspace*{1.25cm}}\begin{pmatrix}
     \dot{x}_{11} \\ 0 \\   0 
\end{pmatrix} 
& = 
\begin{pmatrix}
f_1 \\ g_1   \\ g_2 
\end{pmatrix}, &
\begin{pmatrix}
         \dot{x}_{21} \\ 0 \\ 0 
\end{pmatrix} 
&= 
\begin{pmatrix}
f_2 \\ g_1  \\ g_2 
\end{pmatrix}, \\
    \label{DAE_dimred_index2} 
\text{substructure:\hspace*{0.75cm}}    
\begin{pmatrix}
        \dot{x}_{12} \\
        0  
    \end{pmatrix} 
    &= 
    \begin{pmatrix}
        h_1 \\
        k_1 
    \end{pmatrix}, &
    \begin{pmatrix} 
        \dot{x}_{22} \\ 
        0 
    \end{pmatrix} 
    &= 
    \begin{pmatrix}
        h_2 \\
        k_2
    \end{pmatrix}.
\end{align}
\end{subequations}

The following proposition is an immediate consequence of the above observations.

\begin{proposition} \label{prop:dimred_order2}
Let all the assumptions formulated for system \eqref{eq:private_index2_system} be fulfilled. Assume that a given splitting scheme applied to the two index-1 DAEs in \eqref{DAE_dimred_index1} has convergence order $p$. Then, the same splitting scheme applied to the dimension-reducing DAE decomposition \eqref{DAE_decomposition-private2}, working on \eqref{eq:DAE_dimred}, has also order~$p$ for the index-0 and index-1 coupling variables $x_{11}, x_{21}, z_{11}, z_{21}$. Moreover, it has convergence order $p$ for all variables, if the substructures \eqref{DAE_dimred_index2} do not depend on the {index-1} coupling variables $z_{11}$, $z_{21}$ or if the substructure of only one subsystem depends on the index-1 coupling variable and this subsystem is solved in the final stage of the splitting scheme.
\end{proposition}

Note that due to Proposition~\ref{prop:dimred_order}, the splitting scheme for the index-1 DAEs in \eqref{DAE_dimred_index1} has order $p$ if the same splitting scheme has order $p$ for the associated dimension-reducing ODE decomposition, cf.\ \eqref{eq:ODE_dimred_subsystems} and \eqref{eq:ODE_dimred_decomposition}.

\begin{corollary}\label{cor:dimred_order_discrete2}
Under the assumptions of Proposition~\ref{prop:dimred_order2}, if within the given splitting scheme of convergence order $p$ we apply numerical DAE time integrators of minimal order $p$ to \eqref{DAE_dimred_index1}, the discrete numerical approximations of the index-0 and index-1 coupling variables $x_{11},x_{21},z_{11},z_{21}$ are also of order $p$. If numerical DAE time integrators of minimal order $p$ (for index-2 DAEs) are applied to \eqref{eq:DAE_dimred} where the dependencies of the substructures \eqref{DAE_dimred_index2} on the index-1 coupling variables are as specified in Proposition~\ref{prop:dimred_order2}, then the numerical approximations of all variables are of order $p$.
\end{corollary}

In computational practice, it is preferable to exploit the one-sided dependencies between substructure and coupling part in the subproblems. Instead of solving the entire system~\eqref{eq:DAE_dimred} numerically within the splitting scheme and possibly losing convergence order in the substructure variables, it is advisable to integrate only the coupling part \eqref{DAE_dimred_index1} numerically in the splitting. After having determined $x_{11},z_{11},x_{21},z_{21}$ on the full time interval, the private substructures \eqref{DAE_dimred_index2} can then be computed separately with any appropriate numerical integrator. This decoupling strategy allows the specific choice of different integrators with regard to the differentiation index, accounts for the order in all variables of \eqref{eq:private_index2_system} and reduces the computational effort.

Proposition~\ref{prop:dimred_order2} and Corollary~\ref{cor:dimred_order_discrete2} generalize to systems which have multiple subsystems of similar structure with private index-2 variables.

%---------------------------------------------------------------------
\section{Energy-Associated $J$-$R$ Decomposition for Index-1 pH-DAEs}
\label{sec:energy-based-splitting}
%---------------------------------------------------------------------

The focus of this section is on the linear pH-DAE \eqref{eq:phs.ode.transformed} without feed-through and port matrix, i.e.,
\begin{align}
   \label{eq:PHS-General}
    E\dot{x} & = (J- R) x + Bu(t), \qquad x(0) = x_0, \\
    y& = B^\top  x \nonumber
\end{align}
with $J=-J^\top$, $R^\top = R\ge 0$, and $E=E^\top \ge 0$ as well as sufficiently smooth input function $u$.  We assume the existence of a classical solution $x\in\mathcal{C}^1(\mathbb{I},\mathbb{R}^d)$. For the underlying Hamiltonian $\mathcal{H}$, the solution then satisfies the dissipation inequality, i.e., for $t\in\mathbb{I}$ it holds
\begin{align}\label{eq:PHS-General-H}
\frac{d}{dt} \mathcal{H}(x(t))  \leq y(t)^\top \, u(t), \qquad  \mathcal{H}(x)=\frac{1}{2}x^\top E x.
\end{align}
The idea of the energy-associated $J$-$R$ decomposition is the separation of the energy-conserving part, $\mathrm{f_1}=Jx$, and the dissipative part with sources, $\mathrm{f_2}=-Rx+Bu$.
In case of $E=I$ identity, \eqref{eq:PHS-General} becomes an explicit pH-ODE, for which the first-order Lie-Trotter splitting and the second-order Strang splitting preserve the dissipation inequality, cf.\ \cite{Frommer_2023aa}. Higher order splitting schemes ($p\geq 3$) destroy the dissipativity due to the occurrence of negative step sizes \cite{suzuki1991}. A way out of this dilemma can be provided by, e.g., force gradient-based splitting methods, cf.\ \cite{moench2023}. Concerning the numerical flux approximation for schemes of order $p\leq 2$, the implicit midpoint rule ensures the energy conservation on a discrete level. It can be interpreted as Cayley transform, and if the Cayley transform is iteratively computed via Arnoldi Krylov-subspace approximations in a matrix function sense, all iterates preserve energy, too, cf.\ \cite{Frommer_2023aa}.

In this section, we extend the $J$-$R$ decomposition strategy and associated structure-preserving splitting to implicit pH-ODEs and index-1 pH-DAEs. Thereby, the handling of the algebraic constraints deserves special attention. Moreover, we introduce a novel concept of generalized Cayley transforms to uphold energy-conservation.

%%%%%%%%%%%%%%
\subsection{Handling of constraints and structure-preserving splitting}

The matrix pencil $\{E, (J-R)\}$ is assumed to be regular, and the index of the pH-DAE \eqref{eq:PHS-General} is assumed to be less or equal to one.
Let $K_E$ be a projector onto $\ker(E)$ and let $P_E=I-K_E$. Note that $E P_E =E$. Then, these projectors induce a partitioning $x = x_d + x_a$ via $x_d = P_E x$ and $x_a = K_E x$. To ensure a consistent interplay between algebraic constraints and $J$-$R$ decomposition we assume:
\begin{assumption}\label{ass:restriction-index1}
   For the pH-DAE \eqref{eq:PHS-General},  one of the following restrictions holds:
   \begin{enumerate}
       \item[(a)] $K_E^\top R = 0$ and $K_E^\top B = 0$ (i.e., $R$ and $B$ are restricted to the differential part),  and the matrix pencil $\{E,\, J\}$ is non-singular; or
       \item[(b)] $K_E^\top J = 0$ (i.e., $J$ is restricted to the differential part), and the matrix pencil $\{E,\, R\}$ is non-singular. 
   \end{enumerate}
\end{assumption}
If $\ker(E)$ is nontrivial, we have that in case (a) the matrix pencil $\{E,R\}$ is singular, whereas in case (b) $\{E,J\}$ is singular. We regularize the respective matrix pencil by introducing $E_r=E + K_E^\top K_E$ (such that $E_r=E_r^\top >0$) to obtain a consistent splitting scheme.
The $J$-$R$ decomposition of \eqref{eq:PHS-General} then prescribes the following two subproblems according to the two cases in Assumption~\ref{ass:restriction-index1}
\begin{subequations}\label{eq:J-R-subs}
\begin{align}\label{eq:J-R-subs_a}
\text{(a)} \quad E\dot x&=Jx && E_r \dot x=-Rx+Bu\\
\text{(b)} \quad E\dot x&=-Rx+Bu && E_r \dot x=Jx.\label{eq:J-R-subs_b}
\end{align}
\end{subequations}
In particular, in the semi-explicit formulation of \eqref{eq:phs.ode.transformed}, the subsystems are given by
\begin{align*}
&\text{(a)} \,\,\,   \begin{pmatrix} I_1 & \\
                    & 0
   \end{pmatrix}
   \begin{pmatrix}
   \dot{x}_1 \\
   \dot{x}_2
   \end{pmatrix}
   =  \begin{pmatrix} 
                J_{11} & J_{12}
                \\
                J_{21} & J_{22}
            \end{pmatrix}
    \begin{pmatrix}
       x_1 \\ x_2
   \end{pmatrix},
  &&
   \begin{array}{c}
   \dot{x}_1 \\
   \dot{x}_2
   \end{array}
   \begin{array}{l}
    = - R_{11} x_1 + B_1 u\\
    =0
     \end{array} \\
& \text{(b)} \,\,\,   \begin{pmatrix} I_1 & \\
                    & 0
   \end{pmatrix}
   \begin{pmatrix}
   \dot{x}_1 \\
   \dot{x}_2
   \end{pmatrix}
   =  -\begin{pmatrix} 
                R_{11} & R_{12}
                \\
                R_{21} & R_{22}
            \end{pmatrix}
    \begin{pmatrix}
       x_1 \\ x_2
   \end{pmatrix}+
   \begin{pmatrix} 
                B_{1} 
                \\
                B_{2} 
            \end{pmatrix}
  u, \quad
   &&
   \begin{array}{c}
   \dot{x}_1 \\
   \dot{x}_2
   \end{array}
   \begin{array}{l}
    = J_{11} x_1 \\
    =0.
     \end{array} 
\end{align*}
In both cases, the $J$-$R$ decomposition results in an index-1 DAE and an ODE subproblem. Expressing the algebraic variable $x_2$ in terms of the differential variable $x_1$, we can rewrite the subproblems as inherent ODEs for $x_1$, i.e.,
\begin{align*}
&\text{(a)} \,\,\, \dot x_1= (J_{11} + J_{12} J_{22}^{-1} J_{21}) x_1, &&\dot x_1=-R_{11} x_1 + B_1u\\
&\text{(b)} \,\,\, \dot x_1= -(R_{11} + R_{12} R_{22}^{-1} R_{21}) x_1 +(B_1+R_{22}^{-1}B_2)u,  &&\dot x_1=J_{11} x_1. 
\end{align*}
Consequently, the $J$-$R$ decomposition of \eqref{eq:PHS-General} is equivalent to a $J$-$R$ decomposition of a pH-ODE for the differential variable $x_1$---in the form of 
\begin{align*}
\dot x_1 = \underbrace{J_1x_1}_{=:\mathrm{f}_1} + \underbrace{(-R_1x_1 +\tilde B_1u)}_{=:\mathrm{f}_2}, \qquad J_1=-J_1^\top, \quad R_1=R_1^\top \geq 0.
\end{align*} 

\begin{proposition}\label{prop:JR-splitting}
Let the port-Hamiltonian system \eqref{eq:PHS-General}  on $\mathbb{I}=[0,T]$ have at most index~1 and let Assumption~\ref{ass:restriction-index1} be fulfilled.  Given a splitting scheme of order $p$ for (explicit) ODEs, then this splitting scheme applied to the subsystems \eqref{eq:J-R-subs} of the $J$-$R$-decomposition has convergence order $p$. Moreover, the Lie-Trotter splitting and the Strang splitting preserve the dissipation inequality \eqref{eq:PHS-General-H}.
\end{proposition}

\begin{remark} In Proposition~\ref{prop:JR-splitting}, the implicit ODE case (i.e., $\ker(E)=\{0\}$) is also covered. 
%\hfill $\Box$
\end{remark}

\begin{remark} For the convergence of the splitting schemes, an appropriate handling of the algebraic constraints is crucial. Whereas in the dimension-reducing decomposition in Section~\ref{sec:coupled} the constraints are doubled and incorporated in all subproblems, the $J$-$R$ decomposition requires a clear assignment of the constraints to one subproblem, which goes hand in hand with the regularization of the other subproblem. Note that a decomposition of the algebraic constraints is not possible as it destroys the convergence of the splitting scheme. As an example take $E = 0$, $J$ and $R$ non-singular, and input $Bu(t) = b$ constant. Then, $0 = (J-R)x +b$ is a linear system for $x$.  Its solution $x = (J-R)^{-1}b$ can obviously not be approximated by a splitting scheme acting on the subproblems  $Jx=0$ and $Rx+b=0$. See also Example~\ref{counter_example} in Section~\ref{sec:numerics-JR} for a less trivial setting. 
%\hfill $\Box$
\end{remark}

Let Assumption~\ref{ass:restriction-index1} not be fulfilled. To enable a $J$-$R$ decomposition-based splitting consistent with the algebraic constraints also in this case, we make use of perturbation theory and adopt the idea of $\epsilon$-embedding methods \cite{IVP_DAE95}, cf.\ \cite{ECMI_2023}. For \eqref{eq:PHS-General}, we consider the regularized system with regularization parameter $\epsilon$,
\begin{align*}
    E_\epsilon \dot{x_\epsilon} & = (J- R) x_\epsilon + Bu(t), \qquad x_\epsilon(0) = x_0, \qquad y_\epsilon = B^\top  x_\epsilon \\
    E_\epsilon &= E+\epsilon \, K_E^\top K_E, \qquad E_\epsilon=E_\epsilon^\top >0,  \qquad 0<\epsilon\ll 1.
\end{align*}
This is an implicit pH-ODE with Hamiltonian $\mathcal{H}_\epsilon(x)=\tfrac{1}{2}x^T E_\epsilon x$.
Due to regularity it possess a unique solution $x_\epsilon\in \mathcal{C}^1(\mathbb{I},\mathbb{R}^d)$ which converges uniformly to the solution $x$ of the original index-1 pH-DAE as $\epsilon \rightarrow 0$. As for the dissipation inequality, we have
\begin{align*}
\mathcal{H}&(x_\epsilon(t+h))-\mathcal{H}(x_\epsilon(t)) \\
&= 
(\mathcal{H}(x_\epsilon(t+h))-\mathcal{H}_\epsilon(x_\epsilon(t+h)))
+(\mathcal{H}_\epsilon(x_\epsilon(t+h))-\mathcal{H}_\epsilon(x_\epsilon(t)))
+(\mathcal{H}_\epsilon(x_\epsilon(t))-\mathcal{H}(x_\epsilon(t)))\\
&\leq \int_t^{t+h}(y_\epsilon(s))^\top u(s) \,\mathrm{d}s + \frac{\epsilon}{2}\big|  \|K_E \, x_\epsilon(t+h)\|_2^2 -  \| K_E \, x_\epsilon(t)\|_2^2\big|.
\end{align*}
For the first term we have convergence $\int_t^{t+h} y_\epsilon^T u \, \mathrm{d} s \rightarrow \int_t^{t+h} y^T u \, \mathrm{d} s$ as $\epsilon\rightarrow 0$, since the integrand is bounded on $[t,t+h]$ and $y_\epsilon=B^\top x_\epsilon \rightarrow B^\top x= y$. The second term is of order $\mathcal{O}(\epsilon h)$.  Thus, the dissipation inequality \eqref{eq:PHS-General-H} is satisfied in the limit $\epsilon=0$. Applying a splitting scheme of order $p$ onto a $J$-$R$ decomposition of the regularized pH-ODE, we can expect a respective $p$-th order approximation for the original index-1 pH-DAE as $\epsilon \rightarrow 0$.

%%%%%%%%%%%%%%%%%%%%%%%
\subsection{Generalized Cayley transforms for time discretization}
\label{sec:JR-splitting-time-discrete}

In the structure-preserving splitting schemes (Lie-Trotter and Strang) on top of the $J$-$R$ decomposition, the subproblems \eqref{eq:J-R-subs} have to be solved. The numerical flux approximation should thereby preserve the properties of the respective subproblems and the convergence order of the splitting scheme. Our special interest here is the conservation of energy in $E\dot x =Jx$, $E=E^\top\geq 0$, $J=-J^\top$ with regular matrix pencil $\{E,J\}$.

For simplicity of notation, we consider an equidistant time grid for discretization: $0=t_0 < t_1 < \dotsc < t_N=T$ with $h=t_{n+1}-t_{n}$ and $x(t_n)\approx x_n$. The implicit midpoint rule (1-stage Gauss-Runge Kutta method) applied to the energy-conserving subsystem yields 
\begin{align} \nonumber
 & \quad \tfrac{1}{h} E \left( x_{n+1} -x_n \right) = \tfrac{1}{2} J \left( x_{n+1} + x_{n} \right) \\
  \Leftrightarrow & \quad x_{n+1} = C(E,\tfrac{h}{2}J) x_n
 \label{eq:PHS-DAE-midpointrule}
\end{align}
where
\begin{align*}
   C(E,A) := (E-A)^{-1}(E+A)
\end{align*}
denotes the generalized Cayley transform for $E,A\in \mathbb{R}^{n\times n}$, $E-A$ non-singular. 
In Appendix~\ref{app:Cayley}, we prove the following result.

\begin{lemma}\label{lem:Cayleytrafonorm}
Assume that $E^\top =E \ge 0$, $A=-A^\top$, $E,A \in \mathbb{R}^{n \times n}$, 
and that $\{E,A\}$ is a regular matrix pencil, i.e., the nullspaces satisfy $\nullsp(E) \cap \nullsp(A) = \{0\}$. 
Then, the generalized Cayley transform $C(E,A)$ satisfies 
\begin{equation} \label{eq:isometry}
\|C(E,A)x\|_E  = \|x\|_E \quad \text{ for } x \in \mathbb{R}^n 
\end{equation}
with $\|x\|_E = \langle x, x \rangle_E^{1/2}$ and $\langle x, y \rangle_E = y^\top Ex$. 
\end{lemma}

\begin{remark}\label{remark:semi-explicit-trafo}
The generalized Cayley transforms describing the implicit midpoint rule for the skew-symmetric part associated with the system~\eqref{def:pHS-DAE} and the linearly transformed system~\eqref{eq:phs.ode.transformed} in Remark~\ref{rem:QE-ID} are related via
\begin{align*}
C(E,JQ) = C(Q^\top E,Q^\top JQ) = (V^\top \tilde{D}^{-1/2}) C(\bar E,\tfrac{h}{2}\bar J) (V^\top \tilde{D}^{-1/2})^{-1},
\end{align*}
and we have with $\bar x=\tilde{D}^{1/2}Vx$
\begin{align*} 
   x_{n+1}  &=  C(E,\tfrac{h}{2} JQ) x_n  \qquad
   \Leftrightarrow \qquad
     \bar x_{n+1}  =  C(\bar E,\tfrac{h}{2} \bar J) \bar x_n.
\end{align*}
The Cayley transform $C( E,\tfrac{h}{2}JQ) = C( Q^\top E,\tfrac{h}{2}J)$ preserves the energy semi-norm $\| \cdot \|_{Q^\top E}$, i.e., we have $\|C( E,\tfrac{h}{2} J)x\|_{Q^\top E} = \|x \|_{Q^\top E}$, while the Cayley transform $C(\bar E,\tfrac{h}{2}\bar J)$ preserves the energy semi-norm $\|\cdot\|_{\bar E}$, i.e., $\|C(\bar E,\tfrac{h}{2}\bar J) \bar x\|_{\bar{E}} = \| \bar x \|_{\bar{E}}$. The two semi-norms satisfy $\|  x \|_{Q^\top E} =\| \bar x \|_{\bar{E}}$, where $\bar x$ is the linearly transformed vector $\bar x = \tilde{D}^{-1/2}V x$. 
%\hfill $\Box$
\end{remark}

\begin{lemma}\label{lem:midpoint-energy}
Let the $J$-$R$ decomposition of the port-Hamiltonian system \eqref{eq:PHS-General} with at most index~1 and  Assumption~\ref{ass:restriction-index1} be given. The energy-conserving subsystem \eqref{eq:J-R-subs} is in case (a) a DAE of at most index~1, $E\dot{x} = Jx$ with $x(0)=x_0$ (and $E^\top =E\ge 0$, $J = -J^\top$) and in case (b) an ODE.
Then, the implicit midpoint rule is well-defined for any step size $h\ne 0$, i.e.,  $E\pm \frac{h}{2}J$ is invertible for any $h\ne 0$. It is of second order and preserves the energy
\begin{align*}
  \mathcal{H}(x_{n + 1}) =  \mathcal{H}(x_{n}).
 \end{align*}
\end{lemma}

\begin{proof}
As the energy-conserving subsystem has an index of less or equal one, $\{E,J\}$ is a regular matrix pencil. Hence, the midpoint point rule is well-defined for any $h\neq 0$.
Applying Lemma~\ref{lem:Cayleytrafonorm} with  $A=\tfrac{h}{2}J$ yields that $C(E,\tfrac{h}{2}J)$ is a generalized Cayley transform and the energy is conserved in all time points, i.e.,
   \begin{align*}
     2\mathcal{H}(x_{n+1}) & = x_{n+1}^\top E x_{n+1}  = x_{n}^\top (C^\top E C) x_n 
     = 
   x_{n}^\top  E  x_n = 2\mathcal{H}(x_n).
   \end{align*}
 The convergence order of the implicit midpoint rule (1-stage Gauss-Runge-Kutta method) is a classical result from literature, cf., e.g., \cite{hairer2006}. 
 %\hfill $\Box$
\end{proof}

\begin{remark}
 In Lemma~\ref{lem:midpoint-energy}, the implicit ODE case is covered. Note that in this case the energy semi-norm $\|x\|_E=(x^\top E x)^{1/2}$ becomes a norm, i.e., it additionally holds that $\|x\|_E=0$ if and only if $x=0$. 
 %\hfill $\Box$
 \end{remark}
 
 Proceeding from a Strang splitting on top of the $J$-$R$ decomposition, the discretization of both subsystems with the implicit midpoint rule yields a structure-preserving second-order approximation. In particular, the energy in the skew-symmetric subsystem is conserved and the dissipation inequality should generally be satisfied. In certain cases, however, the discretization of the subsystem belonging to $R$ and $B$ with other second-order dissipative schemes (L-stable) is preferable, as we discuss in Section~\ref{sec:numerics-JR}.

%-------------------------------------------
\section{Numerical Results} 
\label{sec:numerics}
%-------------------------------------------

In this section we numerically demonstrate the applicability and performance of our two decomposition strategies in combination with Strang splitting, using port-Hamiltonian benchmark examples from electric circuits. Our focus is on approximations of second order. 

In the examples we use SI units for the electric quantities, i.e., node potentials $e$ are considered in volts (V), currents $j$ in amperes (A), capacitances $C$ in farads (F), resistances $R$ in ohm ($\Omega$), inductances $L$ in henry (H) and time $t$ in seconds (s). For better readability we often suppress the units in the following and in particular use dimensionless step sizes in the visualizations.

For details on the used splitting schemes and numerical integrators we refer to Appendix~\ref{appendix:splitting-approximation}. All computations are performed in MATLAB, version R2023b, on a MacBook Air with Apple M2, 16GB.

%%%%%%%%%%%%%%%
\subsection{\textbf{Dimension-reducing decomposition}}

The dimension-reducing decomposition is computationally advantageous, when dealing with coupled high-dimensional DAEs. In case of DAEs with private substructures, there are additional decoupling strategies that further increase efficiency, as we discuss in the following.

%%%%%%%%%%%%%%EXAMPLE: Coupled index-1 %%%%%%%%%%%%%

\begin{example}[Coupled index-1 DAEs]\label{ex:D_index1} \end{example}
\begin{figure}[b]
\centerline{\includegraphics[scale=0.8]{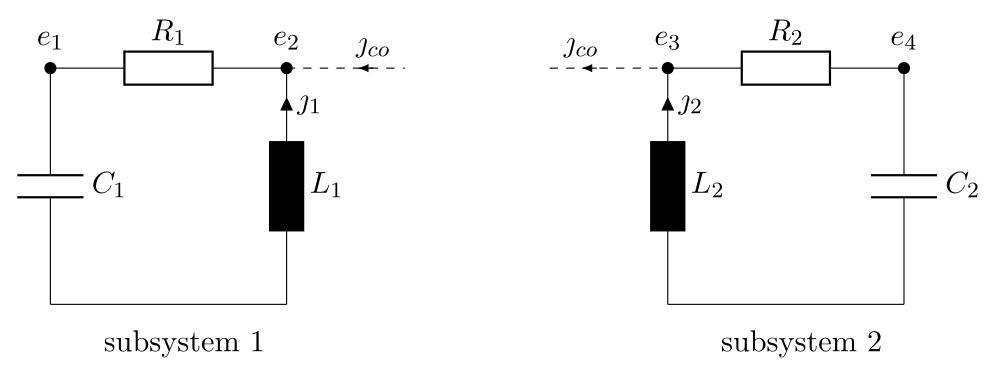}}
 \caption{Example~\ref{ex:D_index1}: LC-oscillator with coupling current $\jmath_{co}$, \cite{Bartel_2023aa}. Parameters:  $C_i= 10^{-5}$~F, $R_i= 10$~$\Omega$ and $L_i= 0.2$~H for $i=1,2$; interval $\mathbb{I}=[0,0.2]$ with the initial values $e_1(0)=e_4(0)=0.1$~V, and $\jmath_1(0)=\jmath_2(0)=1$~A for the dynamic variables; and consistent values for the algebraic variables $e_2(0)=e_3(0)=-9.9$~V and $j_{co}=0$~A.}
\label{fig:LC}
\end{figure} 
The coupled LC-oscillator with damping visualized in  Fig.~\ref{fig:LC} can be described by a model for the node potentials $e_i$, $i=1,...,4$, the currents through inductances $\jmath_1$, $\jmath_2$ as well the coupling current $\jmath_{co}$, where both subsystems and the overall system are of index~1, i.e.,
    \begin{align} 
    \label{eq:LC_osc}
      \text{subsystem 1: } & &  \text{subsystem 2:} \nonumber \\
       C_1\dot{e}_1 - \tfrac{1}{R_1} (e_2 - e_1) &= 0, 
       & \quad  C_2\dot{e}_4 + \tfrac{1}{R_2} (e_4 - e_3) &= 0, \nonumber \\
       \tfrac{1}{R_1} (e_2 - e_1) + \jmath_{1} + \jmath_{co} &= 0, 
       & \quad  -\tfrac{1}{R_2} (e_4 - e_3) + \jmath_{2} - \jmath_{co} 
             &= 0, 
       \\
       L_1 \dot{\jmath}_{1} - e_2 &= 0,  
       & \quad  L_2 \dot{\jmath}_{2} - e_3 &= 0, 
       \nonumber \\
       & 
       & \quad  e_2 - e_3 &= 0, \nonumber
    \end{align}
with capacitances $C_i$, inductances $L_i$ and resistances $R_i$, $i=1,2$.
The analytic solution of the differential variables $x_d^\top=(e_1,e_4,\jmath_{1} , \jmath_{2})$ is given by
\begin{align*}
    x_d(t) &= \text{exp} (M^{-1}At)x_{d,0},\\
    M&=\begin{pmatrix}
    C_1 & & &\\
    & C_2 & &\\
    & & L_1 &\\
    & & & L_2
    \end{pmatrix}, \quad
    A=-\frac{R_1R_2}{R_1+R_2}\begin{pmatrix}
    \phantom{-}\tfrac{1}{R_1R_2} &-\tfrac{1}{R_1R_2}&\tfrac{1}{R_1} &\tfrac{1}{R_1}\\
    -\tfrac{1}{R_1R_2}& \phantom{-}\tfrac{1}{R_1R_2} &\tfrac{1}{R_2} &\tfrac{1}{R_2}\\
    -\tfrac{1}{R_1}& -\tfrac{1}{R_2}& 1 & 1\\
    -\tfrac{1}{R_1}& -\tfrac{1}{R_2}& 1 & 1    
    \end{pmatrix},
\end{align*}
and the algebraic variables $x_a^\top=(e_2,e_3,\jmath_{co})$ are
\begin{align*}
    e_2 &= e_3 = \tfrac{R_1R_2}{R_1+R_2}\big( \tfrac{1}{R_1}e_1 + \tfrac{1}{R_2} e_4 -\jmath_{1} - \jmath_{2} \big), \qquad
    \jmath_{co} =  \tfrac{1}{R_1} (e_1-e_2) - \jmath_{1}.
\end{align*}

\begin{remark}
In the port-Hamiltonian framework the coupled closed system~\eqref{eq:LC_osc} can be formulated  as
\begin{align*}
E^{(i)} \dot x_i &= (J^{(i)}-R^{(i)})x_i + B^{(i)}u_i(t), \quad y_i =(B^{(i)})^\top x_i, \quad i=1,2\\
u_1&=-y_2=j_{co}, \qquad u_2=y_1=-e_2
\end{align*}
for $x_1^\top=(e_1,e_2,\jmath_1)$ and $x_2^\top=(e_3,e_4,\jmath_2,\jmath_{co})$ with
\begin{align*}
& J^{(1)}=\begin{pmatrix}
    	0 \\
    	  & 0 & -1\\
    	  & 1  & 0 
    \end{pmatrix},
 \quad
 R^{(1)}=   \begin{pmatrix}
   	\tfrac{1}{R_1} & -\tfrac{1}{R_1} \\
    -\tfrac{1}{R_1}	& \tfrac{1}{R_1} & \\
   	&   & 0 \\
   \end{pmatrix},
\quad
E^{(1)}=\begin{pmatrix}
 	C_1 \\
 	   & 0 \\
 	   &   &  L_1
 \end{pmatrix},\\
 &J^{(2)}= \begin{pmatrix}
   0& 0  & -1 & 1 \\
   0& 0 &    &\\
   1&   & 0  & \\
   -1&   &    &  0
\end{pmatrix},
\quad
R^{(2)}= \begin{pmatrix}
	\tfrac{1}{R_2} & -\tfrac{1}{R_2} \\
	-\tfrac{1}{R_2}	& \tfrac{1}{R_2} & \\
	&   & 0 \\
	&   &   & 0 \\
\end{pmatrix},
\quad
E^{(2)}=\begin{pmatrix}
	0 \\
	& C_2 \\
	&   &  L_2 \\
	&   &     & 0
\end{pmatrix},
\end{align*}
$(B^{(1)})^\top= ( 0, -1, 0)$ and $(B^{(2)})^\top=(0,0,0,-1)$; or in summarized form as
$$
\begin{pmatrix}
    	E^{(1)} \\
    	   &E^{(2)}  
    \end{pmatrix}
\begin{pmatrix}
    	\dot x_1 \\
    	  \dot x_2  
    \end{pmatrix}
    =
 \left[ \begin{pmatrix}
    	J^{(1)} & -B^{(1)}(B^{(2)})^\top \\
    	B^{(2)}(B^{(1)})^\top    &J^{(2)}  
    \end{pmatrix}  
    -
    \begin{pmatrix}
    	R^{(1)} &\\
    	  &R^{(2)}  
    \end{pmatrix}  \right]
    \begin{pmatrix}
    	x_1 \\
    	 x_2  
    \end{pmatrix}.
    $$
The respective quadratic Hamiltonian $\mathcal{H}$ reads
$$
\mathcal{H}(x)=\tfrac{1}{2} x_1^\top E^{(1)} x_1
     + \tfrac{1}{2} x_2^\top E^{(2)} x_2
     =\tfrac{1}{2} \left( C_1 e_1^2 + C_2 e_4^2 
      + L_1 \jmath_1^2 + L_2 \jmath_2^2 \right). 
$$
%\quad \hfill $\Box$
\end{remark}

\begin{figure}[ht]
    \includegraphics[width=0.525\textwidth]{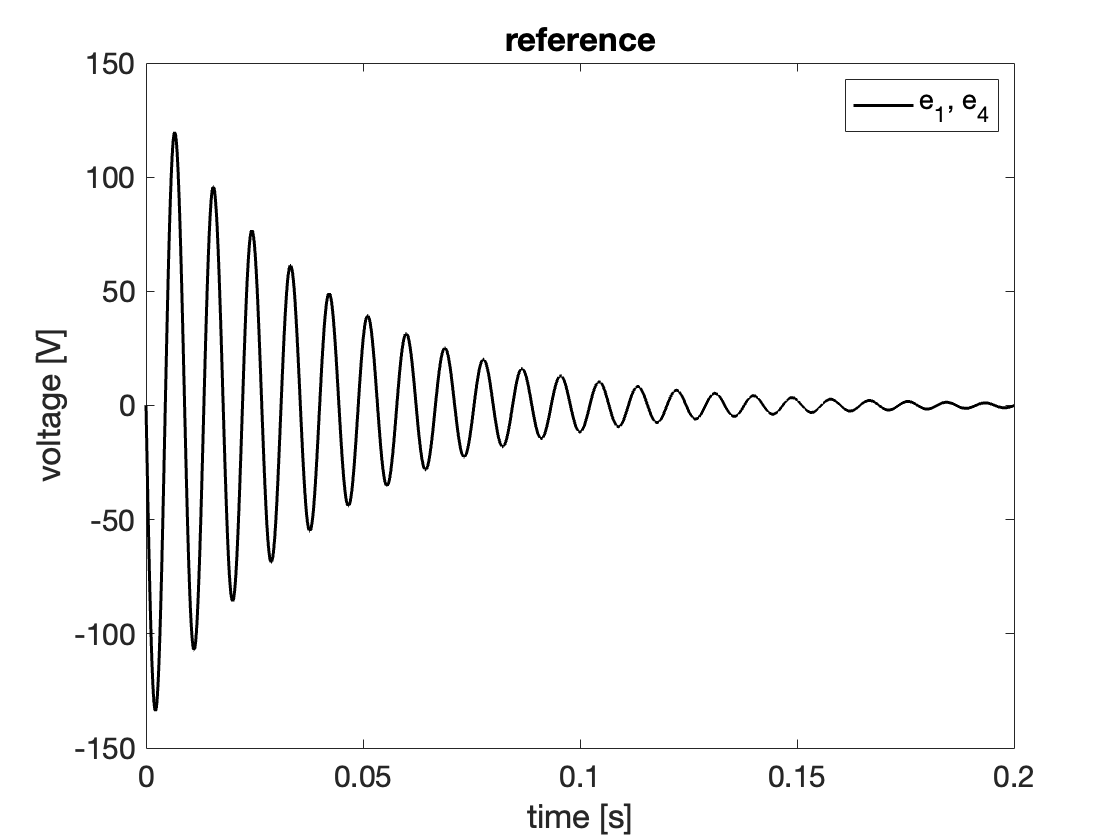}\hspace*{-0.5cm}
     \includegraphics[width=0.525\textwidth]{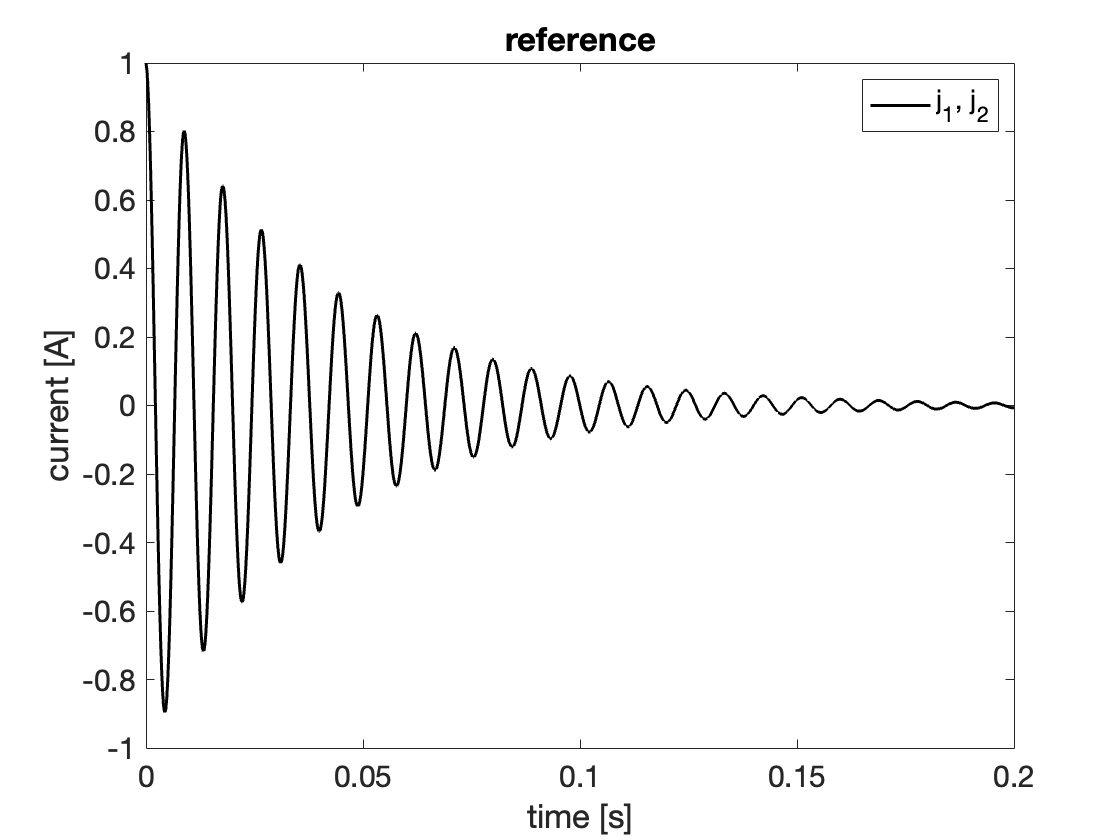}\\
     \includegraphics[width=0.525\textwidth]{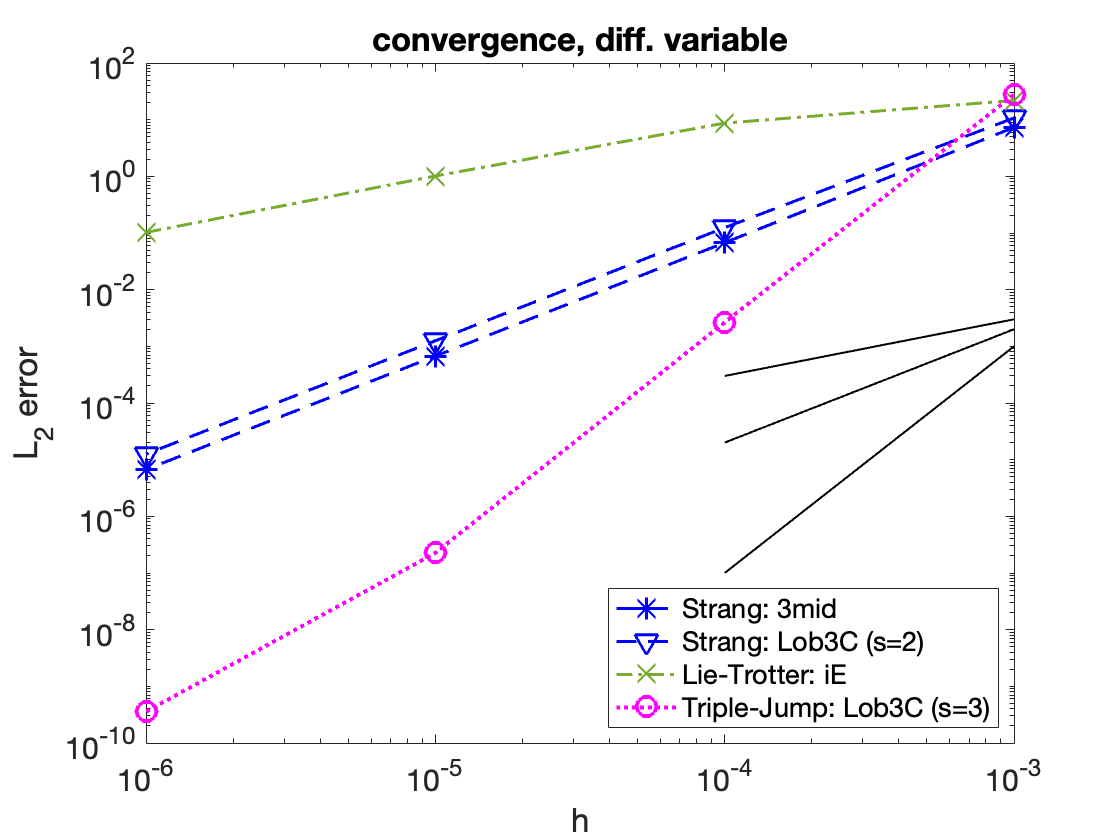}\hspace*{-0.5cm}
     \includegraphics[width=0.525\textwidth]{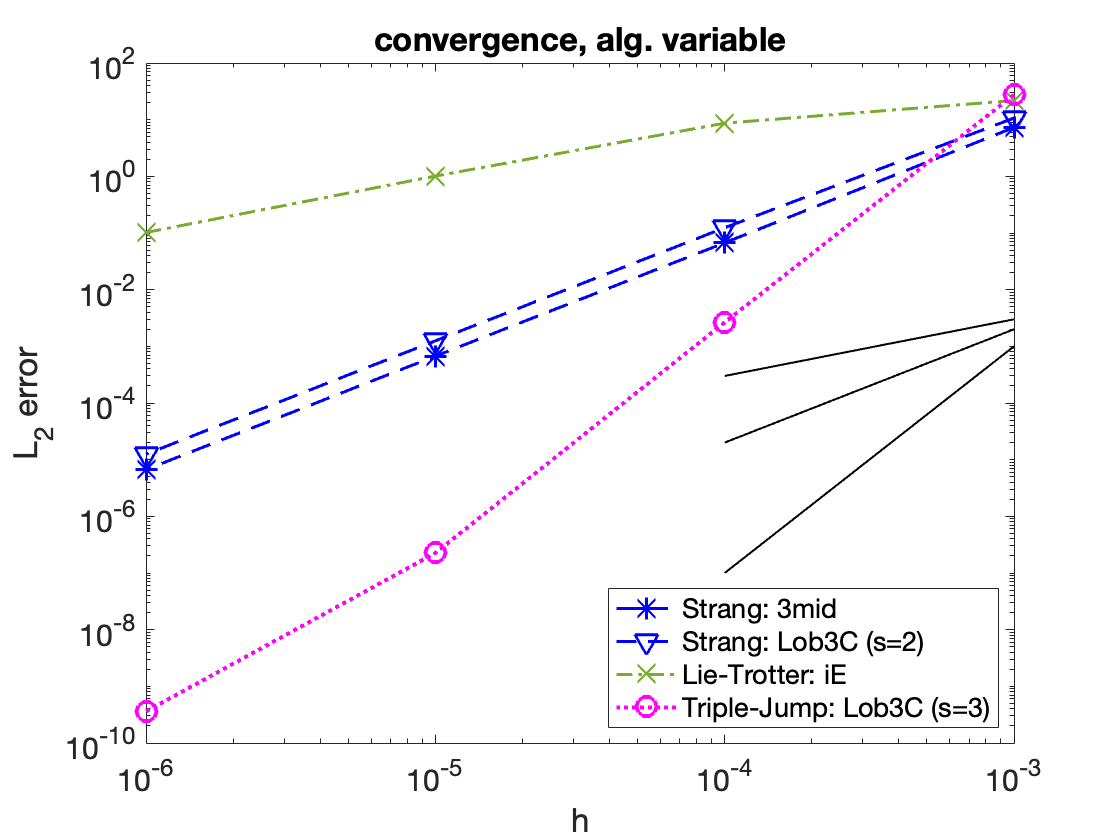}\\
     \includegraphics[width=0.525\textwidth]{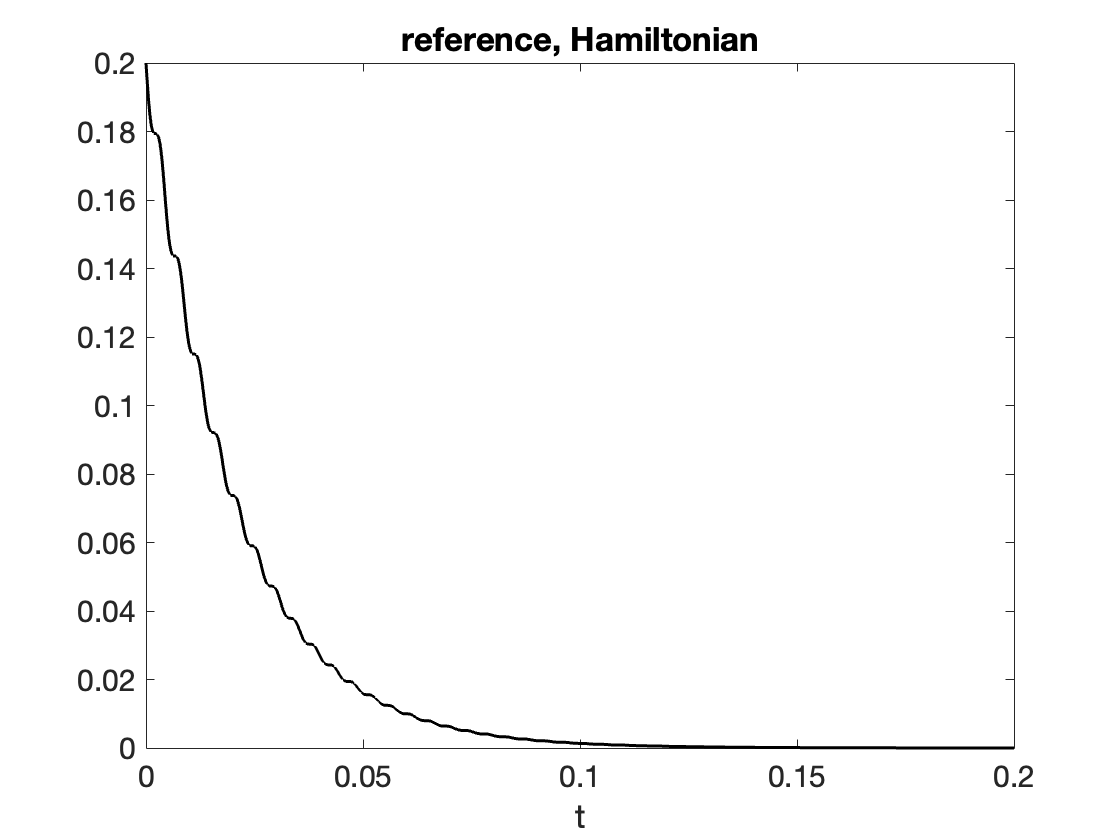}\hspace*{-0.5cm}
      \includegraphics[width=0.525\textwidth]{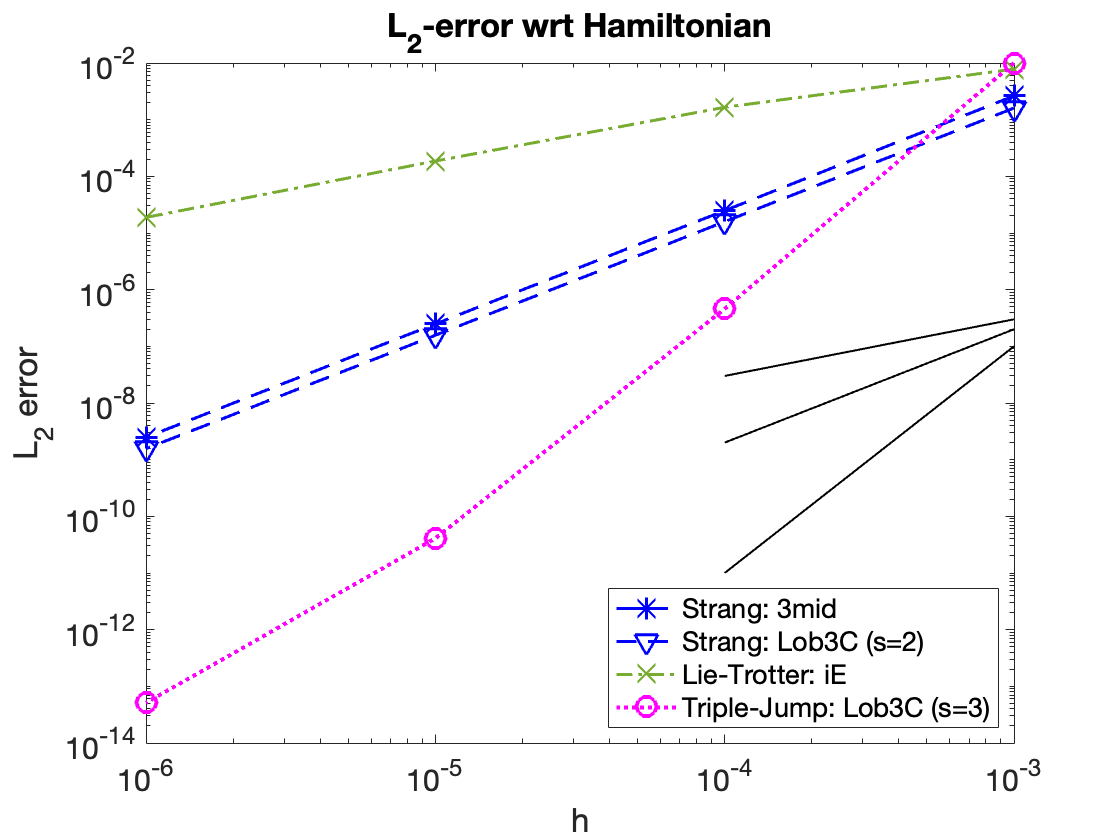}
    \caption{Example~\ref{ex:D_index1}. Top: Reference (analytic) solution for node potentials $e_1$, $e_4$ (left) and currents $\jmath_1$, $\jmath_2$ (right). Middle: Convergence behavior of differential (left) and algebraic variables (right) in some splitting schemes (Lie-Trotter, Strang, Triple-Jump) with appropriate flux approximation for relative step size $h$. Bottom: Hamiltonian (left) and the respective $L^2$-error $\|\mathcal{H}(x_{\mathrm{ref}})-\mathcal{H}(x_{\mathrm{apx},h})\|_{L^2(\mathbb{I})}$ for different numerical approximations (right).}
    \label{fig:exD1_results}
\end{figure}

Figure~\ref{fig:exD1_results} illustrates the analytic solution $x_\mathrm{ref}$ of the node potentials $e_1$, $e_4$ and the currents $\jmath_1$, $\jmath_2$ as well as of the Hamiltonian $\mathcal{H}(x_\mathrm{ref})$ for the setting specified in Fig.~\ref{fig:LC}. Here, the coupling current is zero, and both LC-oscillators show the same damping behavior.

The dimension-reducing decomposition doubles the algebraic equations and considers the respective subsystems as subproblems for the splitting. The splitting scheme applied to the index-1 DAEs inherits the order of convergence from its application to the associated inherent ODEs. An order-preserving numerical discretization of the subproblems is crucial. Figure~\ref{fig:exD1_results} (middle) shows the performance of the second-order Strang splitting where both subproblems are solved with the implicit midpoint rule, i.e.,  A-stable 1-stage Gauss-Runge-Kutta method (3mid), on the one hand, and with the L-stable 2-stage Lobatto-IIIC rule (Lob $s=2$) on the other hand. Both numerical integrators are suitable for index-1 DAEs and yield approximations of second order for the differential and the algebraic variables. The Lie-Trotter splitting combined with the implicit Euler method provides first-order results, and the Triple-Jump splitting combined with the 3-stage Lobatto-IIIC rule gives even fourth-order results, for details on the schemes see Appendix~\ref{appendix:splitting-approximation}. So, the numerical results confirm our theoretical findings.

Let $x_\mathrm{apx}$ be an approximation of order $p$, then the Hamiltonian satisfies
\begin{align*}\|\mathcal{H}(x_\mathrm{ref})-\mathcal{H}(x_\mathrm{apx})\|_{\mathcal{L}^2(\mathbb{I})} &=\|\tfrac{1}{2} (x_\mathrm{ref}-x_\mathrm{apx})^\top E (x_\mathrm{ref}+x_\mathrm{apx})\|_{\mathcal{L}^2(\mathbb{I})}\\
&\leq C \, \|x_\mathrm{ref}-x_\mathrm{apx}\|_{\mathcal{L}^2(\mathbb{I})} =\mathcal{O}(h^p),
\end{align*}
where $h$ denotes the relative step size.
For the setting at hand we numerically observe an error constant $C\sim \mathcal{O}(10^{-4})$ in all approaches, cf.\ Fig.~\ref{fig:exD1_results} (bottom right).

%%%%%%%%%%%%%%%%%%%%%% EXAMPLE 2 %%%%%%%%%%%%%%%%%%%%%%%%

\begin{example}[Coupled DAEs with private index-2 variables] \label{ex:D_index2} \end{example}   

\begin{figure}[b]
\centerline{\includegraphics[scale=0.78]{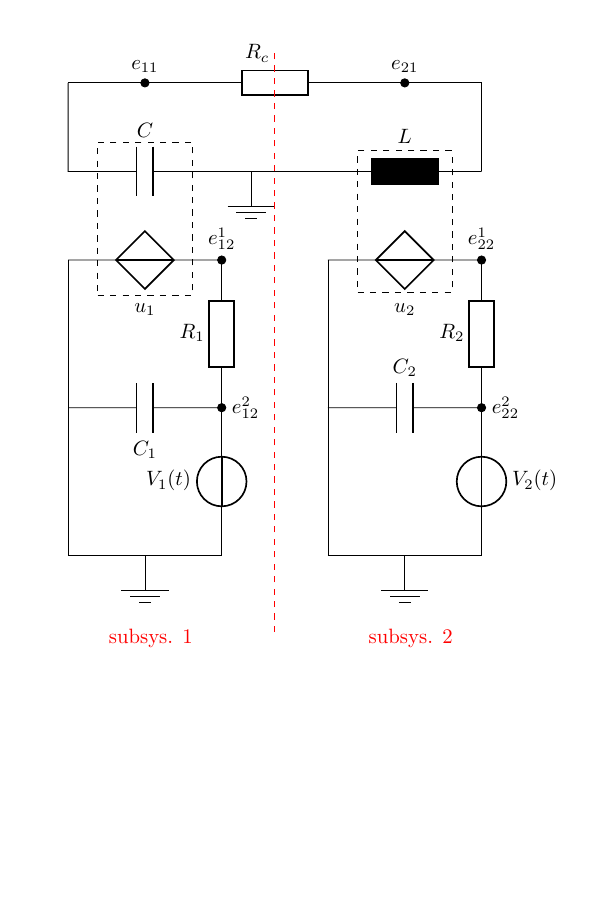}}
        \vspace*{-3cm}
\caption{Example~\ref{ex:D_index2}:
Circuit for two (short) transmission lines with controlled $u_1(t,e_{11})=0.5e_{11}\sin(2\cdot 10^{3}t)$ and $u_2(t,e_{21})=0.5e_{21}\sin(10^{3}t)$ as well as independent voltage sources $v_1(t)=\sin(10^3 t)$~V and $v_2(t)=\sin(3\cdot 10^{3}t)$~V, $t\in \mathbb{I}=[0,4\cdot 10^{-3}]$ and zero initial values except for $e_{21}(0)=-1$~V and $(\jmath,\jmath_{v1},\jmath_{v2})(0)=(1,\tfrac{1}{5},\tfrac{6}{5})$~A.
 Parameters:  $C= 10^{-5}$~F, $R_c= 1$~$\Omega$ and $L= 10^{-2}$~H for the coupling part, moreover $C_1= 2\cdot 10^{-4}$~F, $C_2= 4\cdot 10^{-4}$~F, $R_1= 10^{2}$~$\Omega$ and $R_2= 5\cdot 10^{1}$~$\Omega$ for the private substructures.} \label{fig:ex_D2}
\end{figure}     

The circuit of two transmission lines visualized in Fig.~\ref{fig:ex_D2} can be described by coupled DAEs with private index-2 variables. The unknowns are the current via inductance $\jmath$ in the coupling part as well as the node potentials $e_{i1}$, $e_{i2}^1$, $e_{i2}^2$ and the currents $\jmath_{ui}$, $\jmath_{vi}$ via the controlled $u_i$ and independent voltage sources $v_i$ in the subsystems $i=1,2$.
The controlled voltage sources $u_1$, $u_2$ particularly realize a one-way coupling to the private substructures.

In the notation of Section~\ref{subsec:coupled_index2} we have $x_{11} =e_{11}$, $z_{11} = [\,]$ and $x_{21}=\jmath$, $z_{21} = e_{21}$ as unknowns in the index-1 coupling part, and $x_{i2} =[\,]$, $z_{i2} =[e_{i2}^1, e_{i2}^2, \jmath_{ui}, \jmath_{vi}]$, $i\in\{1,2\}$, as unknowns in the private index-2 substructures. The currents $\jmath_{vi}$ are index-2 variables.
\begin{align*} 
      \text{subsystem 1: } & &  \text{subsystem 2:}  \\
    \text{coupling: \hspace*{0.2cm}}  C\dot{e}_{11} - \tfrac{1}{R} (e_{21} - e_{11}) &= 0, 
       & \quad  L \dot{\jmath} - e_{21} &= 0,  \\
       &  
       & \quad  \tfrac{1}{R} (e_{21} - e_{11}) + \jmath &= 0, \\[2ex]
    \text{substructure:  \hspace*{0.4cm}}   e_{12}^1 - u_1 (t,e_{11}) &=0, \quad &   e_{22}^1 - u_2(t,e_{21})&=0,  \\
     e_{12}^2 - v_1(t) &=0, & \quad    e_{22}^2 - v_2(t)&=0,\\
    \jmath_{u1} - \tfrac{1}{R_1} (e_{12}^1 - e_{12}^2) & = 0,  
       & \quad \jmath_{u2} - \tfrac{1}{R_2} (e_{22}^1 - e_{22}^2)& = 0, \\
    C_1\dot e_{12}^2 - \tfrac{1}{R_1} (e_{12}^1 - e_{12}^2)  - \jmath_{v1} &=0, & \quad
    C_2\dot e_{22}^2 - \tfrac{1}{R_2} (e_{22}^1 - e_{22}^2)  - \jmath_{v2} &=0. 
    \end{align*}

\begin{figure}[ht]
    \includegraphics[width=0.525\textwidth]{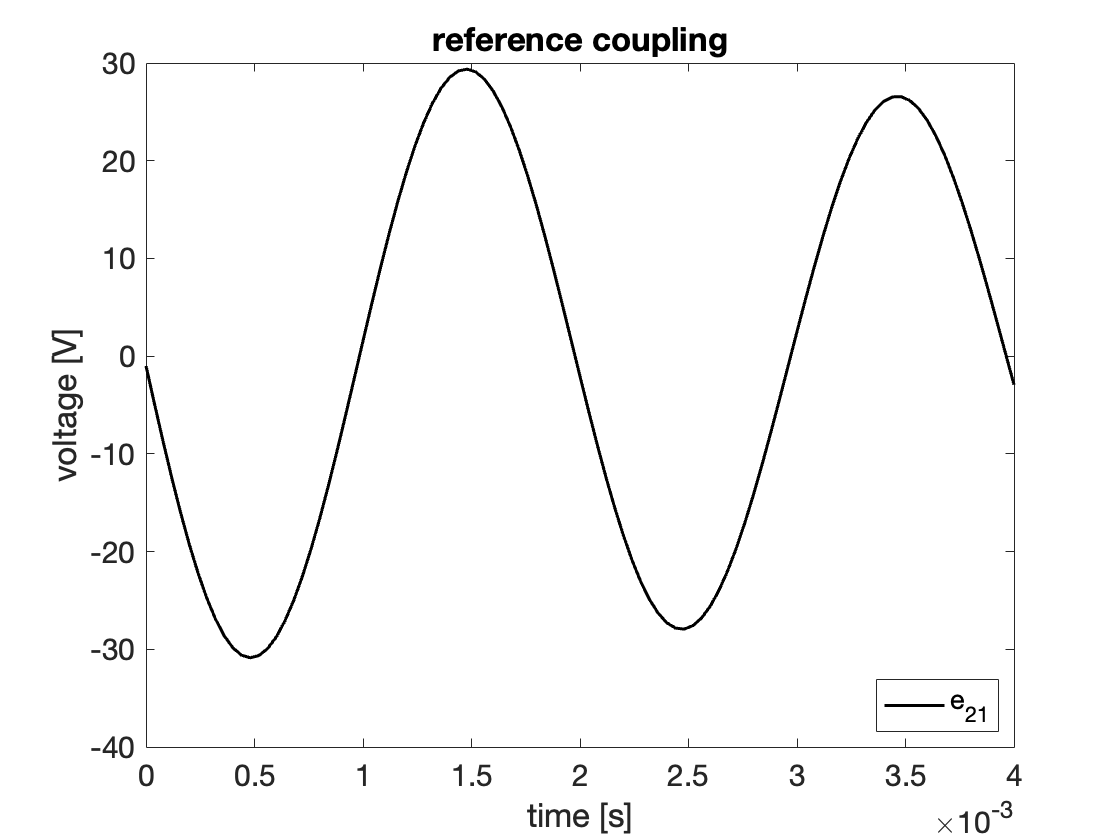}\hspace*{-0.5cm}
     \includegraphics[width=0.525\textwidth]{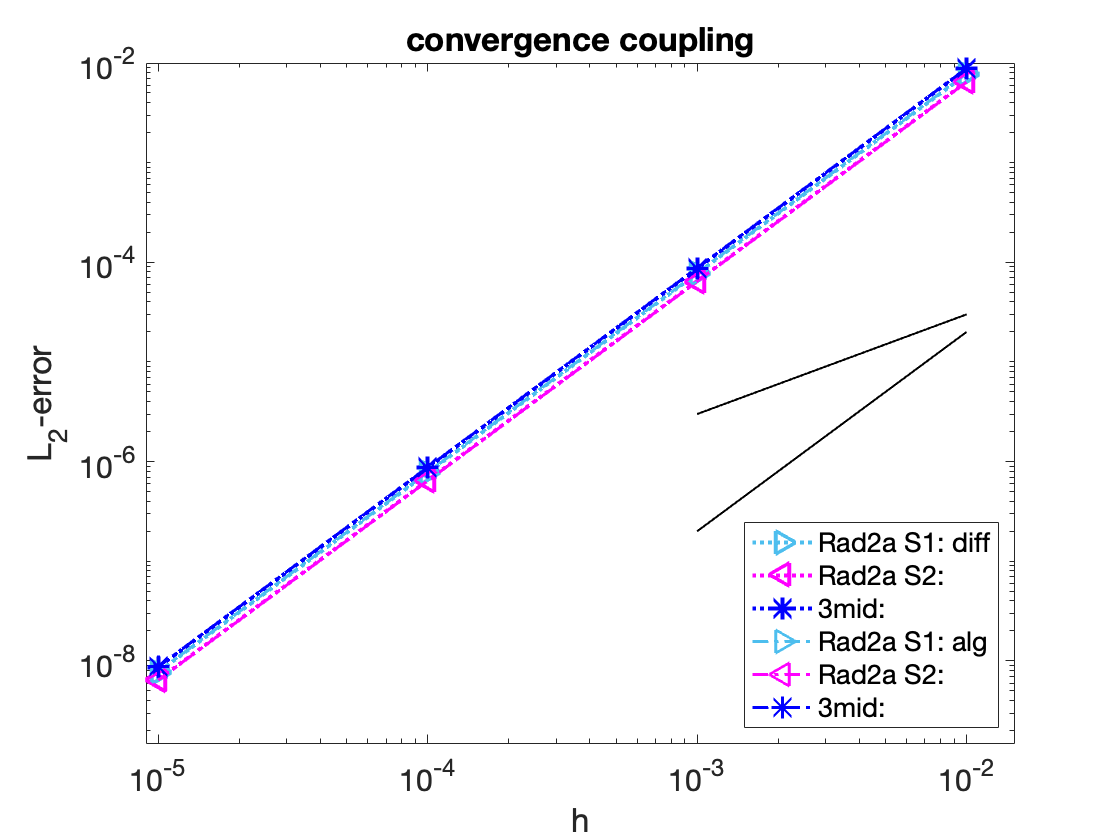}\\
     \includegraphics[width=0.525\textwidth]{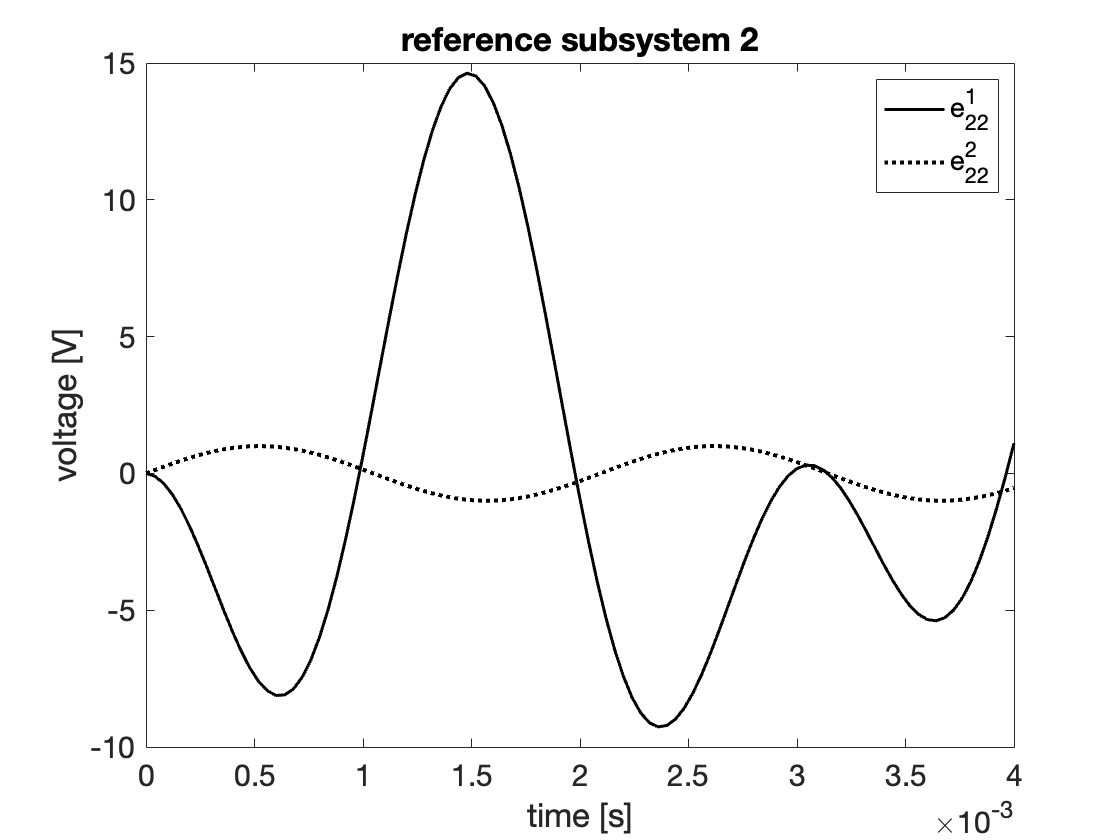}\hspace*{-0.5cm}
     \includegraphics[width=0.525\textwidth]{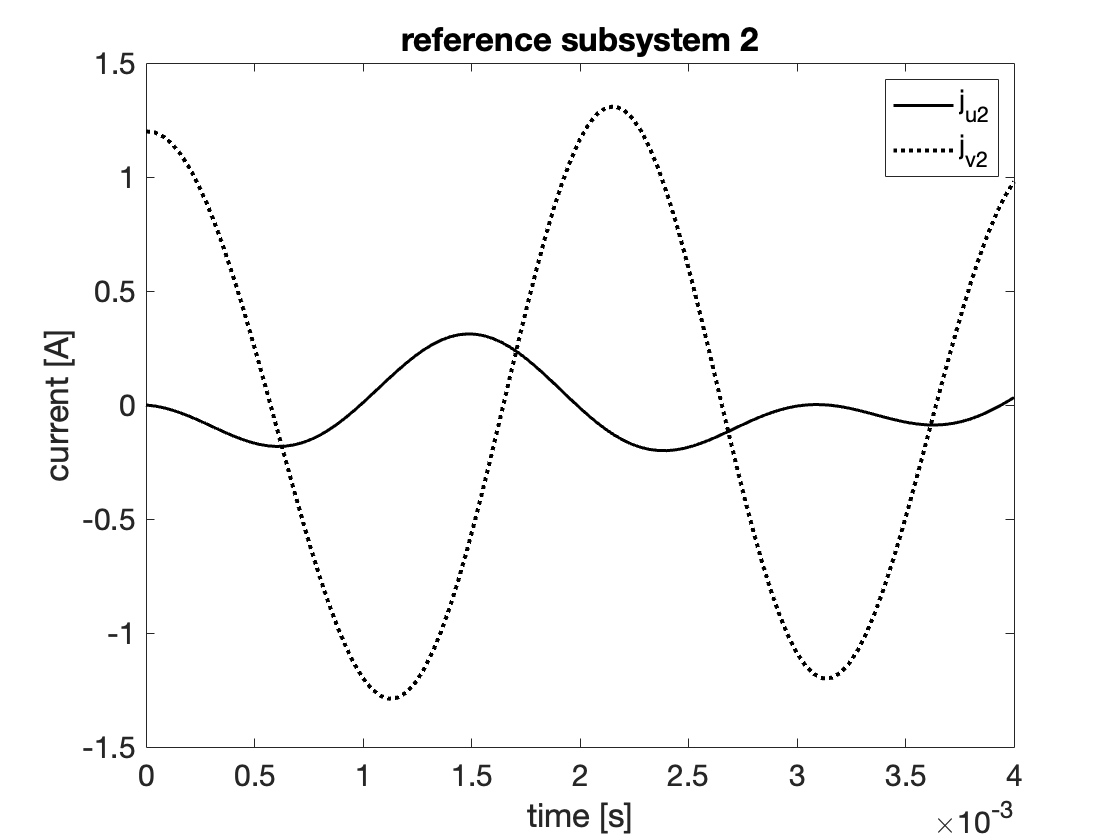}\\
     \includegraphics[width=0.525\textwidth]{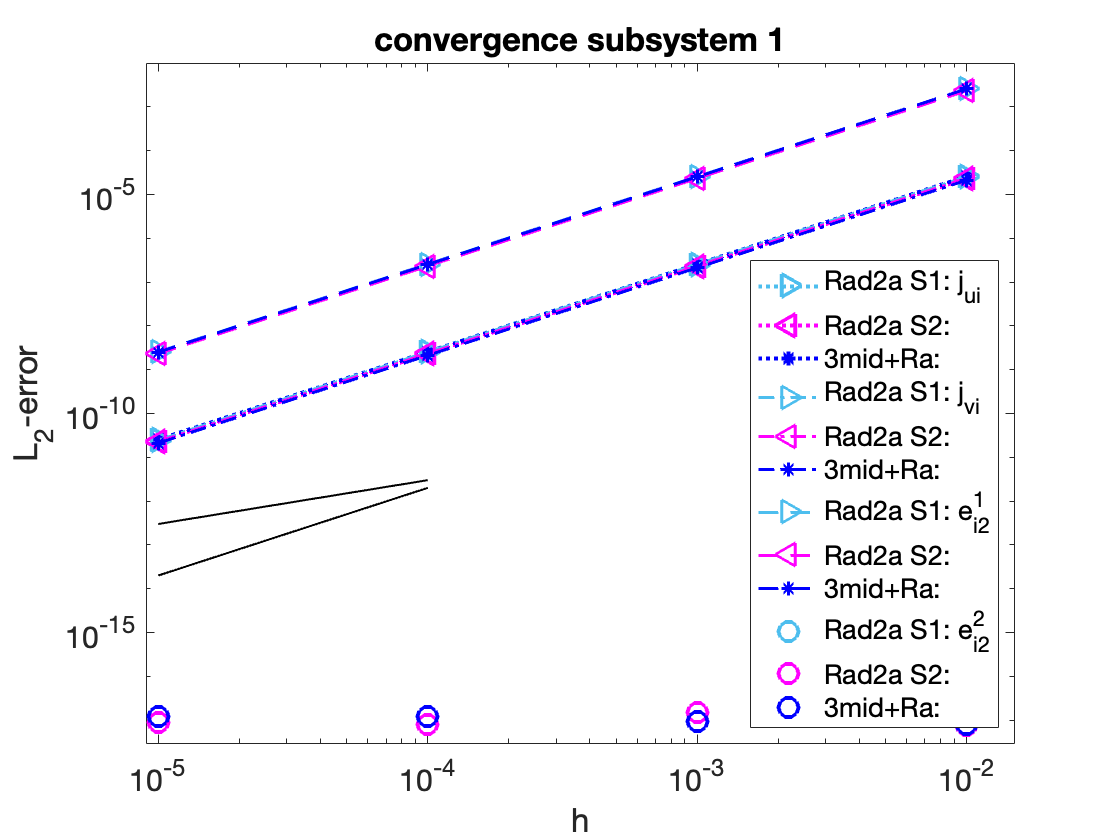}\hspace*{-0.5cm}
      \includegraphics[width=0.525\textwidth]{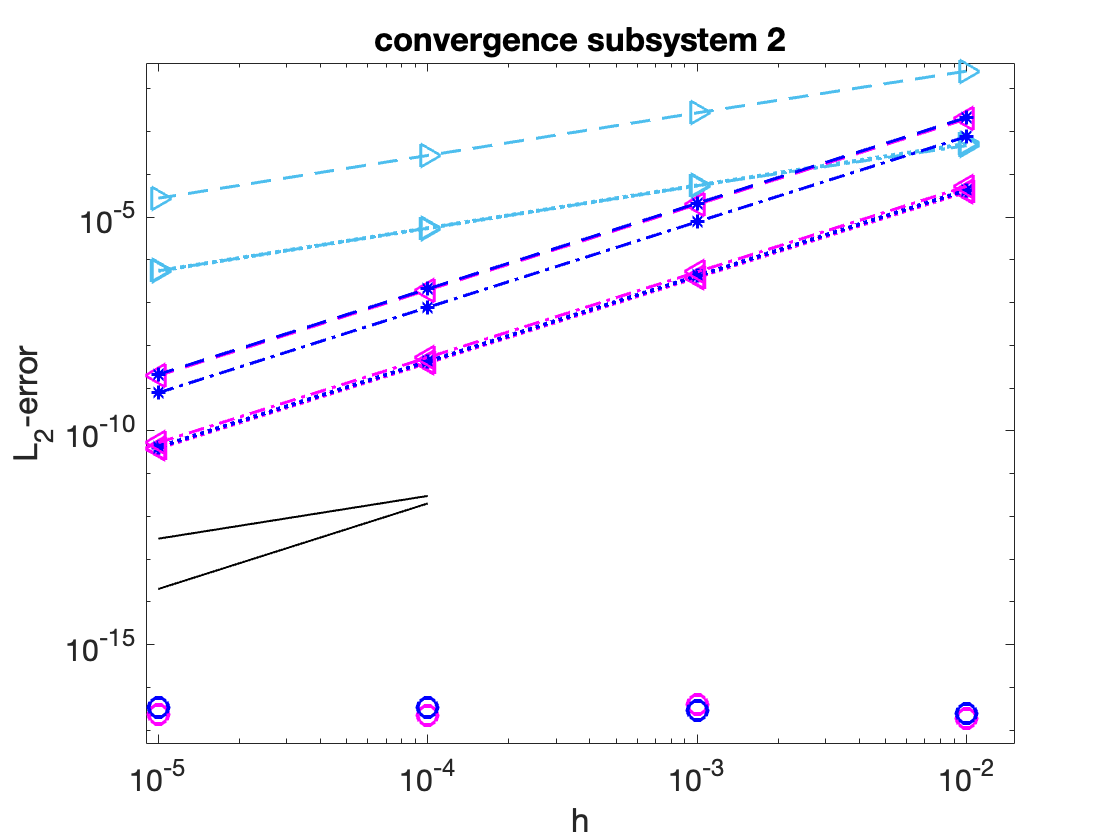}
    \caption{Example~\ref{ex:D_index2}. Top: Reference for node potential $e_{21}$ (index-1 coupling variable) and convergence behavior of coupling variables. Middle: Reference for substructure variables of subsystem~2: node potentials (left), currents (right). Bottom: Convergence behavior of substructure variables for subsystem~1 (left), subsystem~2 (right). Approaches based on Strang splitting: Radau-IIA method (s=2) with subproblem $i$ (Rad2A S$i$) in last splitting step; midpoint rule for coupling part (3mid), on top Radau-IIA for substructures (3mid+Ra).}
    \label{fig:exD2_results}
\end{figure}

For the setting specified in Fig.~\ref{fig:ex_D2}, the reference solutions for the index-1 coupling variable (node potential) $e_{21}$ and for the substructure variables associated to subsystem~2 are illustrated in Fig.~\ref{fig:exD2_results}. The reference is computed for the original coupled system of index-2 DAEs with the 2-stage Radau-IIA rule with relative step size $h=10^{-6}$, i.e., $h=h_T/T$. With focus on second-order approximations we investigate a Strang splitting with different numerical integration strategies. Applying the dimension-reducing decomposition, we observe that the substructure of subsystem~2 depends on the index-1 coupling variable $e_{21}$, whereas the substructure of subsystem~1 relies on the index-0 coupling variable $e_{11}$. Thus, using the 2-stage Radau-IIA rule for both index-2 subproblems and the splitting sequence $S_2$--$S_1$--$S_2$ (where $S_i$ stands for the subproblem $E_i\dot x=\mathrm{f}_i$), we obtain second-order convergence in all variables in accordance with our theoretical results, see (Rad2A S2) in Fig.~\ref{fig:ex_D2}. The substructure variables of subsystem~2 loose convergence order if the splitting sequence is changed into $S_1$--$S_2$--$S_1$, cf.\ (Rad2A S1). The loss is caused by the index-1 coupling variable which is updated in each splitting step and only reaches its accuracy in the last step. In this step, however, $S_1$ is solved and there is no longer any influence on $S_2$. 

Instead of solving the full (index-2) subproblems in the splitting, it is computationally advantageous to consider only the index-1 coupling part of dimension $d_c$ in the splitting (here $d_c=3$) and integrate the two index-2 substructures of lower dimension (here $d_{s1}=d_{s2}=4$) separately afterwards. This decoupling procedure allows for a flexible choice of numerical integrators and maintains the accuracy in all variables of \eqref{eq:private_index2_system} with lower computational effort. We combine here the implicit midpoint rule ($s=1$) for the index-1 coupling part with the Radau-IIA ($s=2$) rule for the substructures, cf.\ (3mid+Ra) in Fig.~\ref{fig:ex_D2}. In this example, the decoupling has the additional highly desirable effect of removing the explicit time dependencies from the splitting, which further increases efficiency.

 %%%%%%%%%%%%%%%%%%%
\subsection{\textbf{Energy-associated $J$-$R$ decomposition}}
\label{sec:numerics-JR}

Let the matrix pencil  $\{E,J-R\}$ be regular, Assumption~\ref{ass:restriction-index1} be fulfilled, and assume $x$ to be a classical solution of the pH-DAE
\begin{align*}
    E\dot{x} & = (J- R) x + Bu(t), \qquad x(0) = x_0, \qquad  y = B^\top  x. 
\end{align*}
According to Proposition~\ref{prop:JR-splitting}, the Strang splitting with step size $h$ given by
\begin{subequations} \label{eq:splittingJR}
\begin{align}
    E_R\,\dot{x}_a & = -R\, x_a + Bu(t), &&x_a(0) = x_0,\\ 
    E_J\,\dot{x}_b & = J\,x_b,             && x_b(0) = x_a(h/2), \\
    E_R\,\dot{x}_c & = -R \,x_c + Bu(t), &&x_c(h/2) = x_b(h),\\
    y_c(h) &= B^\top  x_c(h), \nonumber 
\end{align}
\end{subequations}
provides an approximation of second order, i.e., $x_c(h)=x(h)+\mathcal{O}(h^3)$, 
where $E_A=E$ if $\{E,A\}$ is regular, otherwise $E_A=E+K_E^\top  K_E$ with $K_E$ a projector onto $\mathrm{ker}(E)$ for $A\in \{J,R\}$.  The sequential order of the subproblems (energy-conserving subproblem and dissipative subproblem) could be certainly also swapped. However, the sequential order in \eqref{eq:splittingJR} is advantageous for the numerical discretization as we will see in the following examples, using different schemes depending on the properties of the respective subproblems. Moreover, we study the case when Assumption~\ref{ass:restriction-index1} is not fulfilled.

The selected examples present stiff implicit pH-ODEs as well as {index-1} pH-DAEs. The reference solution to the original system is always computed using the implicit midpoint rule (1-stage Gauss-Runge-Kutta method) with $h=10^{-8}$, ensuring structure-preservation and second-order convergence for differential as well as algebraic variables. Here, $h$ denotes the (dimensionless) relative step size with respect to the considered interval $\mathbb{I}=[0,T]$; i.e., $h=h_T/T$.

%%%%%%%%%%%%%%%%%%%%%% EXAMPLE $J$-$R$ implicit pH-ODE %%%%%%%%%%%%%%%%%

\begin{example}[Implicit pH-ODE] \label{ex:JR-ODE}\end{example}
\begin{figure}[t]
\vspace*{1.75cm}
\begin{tikzpicture} \hspace*{-1cm}
\begin{circuitikz}[scale=.685]
	\draw[fill=black] (0,0) ellipse (.05 and .05) node[left,xshift=-0.05cm]{$e_{1}$};
	\draw (0,3) to[I, l_=$\imath(t)$,i<^=$ $] (0,0);
	\draw[fill=black] (0,3) ellipse (.05 and .05) node[left,xshift=-0.05cm]{$e_{4}$};
	\draw (0,3) to[C, l=$C_R$] (0,6);
	\draw (0,-3) to[C, l=$C_R$] (0,0);
	\draw (0,3) to[R, l=$R_0$] (3,3);
	\draw (3,3) to[L, l=$L_2$,i>^=$\jmath_{2}$] (6,3);
	\draw (0,0) to[R, l=$R_0$] (3,0);
	\draw (3,0) to[L, l=$L_1$,i>^=$\jmath_{1}$] (6,0);
	\draw[fill=black] (3,0) ellipse (.05 and .05) node[above,yshift=0.1cm]{$e_{2}$};
	\draw (3,-3) to[C, l=$C_R$] (3,0);
	\draw (3,3) to[C, l=$C_R$] (3,6);
	\draw[fill=black] (3,3) ellipse (.05 and .05) node[below,yshift=-0.1cm]{$e_{5}$};
	\draw (6,3) to[C, l=$C$] (6,0);
	\draw (8,3) to[R, l=$R_L$] (8,0);
   	 \draw (0,-3) -- (10,-3) -- (10,6) -- (0,6);
    	\draw (6,0) -- (8,0);
   	 \draw (6,3) -- (8,3);
	\draw[fill=black] (7,0) ellipse (.05 and .05) node[above,yshift=0.1cm]{$e_{3}$};
	\draw[fill=black] (7,3) ellipse (.05 and .05) node[below,yshift=-0.1cm]{$e_{6}$};
	 \draw (7,-3) to[C, l=$C_R$] (7,0);
	 \draw (7,3) to[C, l=$C_R$] (7,6);
	 \begin{scope}[yshift=-3cm]
    		\draw (3.,0) -- (3.,-0.6);
    		\draw (2.6,-.6) -- (3.4,-.6);
   	 	\draw (2.75,-.7) -- (3.25,-.7);
    		\draw (2.9,-.8) -- (3.1,-.8);
    	\end{scope}
    	\draw[white] (10.0,7.5) ellipse (.01 and .01);
\end{circuitikz}
\end{tikzpicture}
\vspace{-2.5cm}
\caption{Example~\ref{ex:JR-ODE}: Circuit for two (short) transmission lines with current through voltage source $\imath(t)=0.5 \sin(2\cdot10^7 t)$~A, $t\in \mathbb{I} = [0,10^{-7}]$ and initial value $x(0)=0$. Parameters: $C_R = 10^{-10}$~F, $C =10^{-9}$~F,  $R_0=0.1\,\Omega$,  $R_L=10\,\Omega$, $L_1= 10^{-6}$~H and $L_2=5\cdot 10^{-7}$~H. \label{fig:transmission-circuit}}
\end{figure}

Consider an electric circuit of two transmission lines with crosstalk, Fig.~\ref{fig:transmission-circuit}. The unknown state variables are the potential nodes $e_i$, $i=1,...,6$, and the currents $\jmath_1$ and $\jmath_2$. The input is the current via the voltage source $\imath$. The system behavior can be described by an implicit pH-ODE for  $x^\top= (e_1,e_2,e_3,e_4,e_5,e_6,\jmath_1,\jmath_2)$ with $u(t)=\imath(t)$. The system matrices contain the informations about the interconnectivity in $J$, the resistance in $R$ as well as the capacitance and inductance in $E$, i.e., 
\begin{align*}
J&=\begin{pmatrix}
	0  &    &   &     &   &    &     & \\
	   & 0 & 0 & 0 & 0 & 0 & -1 & 0 \\
	   &   & 0 & 0 & 0 & 0 &  1 & 0 \\
	   &   &   & 0 &   &   &   & \\
	   &   &   &   & 0 &   &   & -1\\
	   &   &   &   &   & 0 &   &  1\\
	   & 1 & -1&   &   &   & 0 &  \\
	   &   &   &   & 1 & -1 &  & 0 \\
\end{pmatrix},
\,\,
R=\begin{pmatrix}
	\tfrac{1}{R_0} & -\tfrac{1}{R_0} & & & & & &\\
	-\tfrac{1}{R_0}& \tfrac{1}{R_0}  & & & & & &\\
	 &  & \tfrac{1}{R_L} & 0& 0 & -\tfrac{1}{R_L} & &\\
	 &  &              & \frac{1}{R_0} & -\tfrac{1}{R_0} & & &  \\
	 &  &              & -\frac{1}{R_0} & \tfrac{1}{R_0} & & & \\
	 &  &-\tfrac{1}{R_L} & 0& 0 & \tfrac{1}{R_L} & & &\\
	 & & & & & & 0 & \\
	 & & & & & &  & 0 \hspace*{-0.2cm}\\
\end{pmatrix},\\
E&=\begin{pmatrix}
	C_R & & & & & & &\\
	& C_R  & & & & & &\\
	&  & C_R+C &    &  & -C  & &\\
	&  &      & C_R & & & &\\
	&  &      &  & C_R & & & \\
	&  & -C   &  &     & C & & \\
	&  &      &  &     &    & L_1 &\\ 
	&  &      &  &     &    &   & L_2\\ 
\end{pmatrix},
\qquad
B=\begin{pmatrix}
	-1 \\ 0 \\ 0 \\ 1 \\ 0 \\ 0 \\ 0 \\ 0 
\end{pmatrix}.
\end{align*}

\begin{figure}[ht]
    \includegraphics[width=0.525\textwidth]{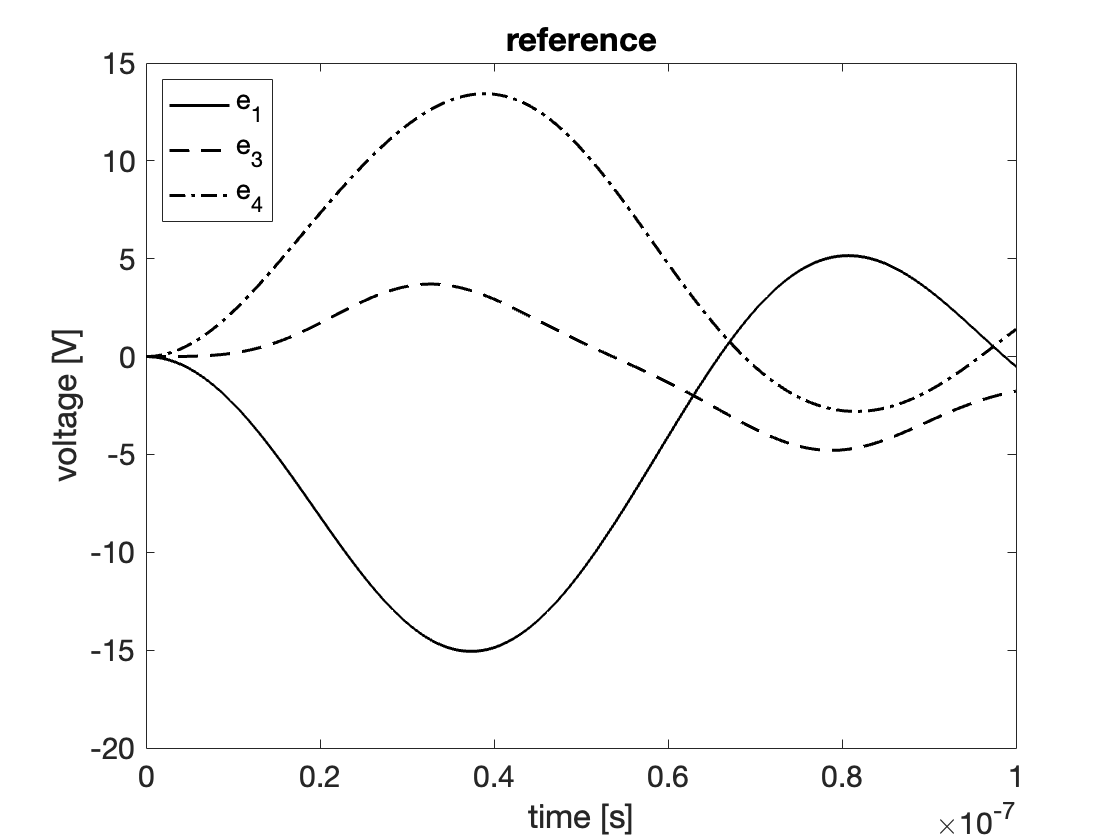}\hspace*{-0.5cm}
     \includegraphics[width=0.525\textwidth]{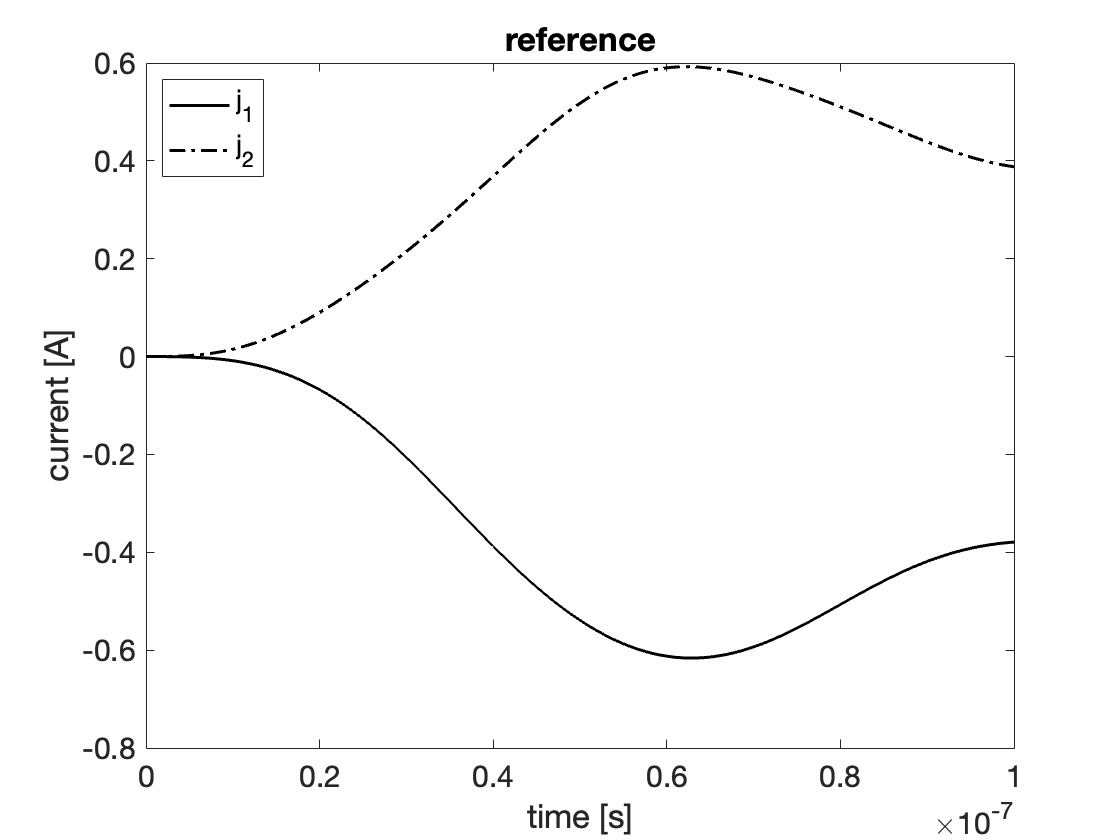}\\
     \centerline{ \includegraphics[width=0.525\textwidth]{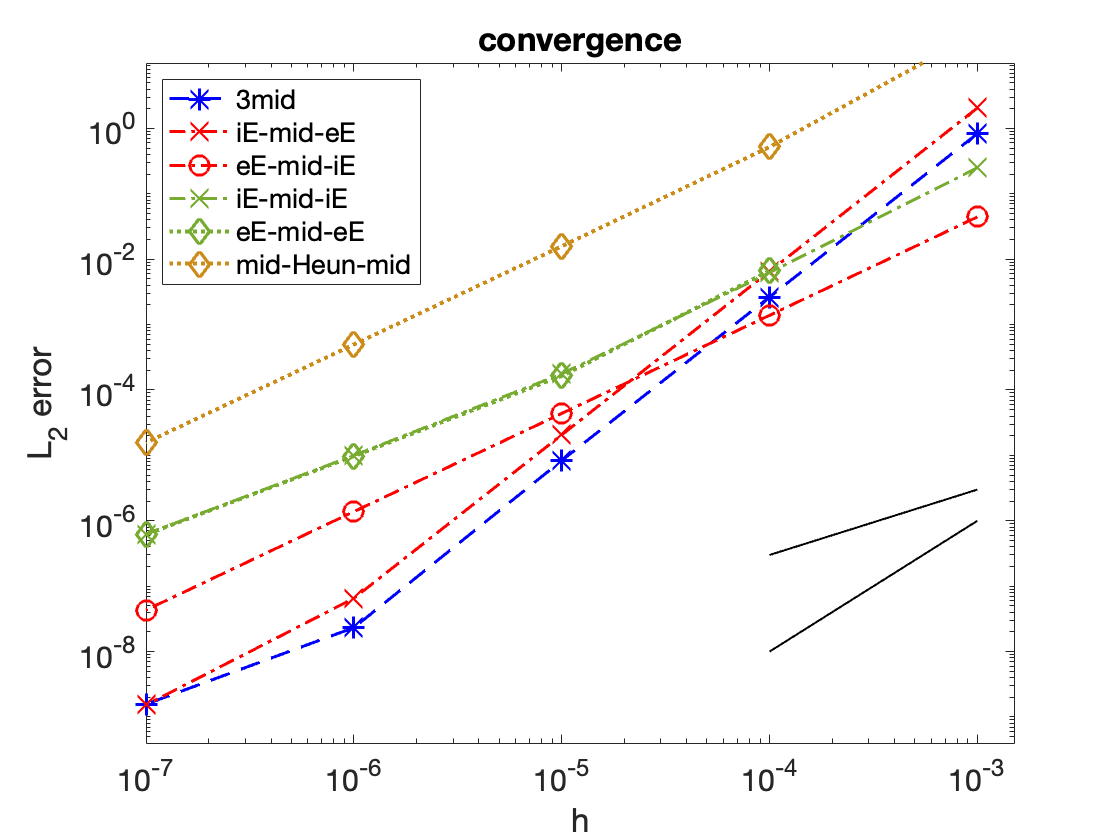}}
     \includegraphics[width=0.525\textwidth]{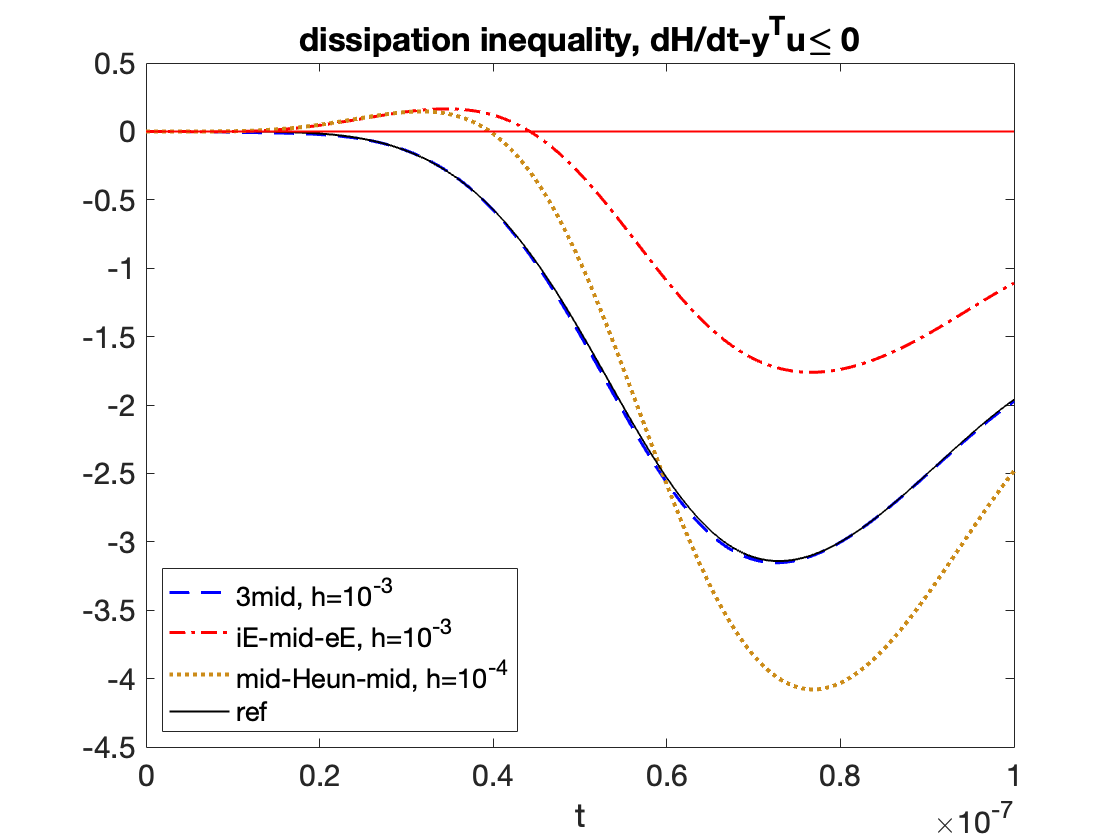}\hspace*{-0.5cm}
      \includegraphics[width=0.525\textwidth]{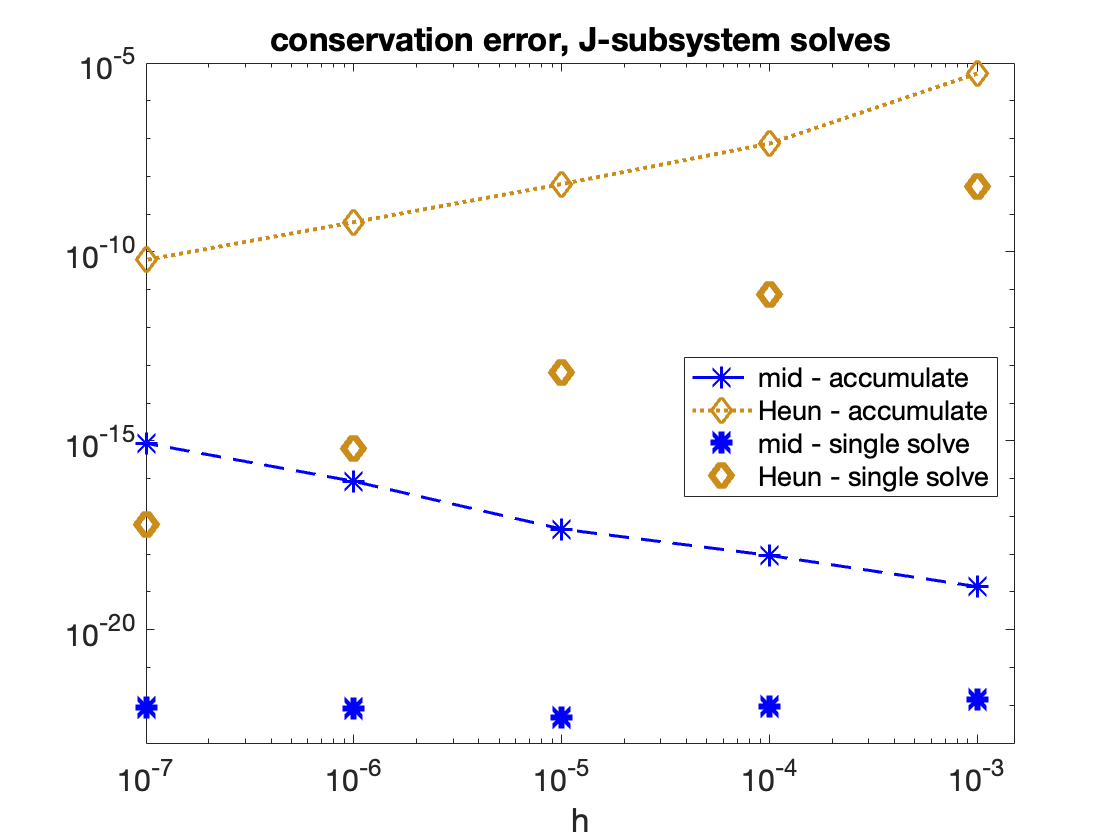}
    \caption{Example~\ref{ex:JR-ODE}. Top: Reference solution for node potentials $e_1$, $e_3$, $e_4$ (left) and currents $\jmath_1$, $\jmath_2$ (right). Middle: Convergence behavior of Strang splitting with different flux approximations for relative step size $h$. Bottom: Evaluation of dissipation inequality for Strang splitting with some second-order flux approximations of fixed step size (left); conservation error in solves of $J$-associated subsystem for implicit midpoint rule and Heun rule: error of single solve (single time step) vs accumulated error over $\mathbb{I}$ for relative step size $h$ (right).}
    \label{fig:exJR-ODE_results}
\end{figure}

In the following study we use a rather academic setup, see Fig.~\ref{fig:transmission-circuit} for details. The reference solution of some node potentials $e_1$, $e_3$, $e_4$ as well as both currents $\jmath_1$, $\jmath_2$ is visualized in Fig.~\ref{fig:exJR-ODE_results} (top).
The symmetric Strang splitting is structure-preserving and second-order convergent. To keep these properties on the discrete level, appropriate numerical discretizations have to be chosen for the subproblems. Applying the implicit midpoint rule to both subsystems (3mid) shows by far the best performance: The convergence order and the dissipation inequality are preserved even for moderate step sizes, see Fig.~\ref{fig:exJR-ODE_results} (middle and bottom left). Note that $h$ denotes the relative step size with respect to the considered interval $\mathbb{I}=[0,T]$, i.e., $h=h_T/T$, here $T=10^{-7}$. Moreover, the numerical conservation error caused by the implicit midpoint rule in the solves of the energy-conserving subsystem associated with $J$ is negligibly small, i.e., we find $|\|x_{n+1}\|_E-\|x_{n}\|_E|\approx \mathcal{O}(10^{-22})$ independently of $h$. This corresponds to a single solve in one time step. The accumulated error for the whole interval $\mathbb{I}$ scales with number of time steps $h^{-1}$, cf. Fig.~\ref{fig:exJR-ODE_results} (bottom right).

The convergence order can also be achieved with other second-order discretizations. Since $E$ is regular, even explicit schemes, such as Runge or Heun methods, or compositions with adjoint methods can be used in \eqref{eq:splittingJR}. Combining a first-order Euler method with the midpoint rule and the respective adjoint Euler method yields a second-order approximation due to symmetry. The two adjoint variants, explicit Euler--midpoint rule--implicit Euler (eE-mid-iE) and implicit Euler--midpoint rule--explicit Euler (iE-mid-eE), are visualized in Fig.~\ref{fig:exJR-ODE_results} (middle). With explicit components, however, the performance is generally worse at moderate step sizes (here $h=10^{-3}$). In some cases, instabilities or violation of the dissipation inequality occurs as a result of the numerical energy supply, see, e.g., the combination of midpoint rule--Heun method--midpoint rule (mid-Heun-mid). Applying the Heun method to the $J$-associated subsystem only ensures energy-conservation, as $h\rightarrow 0$, cf.\ Fig.~\ref{fig:exJR-ODE_results} (bottom right). 
Figure~\ref{fig:exJR-ODE_results} also shows the loss of order due to the use of first-order components without adjoints, see implicit Euler or explicit Euler method for the $R$-associated subsystem with midpoint rule for the $J$-associated subsystem (iE-mid-iE and eE-mid-eE).

%%%%%%%%%%%%%%%%%%%%%% EXAMPLE JR-INDEX 1, DAE for J %%%%%%%%%%%%%%%

\begin{example}[Index-$1$ pH-DAE, Case (a)] \label{ex:JR-DAE-J} \quad \end{example}

\begin{figure}[ht]
\centerline{\includegraphics[width=0.525\textwidth]{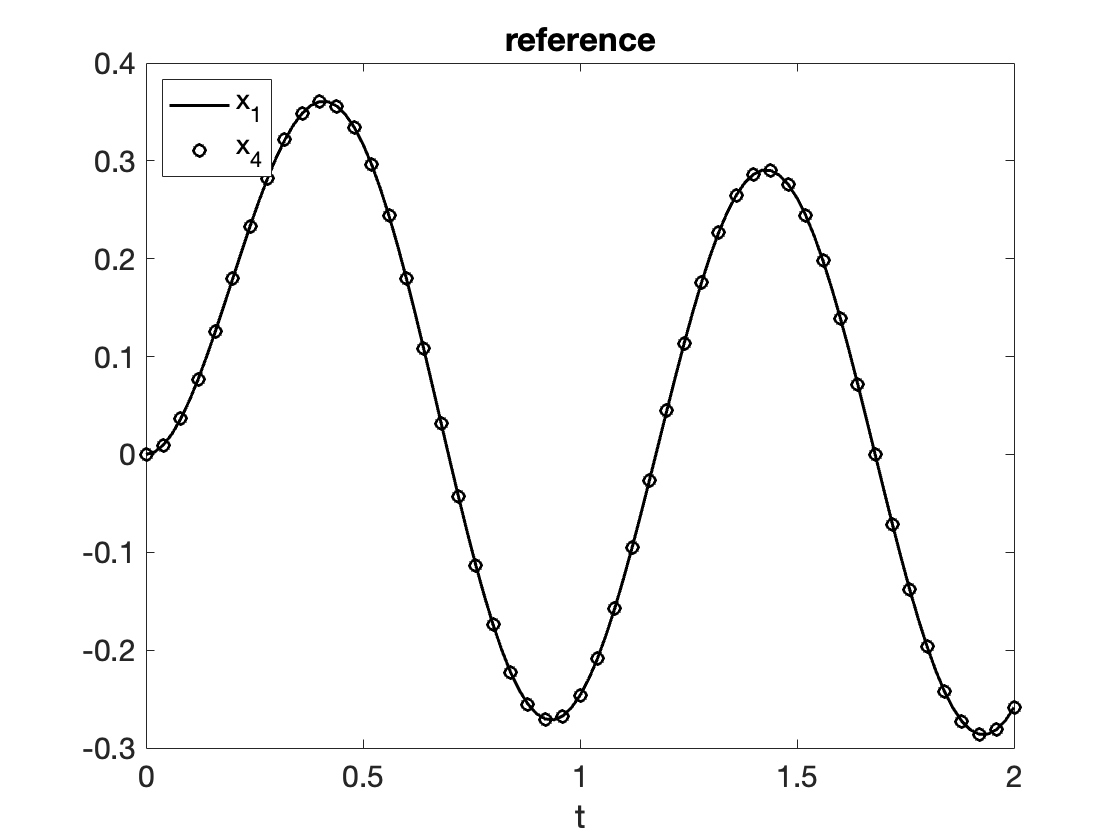}}
     \includegraphics[width=0.525\textwidth]{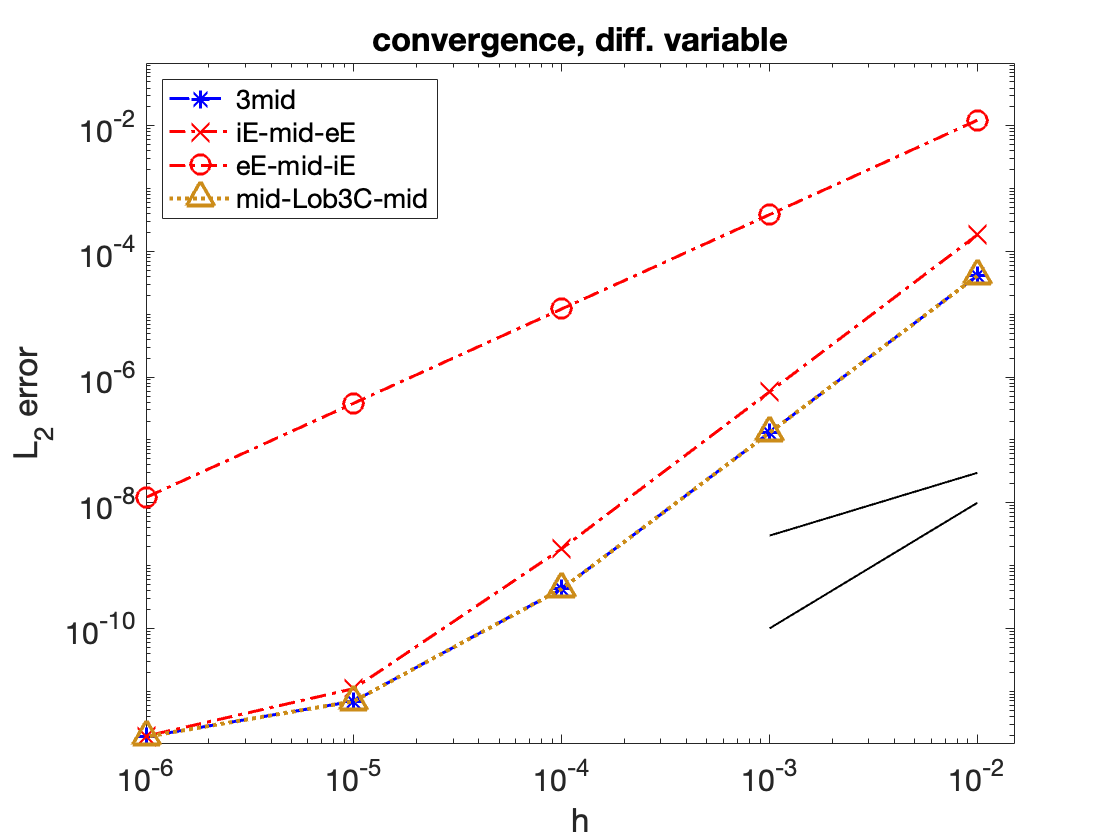}\hspace*{-0.5cm}      	  
      \includegraphics[width=0.525\textwidth]{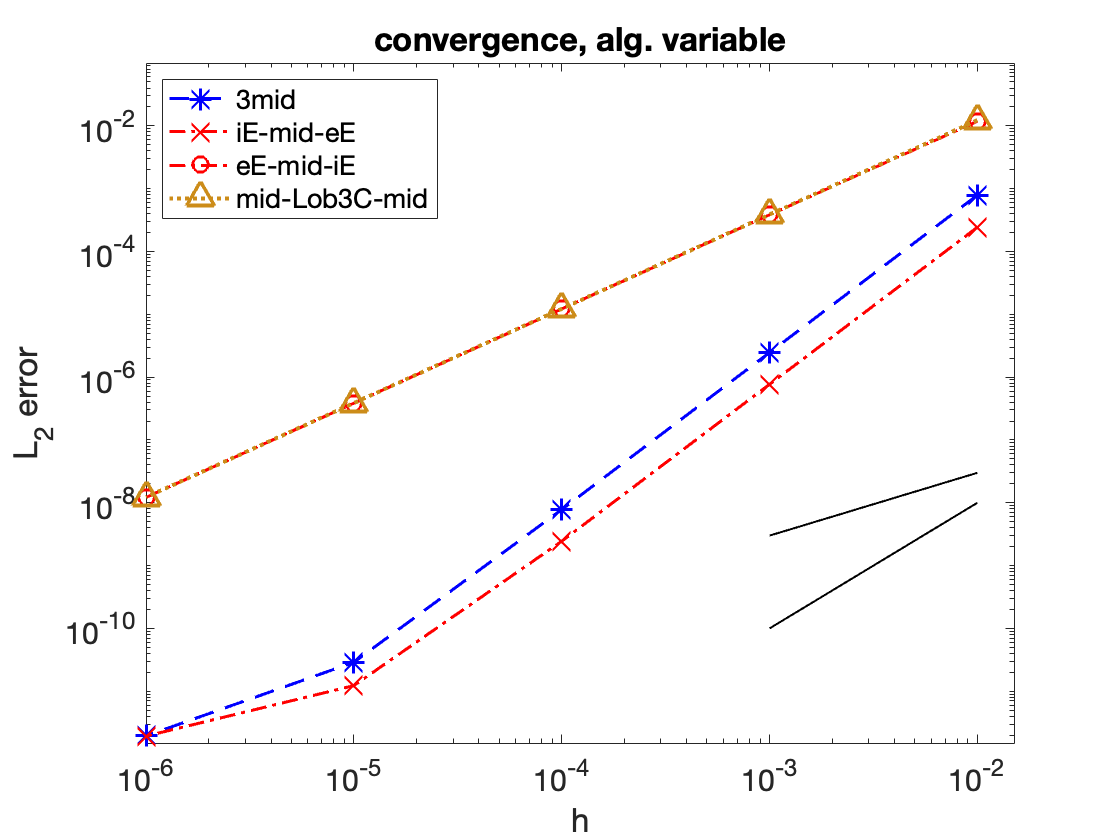}\\
      \includegraphics[width=0.525\textwidth]{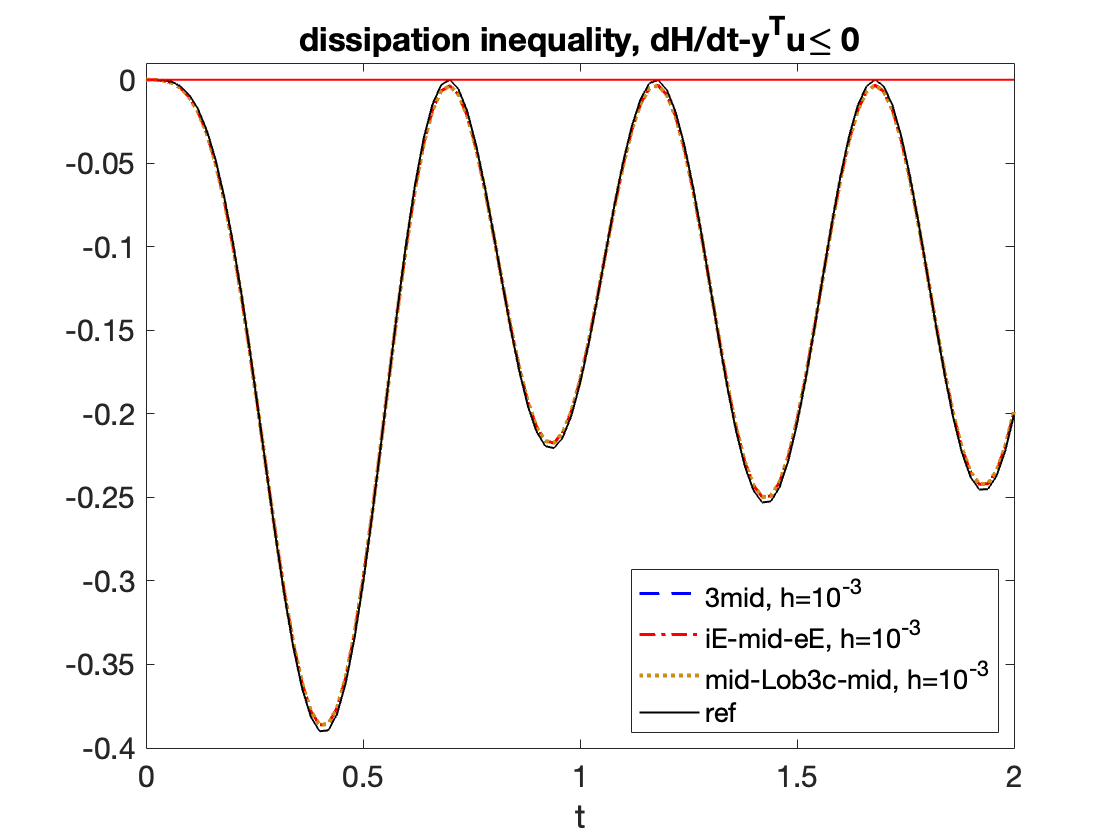}\hspace*{-0.5cm}
       \includegraphics[width=0.525\textwidth]{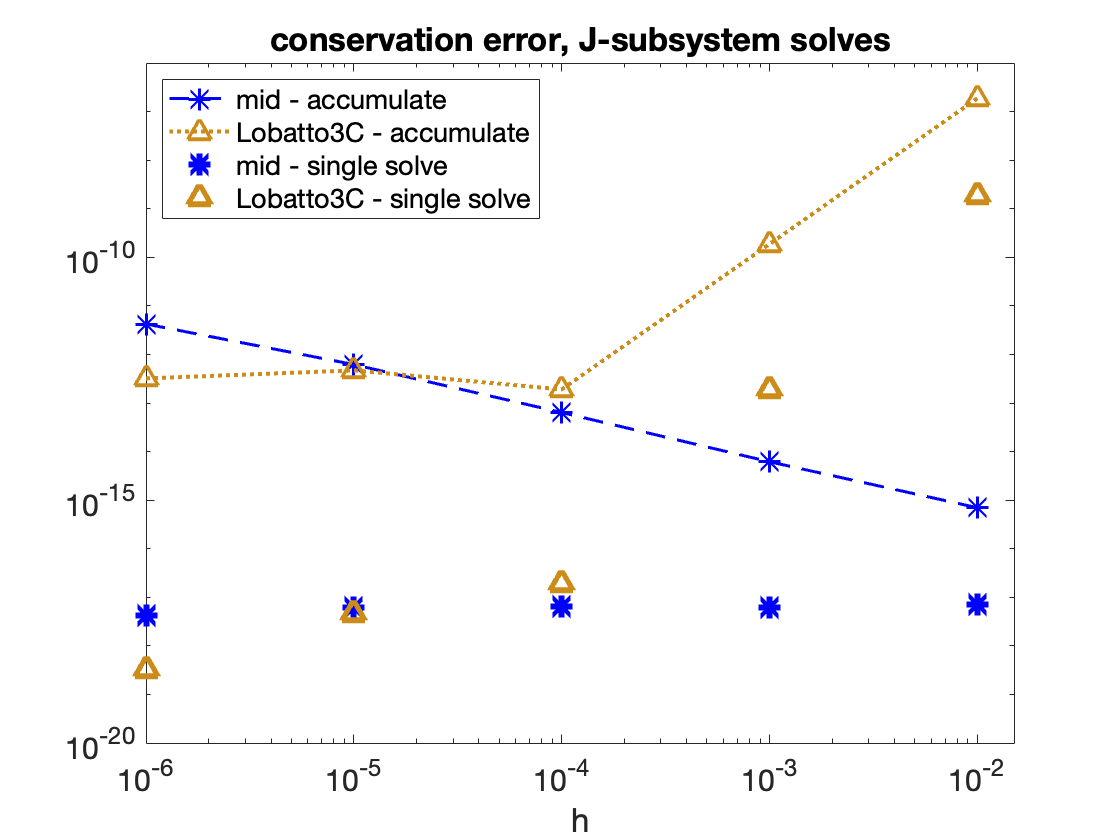}
    \caption{Example~\ref{ex:JR-DAE-J}. Top: Reference solution for differential $x_1$ and algebraic variable $x_4$. Middle: Convergence behavior of differential variables (left) and algebraic variables (right) in Strang splitting with different flux approximations for relative step size $h$. Bottom: Evaluation of dissipation inequality for Strang splitting with some second-order flux approximations of fixed step size (left); conservation error in solves of $J$-associated subsystem for implicit midpoint and 2-stage Lobatto-IIIC rules (right). }
    \label{fig:exJR-DAE-J_results}
\end{figure}

Consider an index-$1$ pH-DAE with the following system matrices,
 \begin{align*}
J=\left(\begin{array}{rrrr}
   0 &   & -1  & \\
     &  0&  1  & -1 \\
   1 & -1&     & -1 \\  
     &  1&  1  &    
\end{array}\right), \quad
R=\begin{pmatrix}
   \phantom{-}3 & -1 \\
   -1& \phantom{-}3\\
     &   & 0 \\
     &   &   & 0
\end{pmatrix}, \quad
E=\begin{pmatrix} 1 \\
                    & 1 \\
                    &   & 0 \\
                    &   &   & 0
  \end{pmatrix}, \quad
B=\begin{pmatrix}
    1 \\ 0 \\ 0 \\ 0
\end{pmatrix},
\end{align*}
input $u(t)=2\sin(2\pi t)$, $t\in \mathbb{I}=[0,2]$ and consistent initial value $x_0 =0$.
The constraint occurs only in the energy-conserving part, hence the matrix pencil $\{E,J\}$ is regular, whereas $\{E,R\}$ is singular and needs to be regularized for the splitting. We get $E_J=E$ and $E_R=E+K_E^\top K_E$ in \eqref{eq:splittingJR}.

The reference solutions for the differential variable $x_1$ and the algebraic variable $x_4$ are illustrated in  Fig.~\ref{fig:exJR-DAE-J_results} (top). Using the Strang splitting with the implicit midpoint rule for both subsystems (3mid) yields second-order approximations that preserve the dissipation inequality of the original pH-DAE and conserve the energy  of the $J$-associated subsystem (index-1 DAE) even for moderate step sizes, see Fig.~\ref{fig:exJR-DAE-J_results} (middle and bottom). Since the $R$-associated subsystem is a (regularized) ODE and solved in the first and last splitting step of the scheme \eqref{eq:splittingJR}, also compositions with adjoint methods can be used as in Example~\ref{ex:JR-ODE}. The  implicit Euler--midpoint rule--explicit Euler (iE-mid-eE) variant shows a similar convergence behavior as 3mid with a significantly smaller error constant than the other adjoint variant, eE-mid-iE. 
The use of a two-stage Lobatto-IIIC rule for the $J$-associated subsystem, the index-1-DAE, and the midpoint rule for the $R$-associated subsystem (mid-Lob3C-mid) also leads to second-order convergence for differential and algebraic variables. The differential variables perform here even slightly better than the algebraic variables.
The Lobatto-IIIC methods are L-stable, B-stable and suitable for stiff problems and DAEs. However, they are not symmetric and ensure the energy-conservation of the $J$-associated subsystem only as  $h\rightarrow 0$, Fig.~\ref{fig:exJR-DAE-J_results} (bottom right). 

\begin{remark}
When modeling electric networks via modified nodal analysis, there is no index-1 configuration that allows a $J$-$R$ decomposition according to Assumption~\ref{ass:restriction-index1}, case (a), see Appendix~\ref{appendix:circuits-and-JRsplitting} for details. 
%\hfill $\Box$
\end{remark}

%%%%%%%%%%%%%%%%%%%%%% EXAMPLE JR-INDEX 1, DAE for R %%%%%%%%%%%%%%%

\begin{example}[Index-$1$ pH-DAE, Case (b)] \label{ex:JR-DAE-R} \quad \end{example}
 \begin{figure}[t]
	 	\begin{center}
		 \begin{tikzpicture}
		 \begin{circuitikz}
			    	\draw (0,0) to[L, l=${L}$,i>^=$\jmath$] (0,3); 
			    	\draw (0,3) -- (2,3);
			 	\draw (2,0) to[C, l=$C$] (2,3);
			 	\draw[fill=black] (2,3) ellipse (.06 and .06) node[above,yshift=0.1cm]{$e_{1}$};
			 	\draw (2,3) to[R, l=$R_1$] (4.5,3);
			 	\draw (4.5,3) to[R, l=$R_2$] (4.5,0);
			 	\draw[fill=black] (4.5,3) ellipse (.06 and .06) node[above,yshift=0.1cm]{$e_{2}$};
			    	\draw (4.5,3) -- (6.5,3);
			 	\draw (6.5,3) to[I, l=$\imath(t)$] (6.5,0);
			    	\draw (0,0) -- (6.5,0);
		 \end{circuitikz}
		 \end{tikzpicture}
	    	\end{center}
\caption{Example~\ref{ex:JR-DAE-R}: Electric circuit with current through voltage source $\imath (t) = 5 \sin{(10^2 t)}$~A, $t\in \mathbb{I}=[0,1]$ and initial value $x(0)=0$. Parameters: $C= 10^{-4}$~F, $R_1 =R_2 = 1~\Omega$ and $L = 0.2$~H.} \label{fig:JR-network-index1}
 \end{figure}
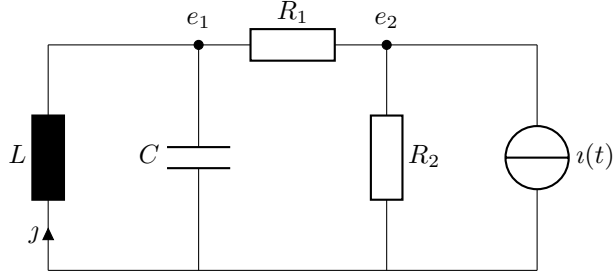
 
 Modeling the electric network from Fig.~\ref{fig:JR-network-index1} via modified nodal analysis yields a pH-DAE of index-1 for $x^\top=(e_1, \jmath, e_2)$ with the system matrices
 \begin{align*}
 J&= \begin{pmatrix}
             0 & -1 &\\
              1&  \phantom{-}  0 \\
                &   & 0\\
           \end{pmatrix}, 
           \quad
 R= \begin{pmatrix}
              \phantom{-} \tfrac{1}{R_1} & & - \tfrac{1}{R_1}
              \\
              &\phantom{.} 0 \phantom{.}& 
              \\
              -\tfrac{1}{R_1} & & \tfrac{1}{R_1}+\tfrac{1}{R_2}
           \end{pmatrix},
           \quad
E=\begin{pmatrix}
    C\\ & L \\ & & 0 
	\end{pmatrix},
	\quad
B=   \begin{pmatrix}
            0 \\ 0 \\ 1
        \end{pmatrix}
\end{align*}
and input $u(t)=\imath(t)$, $t\in \mathbb{I}=[0,T]$. The nodal potential $e_1$ and the current $\jmath$ are differentiable variables, and the nodal potential $e_2$ is the algebraic variable to the constraint that only occurs in the dissipative subsystem with the sources. Hence, the matrix pencil $\{E,R\}$ is regular, whereas $\{E,J\}$ is singular and needs to be regularized for the splitting. With $E_R=E$ and  $E_J=E+K_E^\top K_E$ in \eqref{eq:splittingJR}, the $R$-associated subsystem is an index-1 DAE and the $J$-associated subsystem is a (regularized) ODE. 

\begin{figure}[t]
    \includegraphics[width=0.525\textwidth]{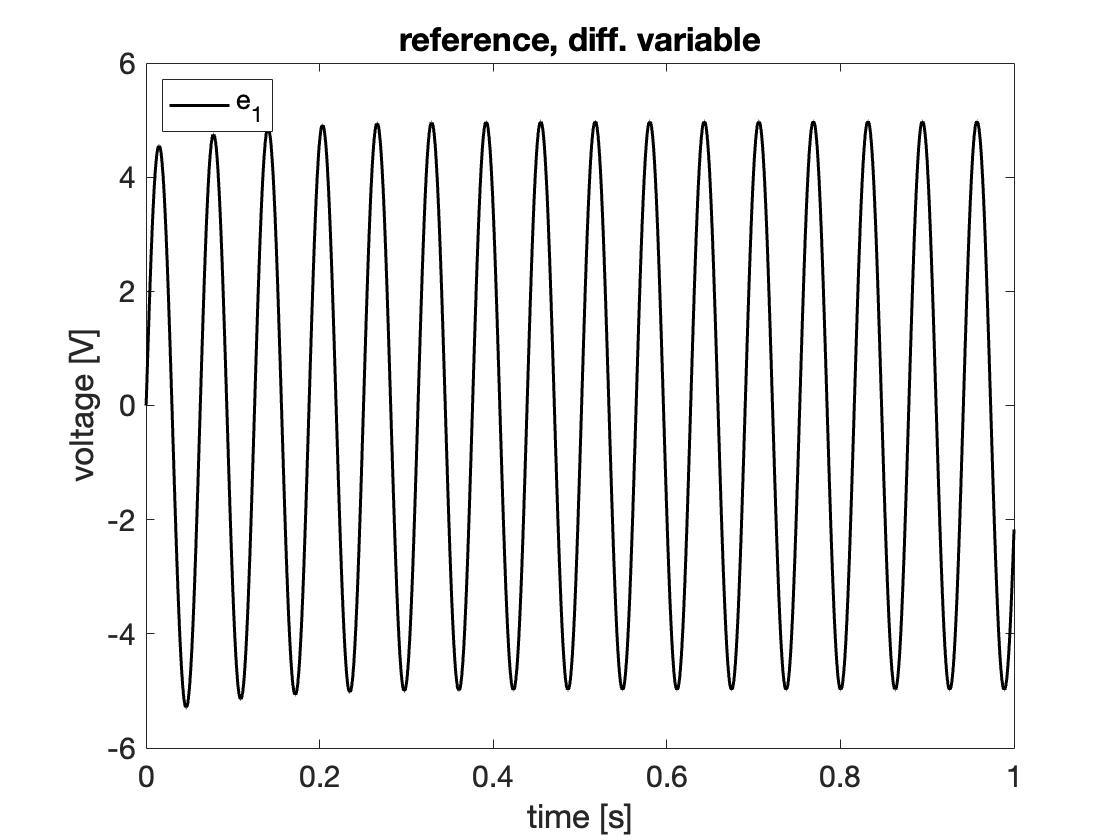}\hspace*{-0.5cm}
     \includegraphics[width=0.525\textwidth]{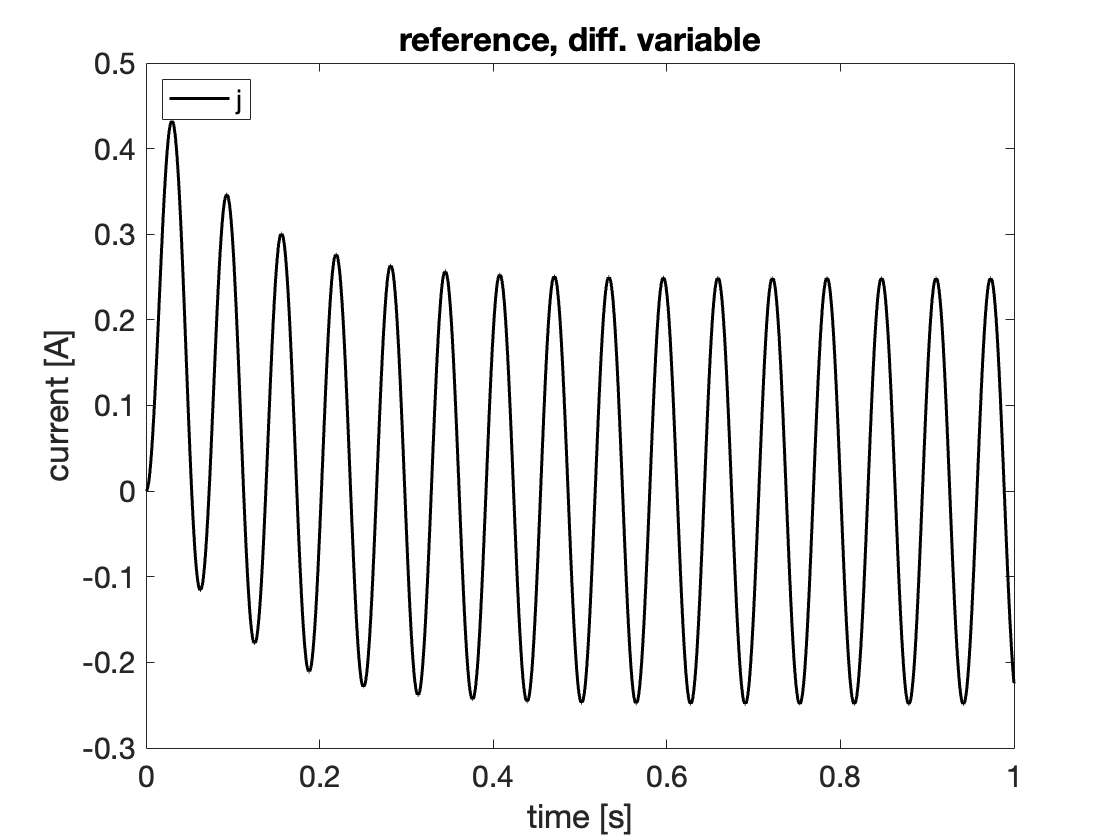}\\
      \includegraphics[width=0.525\textwidth]{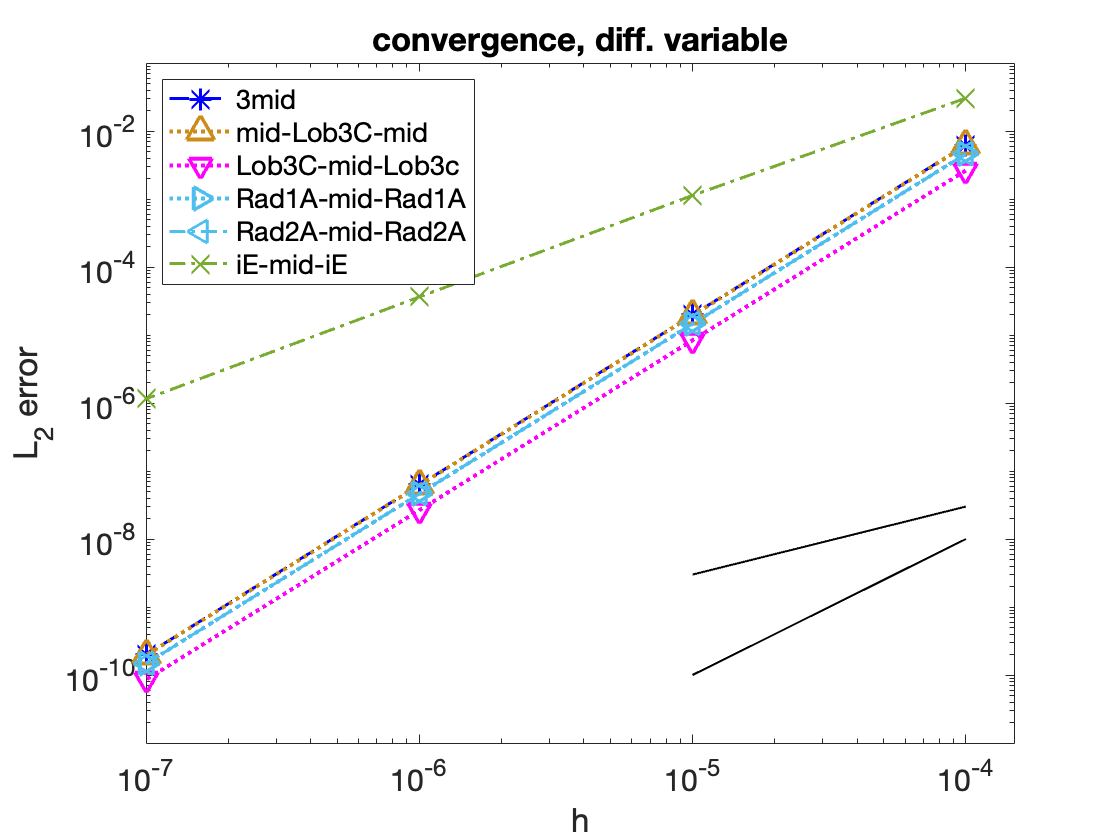}\hspace*{-0.5cm}
      \includegraphics[width=0.525\textwidth]{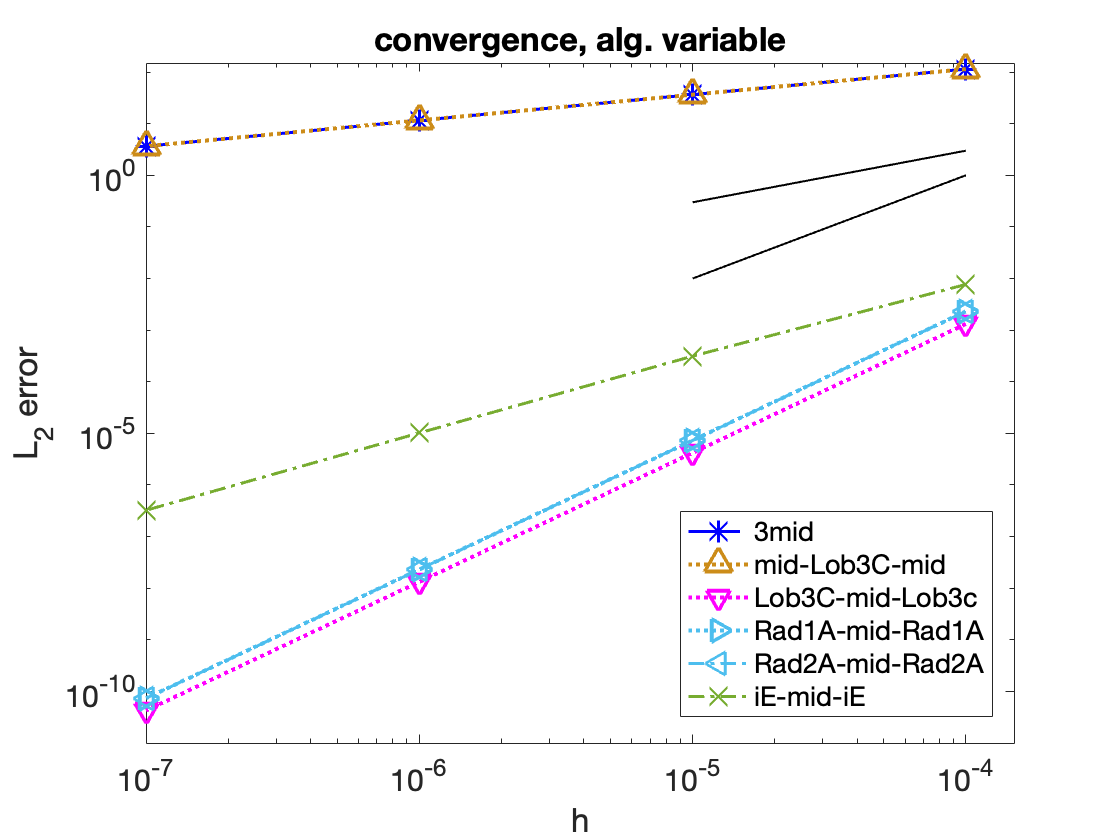}\\
     \includegraphics[width=0.525\textwidth]{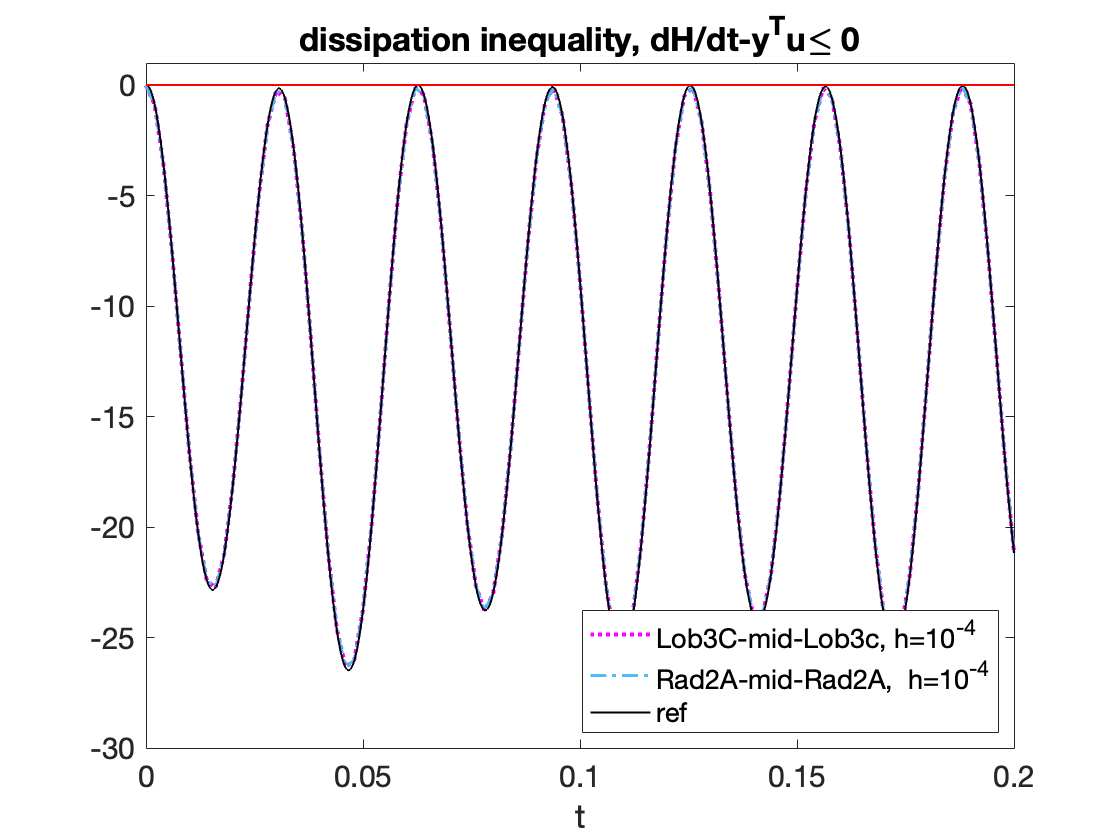}\hspace*{-0.5cm}
     \includegraphics[width=0.525\textwidth]{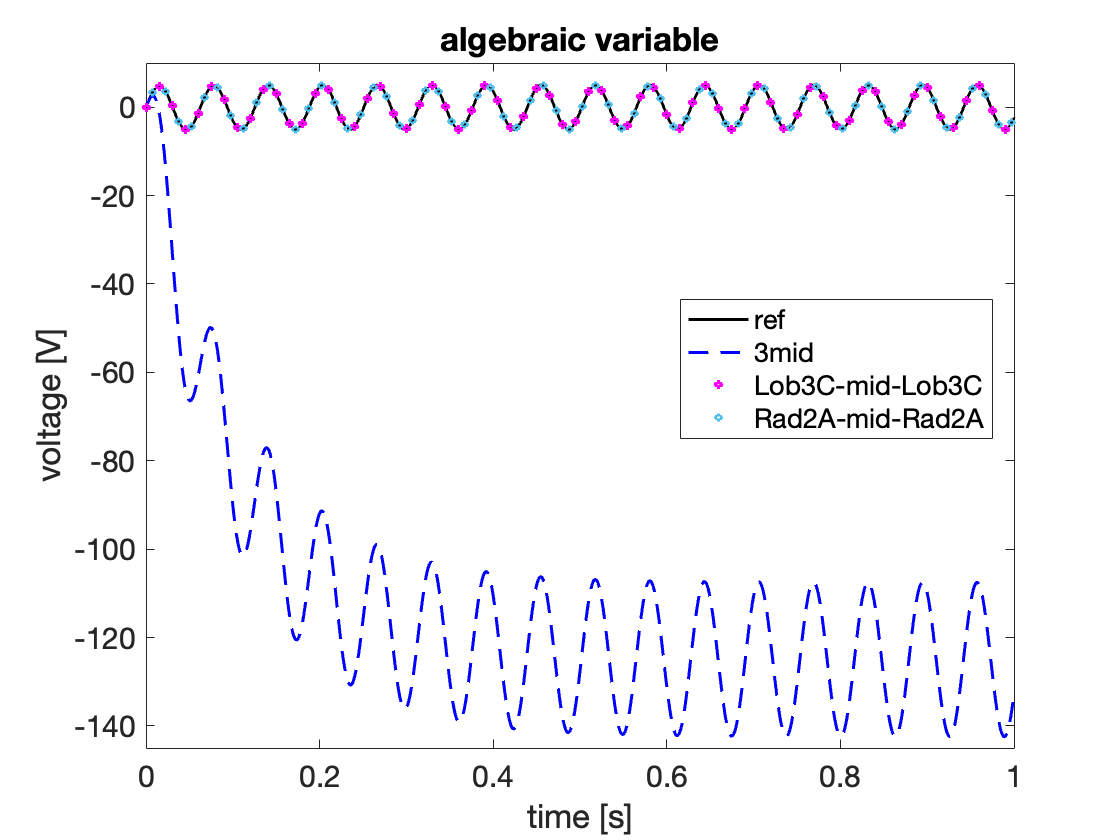}
    \caption{Example~\ref{ex:JR-DAE-R}. Top: Reference solution for differential variables $e_1$ (left) and $\jmath$ (right). Middle: Convergence behavior of differential variables (left) and algebraic variables (right) in Strang splitting with different flux approximations for relative step size $h$. Bottom: Evaluation of dissipation inequality for Strang splitting with some second-order flux approximations of fixed step size (left); reference solution and simulation results for algebraic variable $e_2$ (right).}
    \label{fig:exJR-DAE-R_results}
\end{figure}

In the following study, we use the parameters given in Fig.~\ref{fig:JR-network-index1}. As already shown in the previous examples, the energy-conserving implicit midpoint rule is an excellent scheme for the $J$-associated subsystem. As a 1-stage Gauss-Runge-Kutta method it theoretically has convergence order two for differential and algebraic variables for index-1 DAEs; for index-2 DAEs, convergence of the algebraic variables is no longer guaranteed, \cite{hairer1989}. 
The theoretical convergence order is also achieved for the original pH-system, but it turns out that the A-stable symmetric symplectic integrator is not suitable for the stiff $R$-associated subsystem in the splitting procedure. 
We observe not only a loss of order in the algebraic variable, but also an undesirably high error constant, cf.\ 3mid in Fig.~  \ref{fig:exJR-DAE-R_results}. As known from literature, the L-stable Lobatto-IIIC rules as well as the Radau family (Radau-IA and Radau-IIA rules) are very well suited for stiff systems and DAEs, the 2-stage variants analytically also promise second order in the algebraic variables. Applied to the dissipative $R$-associated subsystem, we obtain the desired second-order convergence in all variables in the Strang splitting and confirm our analytical findings, cf.\ Lob3c-mid-Lob3c, Rad1A-mid-Rad1A and Rad2A-mid-Rad2A in Fig.~  \ref{fig:exJR-DAE-R_results}. The implicit Euler method also belongs to the Radau family. Used for the dissipative subsystem, it leads to a convergence reduction of the Strang splitting to order one as expected (see iE-mid-iE). Approximations generated with  Lobatto-IIIC or Radau rules for the  $R$-associated subsystem and implicit midpoint rule for the $J$-associated subsystem satisfy the dissipation inequality and conserve the energy of the $J$-associated subsystem even for moderate step sizes.

%%%%%%%%%%%%%% EXAMPLE - Handling of Constraint %%%%%%%%%%%%%%%%%%%%%%%%

\begin{example}[Impact of Assumption~\ref{ass:restriction-index1}, regularized pHS]  \label{counter_example} \quad \end{example}
\begin{figure}[t]
\centerline{\includegraphics[width=0.65\textwidth,trim={0 1.1cm 0 0.5cm},clip]{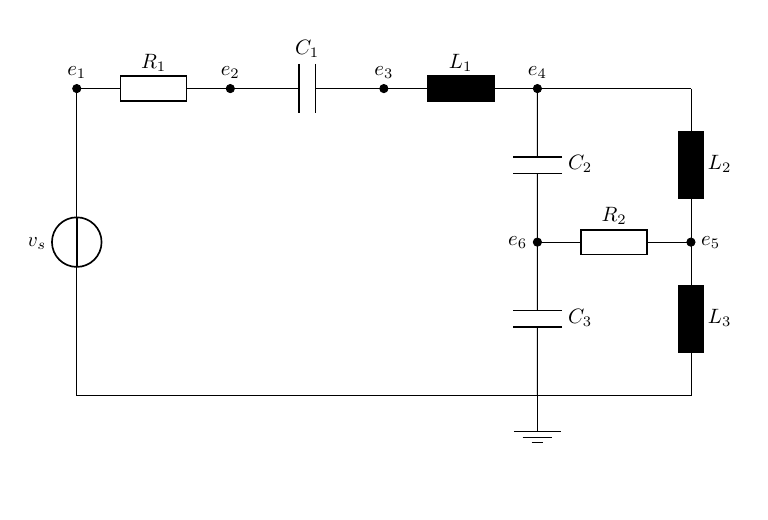}}
 \caption{Example~\ref{counter_example}: RLC-circuit  with voltage source $v_s(t) = \sin(10^{9}\,t)$~V, $t\in \mathbb{I}=[0,10^{-7}]$ and initial value $x(0)=0$. Parameters:  $C_1 = C_3= 10^{-12}$~F,  $C_2 = 5\cdot 10^{-13}$~F,   $R_1=R_2 =2\cdot 10^{-2}~\Omega$ as well as $L_1=L_2=L_3= 5\cdot 10^{-7}$~H; benchmark taken from \cite{OSM_DAE}.
 \label{fig:circuit}}
\end{figure} 

The small RLC-circuit presented in Fig.~\ref{fig:circuit} operates in a GHz regime as often used in chip design. Using the branch-oriented loop-cutset approach it can be modeled as an index-1 pH-DAE for the branch currents and voltages $x^\top=(\imath_{L_1},\imath_{L_2},\imath_{L_3},v_{C_1},v_{C_2},v_{C_3},v_{R_1},v_{R_2})$. The system matrices are
\begin{align*}
    J &= \begin{pmatrix}
         \phantom{-}0 &  \phantom{-}0 &  \phantom{-}0 &  \phantom{-}1 &  \phantom{-}1 &  \phantom{-}1 &  \phantom{-}1 &  \phantom{-}0 \\
         \phantom{-}0   &  \phantom{-}0 &  &  & -1 &  \phantom{-}0 &  \phantom{-}0 & \phantom{-}1 \\ 
           \phantom{-}0 &   &  \phantom{-}0 &  &  & -1 &  \phantom{-}0 & -1 \\ 
        -1 &   &   &  \phantom{-}0 &  &  &  &  \\
        -1 &  \phantom{-}1 &  &  &  \phantom{-}0 &  & &  \\ 
        -1 &  \phantom{-}0 &   \phantom{-}1 &  &  &  \phantom{-}0 &  &  \\
        -1 &  \phantom{-}0 &  \phantom{-}0 &  &  &  &  \phantom{-}0 &  \\
        \phantom{-}0& -1 &  \phantom{-}1 &  &  &  &  &  \phantom{-}0
    \end{pmatrix}, \qquad
    &&R = \begin{pmatrix}
        0 &    &  &  &  &  & &  \\
           & 0 &  &  &  & & &  \\
           &  & 0 &  &  & & &  \\
           &  &   & 0 & & &  &  \\
           &  &   &   & 0 &  &  &  \\
          &  &   &   &   & 0  & &  \\
          &  &   &   &   &   & \tfrac{1}{R_1} & \\
          &  &   &   &  &   &   & \tfrac{1}{R_2}
    \end{pmatrix},\\
     E &= \begin{pmatrix}
        L_1 &    &  &  &  &  & &  \\
           & L_2 &  &  &  & & &  \\
           &  & L_3 &  &  & & &  \\
           &  &   & C_1 & & &  &  \\
           &  &   &   & C_2 &  &  &  \\
          &  &   &   &   & C_3  & &  \\
          &  &   &   &   &   & 0 & \\
          &  &   &   &  &   &   & 0
    \end{pmatrix}, \qquad
    &&B = \begin{pmatrix}
    -1\\0\\0\\0\\0\\0\\0\\0
     \end{pmatrix},
\end{align*}
and the input function is prescribed by the voltage source $u(t)=v_s(t)$, $t\in \mathbb{I}$. The branch voltages $v_{R_1}$, $v_{R_2}$ are the algebraic variables to the constraints that occur here in the energy-conserving matrix  $J$ and in the dissipation matrix $R$. 

\begin{figure}[t]
  \includegraphics[width=0.525\textwidth]{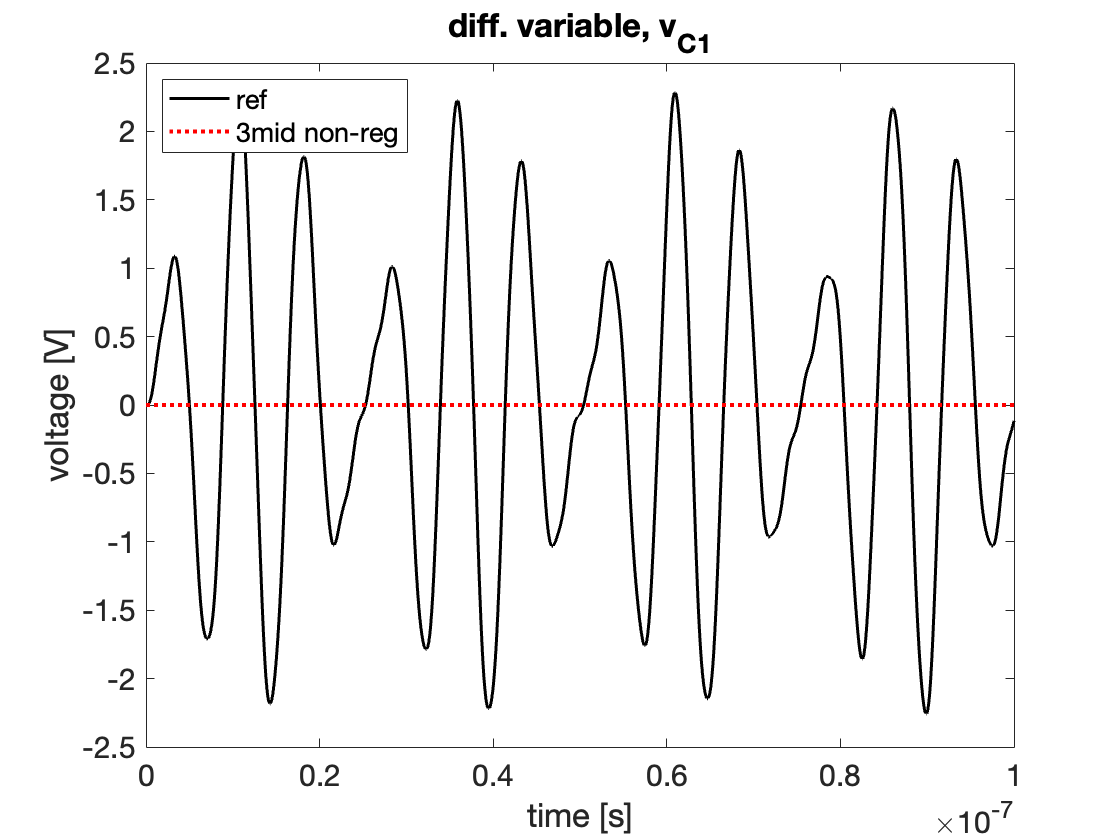}\hspace*{-0.5cm}
     \includegraphics[width=0.525\textwidth]{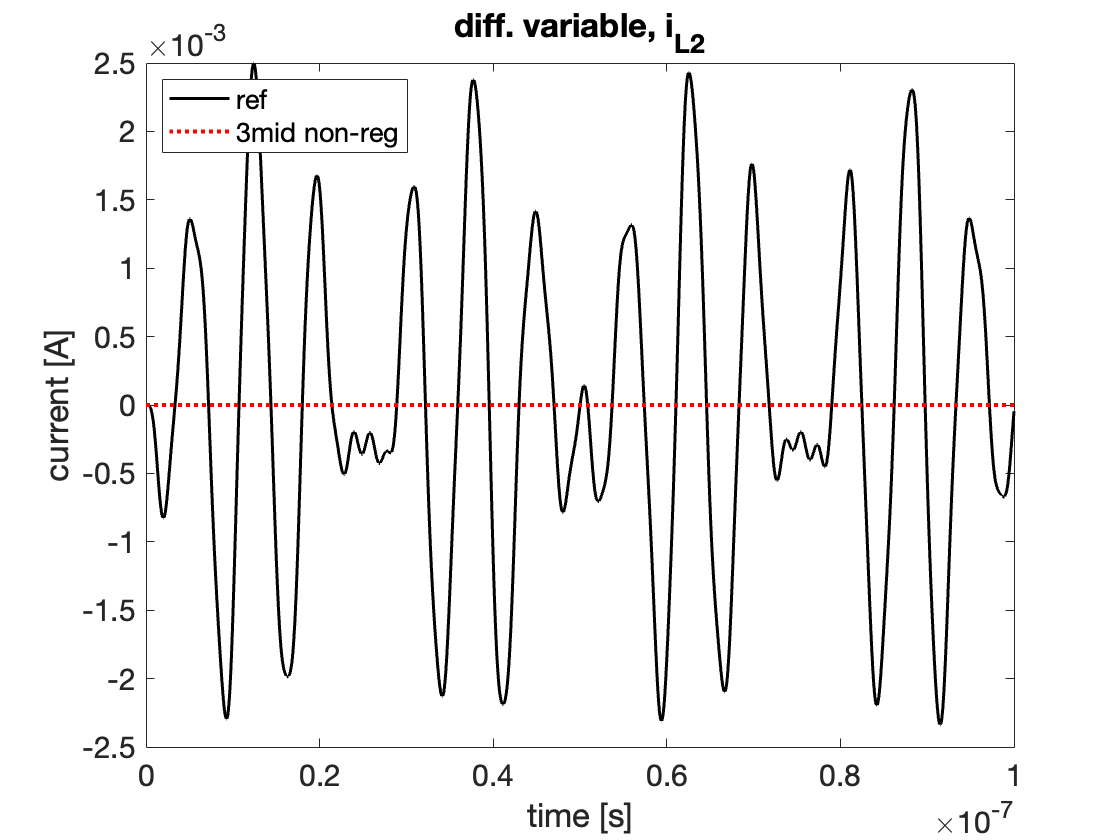}\\
      \includegraphics[width=0.525\textwidth]{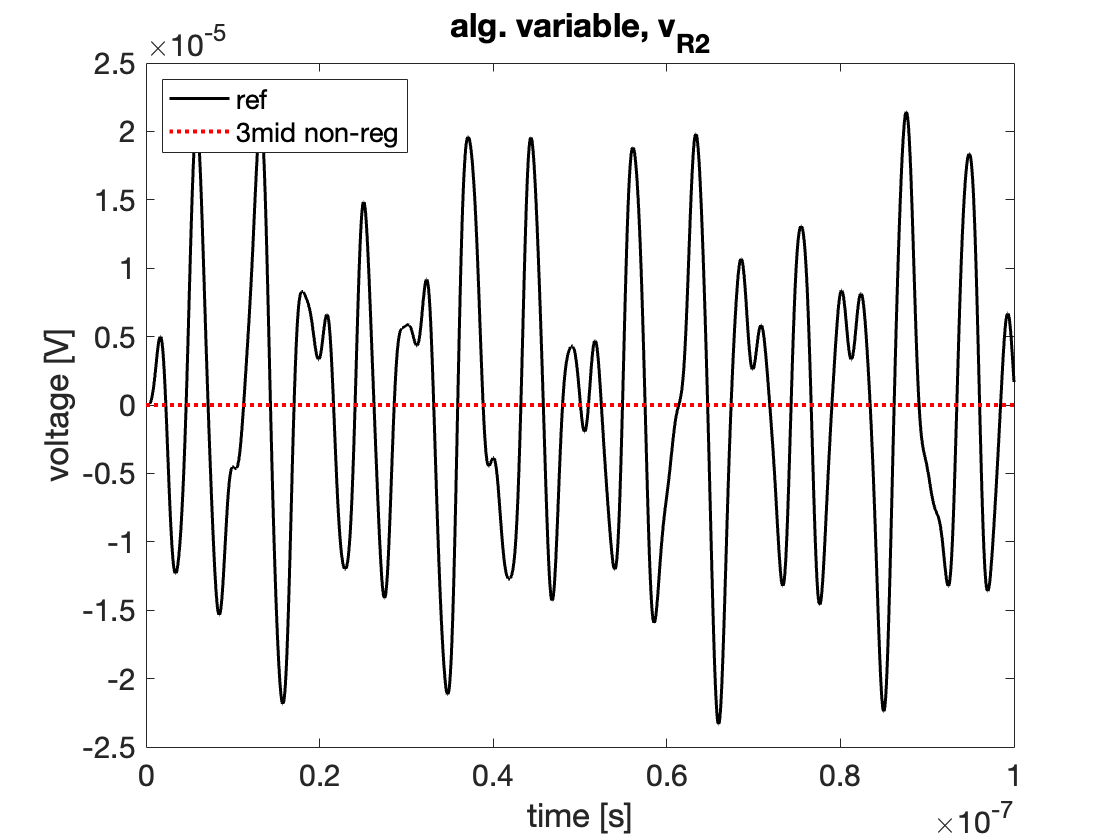}\hspace*{-0.5cm}
     \hfill\\
     \includegraphics[width=0.525\textwidth]{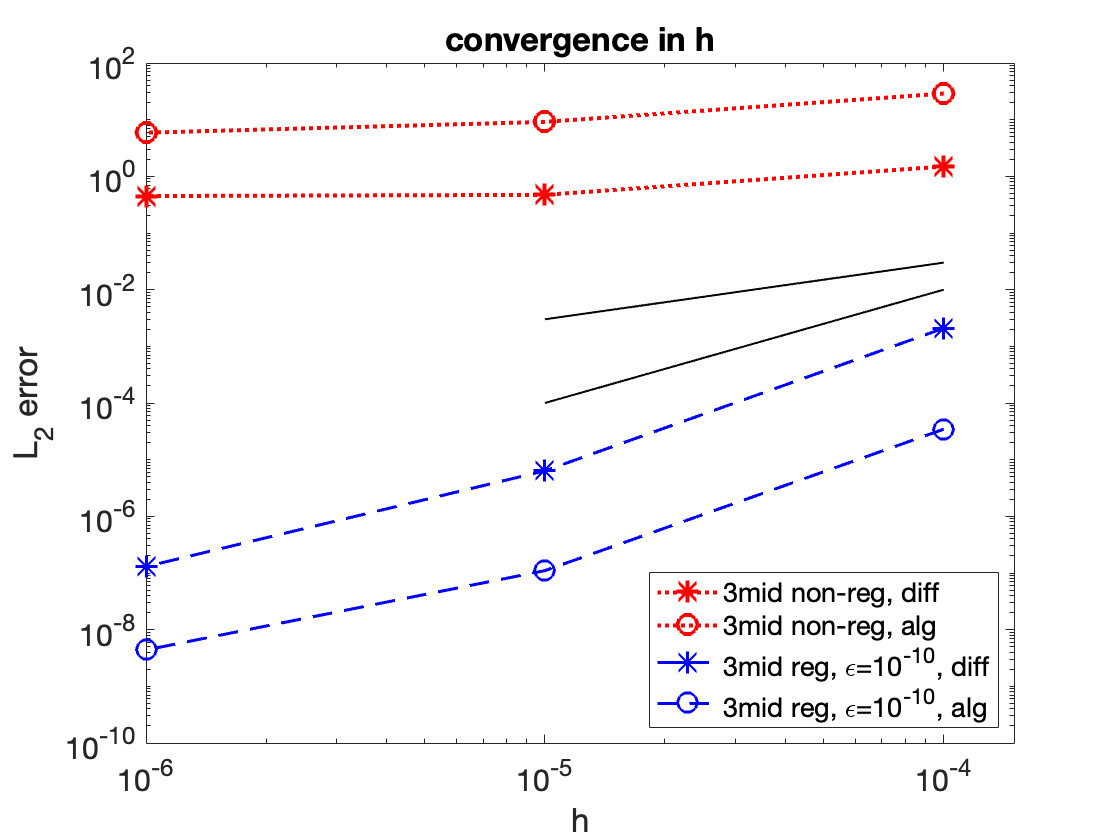}\hspace*{-0.5cm}
     \includegraphics[width=0.525\textwidth]{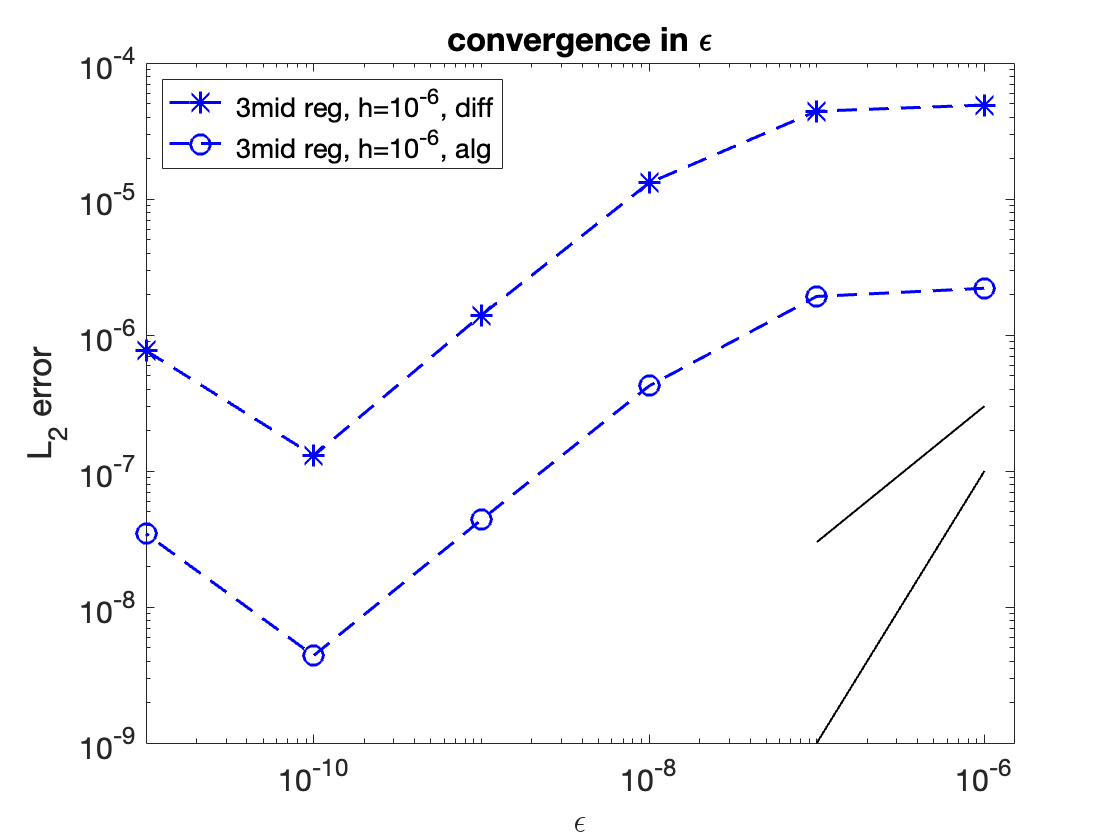}
 \caption{Example~\ref{counter_example}. Top: Reference solution and result of non-regularized splitting for differential variables $v_{C_1}$ (left) and $\imath_{L_2}$ (right). 
 Middle:  Reference solution and result of non-regularized splitting for algebraic variable $v_{R_2}$.
 Bottom: Convergence behavior in $h$ for non-regularized and regularized Strang splitting with $\epsilon=10^{-10}$ (left) and convergence in $\epsilon$ with fixed relative step size $h=10^{-6}$ (right).}
\label{fig:exJR-reg_results}
\end{figure}

The example does not satisfy Assumption~\ref{ass:restriction-index1}.  
If a $J$-$R$ decomposition is carried out regardless of this,  the resulting subsystems are
$$ E\dot x=Jx, \,\,\, \{E,J\} \text{ regular}, \qquad \qquad E\dot x=-Rx+Bu(t), \,\,\, \{E,R\} \text{ regular},$$
which splits the algebraic constraints and destroys the underlying structure. The interplay between decomposition and algebraic constraints lacks logical consistency. As a consequence, no splitting scheme converges, as demonstrated here for the Strang splitting.
The limitation can be overcome by help of regularization. Introducing $E_\epsilon=E+\epsilon K_E^\top K_E$ with $K_E$ a projector onto $\ker(E)$ and regularization parameter $0<\epsilon\ll1$, we apply the splitting scheme to the regularized subsystems 
$$ E_\epsilon \dot x_\epsilon=Jx_\epsilon, \qquad \qquad E_\epsilon \dot x_\epsilon=-Rx_\epsilon+Bu(t), \qquad \qquad E_\epsilon>0.$$
We observe convergence to the reference solution of the original index-1 pH-DAE as $\epsilon \rightarrow 0$.
Figure~\ref{fig:exJR-reg_results} shows the results for Strang splitting with implicit midpoint rule (3mid). We find the analytically predicted second order in differential and algebraic variables. Note that certainly also other discretization schemes can be used for the (regularized) ODE subsystems as discussed in the previous examples, but the outcome is here similar and hence not shown. 

The approximation quality of the regularized $J$-$R$ decomposition-based splitting approach relies on the asymptotic convergence of the solution $x_\epsilon$ of the $\epsilon$-regularized pH-ODE to the solution $x$ of the original index-1 pH-DAE as $\epsilon \rightarrow 0$. We observe a linear convergence in $\epsilon$. According to the regularization theory we find the best regularization parameter for the simulation in the kink of the L-curve (here $\epsilon=10^{-10}$), cf.\ Fig.~\ref{fig:exJR-reg_results} (bottom, right).
Note that the dissipation inequality only holds in the asymptotic result and cannot be expected to hold for fixed $\epsilon$, $\epsilon\neq 0$.

%------------------------------
\section{Conclusions} 
\label{sec:conclusion}
%------------------------------

In this work, we established two novel decomposition strategies in operator splitting for port-Hamil\-tonian systems. Splitting schemes on top of a dimension-reducing decomposition for coupled index-1 DAEs (with and without private index-2 variables) have the same convergence rate as known from the ODE case. The energy-associated $J$-$R$ decomposition is generally applicable to implicit pH-ODEs, and under certain assumptions---about the interactions between the energy parts and the constraints---to index-1 pH-DAEs. The splitting schemes keep thereby their convergence behavior, and the concept of generalized Cayley transforms ensures energy conservation for the numerical discretization. Since many port-Hamiltonian applications imply index-2 DAEs, it seems promising to extend the $J$-$R$ decomposition to this model class in future. This goes hand in hand with the development of structure-preserving higher order splitting schemes and respective Cayley transforms, since numerical schemes generally lose order of convergence in the algebraic variables for index-2 systems.

\appendix 

%%%%%%%%%%%%%%%%%%%%%% Appendix: Generalized Cayley %%%%%%%%

\section{Generalized Cayley Transforms} \label{app:Cayley}
This appendix provides some details to our novel concept of generalized Cayley transforms that is used in Section~\ref{sec:JR-splitting-time-discrete}
The following definition extends the familiar notion of a Cayley transform.

\begin{definition} Given $E,A\in \mathbb{R}^{n\times n} $ with $E^\top =E \ge 0$, and $E-A$ non-singular, we refer to the matrix operator 
\begin{align*}
C(E,A) := (E-A)^{-1}(E+A)
\end{align*}
as \textit{generalized Cayley transform}. 
\end{definition}

\begin{lemma}\label{lem:Gen-Cayleytrafo}
Assume that $E^\top =E \ge 0$, $A^{\top}=-A$, $E,A \in \mathbb{R}^{n \times n}$, 
and that the nullspaces satisfy $\nullsp(E) \cap \nullsp(A) = \{0\}$. 
Then the generalized Cayley transform $C(E,A)$ satisfies 
\begin{align}\label{eq:genCayleyequation}
C(E,A)^\top \cdot E \cdot C(E,A) = E, 
\end{align}
\end{lemma}

\begin{proof} We note that the assumption on the nullspaces implies that $E-A$ is non-singular, so the Cayley transform $C(E,A)$ is indeed defined.

Now, assume first that $ E$ is positive definite, $E > 0$ and let $E^{1/2}$ denote the matrix square root of $E$ which is non-singular because $E$ is. Since $C(E,A)^\top = (E-A)(E+A)^{-1}
$ we have
\begin{align}
C(E,A)^\top \,E\,  C(E,A)& = (E-A)(E+A)^{-1}\,E\,(E-A)^{-1}(E+A) \nonumber \\
&= (E-A)(E+A)^{-1}E^{1/2} \cdot E^{1/2}(E-A)^{-1}(E+A). \label{eq:Cayley_product}
\end{align}
Setting $K = E^{-1/2}AE^{-1/2}$, we obtain
\begin{eqnarray*}
(E-A)(E+A)^{-1}E^{1/2} &=& E^{1/2}(I-K)(I+K)^{-1}, \\
E^{1/2}(E-A)^{-1}(E+A) &=& (I-K)^{-1}(I+K)E^{1/2},
\end{eqnarray*}
which upon substitution into \eqref{eq:Cayley_product} gives
\[
C(E,A)^\top \,E \,C(E,A) = E^{1/2}(I-K)(I+K)^{-1} \cdot (I-K)^{-1}(I+K)E^{1/2} = E,
\]
as rational functions of matrices commute. This proves the assertion for the case $E > 0$. If $E$ is only semi-definite, define $E_\epsilon = E + \epsilon I$, which is positive definite for $\epsilon > 0$. From what we have already shown we have
\[
C(E_\epsilon,A)^\top E_\epsilon C(E_\epsilon,A) = E_\epsilon.
\]    
Letting $\epsilon \to 0$, by continuity, we obtain \eqref{eq:genCayleyequation} 
for the case of positive semi-definite $E$.
%\hfill $\Box$
\end{proof}

As immediate consequence, Lemma~\ref{lem:Gen-Cayleytrafo} yields the isometry in \eqref{eq:isometry} (conservation of energy), see Lemma~\ref{lem:Cayleytrafonorm} in Section~\ref{sec:JR-splitting-time-discrete}.
If the symmetric part $\tfrac{1}{2}(A+A^T)$ of $A$ is negative semi-definite, the generalized Cayley transform decreases the energy.
\begin{lemma}\label{lem:Gen-Cayleytrafo_diss}
Assume $E^\top =E \ge 0$,  $A+A^T \leq 0$ and that the nullspaces satisfy $\nullsp(E) \cap \nullsp(A) = \{0\}$. Then the generalized Cayley transform $C(E,A)$ satisfies 
\[
C(E,A)^\top \cdot E \cdot C(E,A) \leq E, 
\]  
and thus
\begin{align*} 
\|C(E,A)x\|_E  \leq \|x\|_E \quad \text{ for } x \in \mathbb{R}^n.
\end{align*} 
with $\|x\|_E=\langle x,x\rangle^{1/2}_E$ and  $\langle x,y\rangle_E =y^\top Ex$. 
\end{lemma} 

\begin{proof} 
Without further mentioning, our proof will at several places use the basic property of the positive semi-definite ordering which states that for square matrices $F$ and $G$  we have that $F\leq G$ is equivalent to $HFH^T \leq HGH^T$ for any non-singular matrix $H$. 

As in the proof of Lemma~\ref{lem:Gen-Cayleytrafo}, assume first that $ E > 0$. With $K = E^{-1/2}AE^{-1/2}$ we this time obtain 
\begin{align}
C(E,A)^\top E C(E,A)
&= (E+A)(E-A)^{-1}E(E-A)^{-1}(E+A) \nonumber \\
&= (E+A^T)(E-A^T)^{-1}E^{1/2} \cdot E^{1/2}(E-A)^{-1}(E+A) \nonumber \\
&= E^{1/2}(I+K^T)(I-K^T)^{-1} \cdot (I-K)^{-1}(I+K)E^{1/2}. \label{eq:Cayley_product_diss2}
\end{align}
Since $A+A^T \leq 0$, we have $K+K^T \leq 0$ and thus
\begin{align}
    \qquad & I + K^T + K + K^TK \, \leq \,  I - K^T - K + K^TK  \nonumber \\
   \Leftrightarrow \qquad &(I + K^T)(I+K)  \, \leq \,  (I - K^T)(I-K)  \nonumber \\
   \Leftrightarrow \qquad &(I-K^T)^{-1}(I+K^T)(I+K)(I-K)^{-1} \, \leq \, I.  \label{eq:leqI}
\end{align}
 Since rational functions of matrices commute, \eqref{eq:leqI} gives
\begin{align*}
(I+K^T)(I-K^T)^{-1} (I-K)^{-1}(I+K)
=  (I-K^T)^{-1}(I+K^T)(I+K)(I-K)^{-1}  
\, \leq \,  I,
\end{align*}
and inserting into \eqref{eq:Cayley_product_diss2} yields the sought result, 
$C(E,A)^\top \cdot E \cdot C(E,A) \leq E$. If $E$ is only positive semi-definite, we use the same continuity argument as in the proof of Lemma~\ref{lem:Gen-Cayleytrafo}.
%\hfill $\Box$
\end{proof}

%%%%%%%%%%%%%%%%%%%%%% Appendix: Splitting schemes and flux approximation %%%%%%%%
\section{Splitting Schemes and Flux Approximation}
\label{appendix:splitting-approximation}
Splitting is a powerful numerical tool to deal with dynamical systems that consists of subproblems with profoundly different behavior. There exists a large number of methods of different order for autonomous dynamical systems, see, e.g., \cite{hairer2006,mclachlan2002} and reference within. In case of a non-autononous system, the usual trick is to introduce a new time variable for every considered subsystem and then apply the standard schemes on the transformed autonomous equations. In the splitting scheme only the time variable of the considered subproblem is integrated, while the others are kept unchanged.

Consider a decomposition into two subproblems, i.e.,
\begin{align*}
\dot x=f(t,x)=\mathrm{f_1}(t,x)+\mathrm{f_2}(t,x), \qquad x(0)=x_0,
\end{align*}
where $\phi$ represents the exact flux of the original system and 
$\phi_i$, $i=1,2$, the exact (analytical) fluxes of the two subproblems, $\dot x=\mathrm{f}_i(x,t)$.
We use the notation $\phi^{t+h,t}(x)=x(t+h)$ where $\phi^{t,t}(x)=x$.

Well-established splitting schemes are the first-order Lie-Trotter and the second-order Strang splitting. Both preserve the dissipation inequality of port-Hamiltonian systems. The fourth-order Tripel Jump scheme results from a composition of the Strang splitting. It has negative step sizes, which even leave the integration domain $[t,t+h]$. For a non-autonomous system they are given by  \\[0.2ex]

\begin{tabular}{l l r}
Lie-Trotter &  $\Phi^{h,0}(x_0)=\phi_1^{h,0}\circ \phi_2^{h,0}(x_0)$ & $p=1$ \\
Strang  & $\Phi^{h,0}(x_0)=\phi_2^{h,h/2}\circ \phi_1^{h,0}\circ \phi_2^{h/2,0}(x_0)$ &  $p=2$\\
Triple Jump &  
$\Phi^{h,0}(x_0)=\phi_2^{(2\alpha+\beta)h,(3\alpha/2+\beta)h}\circ \phi_1^{(2\alpha+\beta)h,(\alpha+\beta) h} \circ ... $& $p=4$\\
& \multicolumn{2}{l}{\hspace*{1.5cm} $\circ\, \phi_2^{(3\alpha/2+\beta)h,(\alpha+\beta/2)h}\circ \phi_1^{(\alpha+\beta) h, \alpha h} \circ ...$}\\
& \hspace*{1.5cm} $ \circ\, \phi_2^{(\alpha+\beta/2)h,\alpha h/2}\circ \phi_1^{\alpha h,0} \circ \phi_2^{\alpha h/2,0}(x_0)$ &\\
& with $\alpha=1/(2-2^{1/3})$, $\beta=-2^{1/3}/(2-2^{1/3})$, thus $2\alpha+\beta=1$ 
\hspace*{-0.35cm}&
\end{tabular}\\ \quad

The numerical flux approximation (discretization method) has to be chosen with respect to the properties of the respective subproblem and the convergence order $p$ of the underlying splitting scheme. Table~\ref{tab:RK-conv} gives an overview on the convergence order and the stability behavior for the Runge-Kutta methods used in Section~\ref{sec:numerics}.

Note that the prescribed splitting procedure is of interest for non-autonomous systems if the time dependencies are cheap to compute. Otherwise the overall algorithm may be computationally costly, since the subproblems have to evaluated several times (number of stages) per time step. Thus splitting methods for non-autonomous systems were developed on the basis of the Magnus series, which show better efficiency \cite{blanes2006}.

\begin{table}[tb]
\begin{center}
    \begin{tabular}{l|c r|c|c c| c c| c}
        \textbf{method} & 
        \multicolumn{2}{|c|}{\textbf{stages}} & \textbf{ODE} & \multicolumn{2}{|c|}{\textbf{index-1}} 
         & \multicolumn{2}{|c|}{\textbf{index-2}} & \textbf{stability}\\
         & & & & diff. & alg. & diff. & alg. &\\
         \hline
         Gauss & $s$ & odd & $2s$ & $2s$ & $s+1$ & $s+1$ & $s-1$ & A\\
         &  & even & $2s$ & $2s$ & $s$ & $s$ & $s-2$ & A\\
         Radau IA & $s$& &$2s-1$&$2s-1$ & $s$ & $s$&$s-1$& L\\
         Radau IIA & $s$& &$2s-1$&$2s-1$ & $2s-1$ & $2s-1$&$s$& L\\
         Lobatto IIIC & $s$& &$2s-2$&$2s-2$ & $2s-2$ & $2s-2$&$s-1$ & L\\
         \hline
         Euler explicit & 1& & 1 & - & - &- &- &-\\
         Heun & 2 & & 2 & - & - &- &-&-\\
    \end{tabular}
    \end{center}
    \caption{Convergence order of Runge-Kutta methods for ODEs and DAEs in semi-explicit form with regard to differential and algebraic variables. The implicit Euler method corresponds to a Radau IIA method with $s=1$.}\label{tab:RK-conv}
\end{table}

%%%%%%%%%%%%%%%%%%%%%% Appendix: ELECTRIC NETWORK and $J$-$R$ splitting %%%%%%%%
\section{Energy-Associated $J$-$R$ Decomposition and Electric Network Models}
\label{appendix:circuits-and-JRsplitting}

In this appendix we discuss two standard modeling frameworks for electric networks and their appropriateness for the $J$-$R$ decomposition. We show that a suitable description as index-$1$ pH-DAE only exists  if the energy-conserving matrix $J$ does not contribute to the algebraic constraint, i.e., the $J$-associated subsystem is described by an ODE (Assumption~\ref{ass:restriction-index1}, case (b)).

%%%%%%%%%%%%
\subsection{Modified nodal analysis (MNA) modeling}\label{app:mna}

We consider a linear electric network with capacitances, inductances, resistances, and independent current $\imath(t)$, voltage sources $v(t)$, see, e.g., \cite{Guenther1999,Bartel2018}. Branch-node relations are specified via (reduced) incidence matrices, $A=[A_C, A_L, A_R, A_I, A_V] \in \{-1,0,1\}^{n\times b}$ (i.e., the ground node is removed from the incidences). Furthermore, let $A_V$ have full rank such that the circuit is physically sound. 
The linear MNA network equations for the node potentials $e$ and currents through inductances $\jmath_L$ as well as through voltage sources $\jmath_V$ imply the following linear system of pH-DAEs for $x^\top=(e,\jmath_L,\jmath_V)$,
\begin{align*}
\underbrace{\begin{pmatrix}
    A_C C A_C^\top & 0 & 0 \\
    0 & L & 0 \\
    0& 0 & 0 
\end{pmatrix}}_{\displaystyle E:=}
\begin{pmatrix}
    \dot e \\ \dot \jmath_L \\ \dot \jmath_V
\end{pmatrix}
 &=
\Biggl(
\underbrace{\begin{pmatrix}
    0 & - A_L & -A_V \\
    A_L^ \top & 0 & 0 \\
    A_V^\top & 0 &0 
\end{pmatrix}}_{\displaystyle J:=} -
\underbrace{\begin{pmatrix}
    A_R G A_R^\top & 0 & 0\\
    0 & 0 & 0 \\
    0 & 0 & 0
\end{pmatrix}}_{\displaystyle R:=}\Biggl)
\begin{pmatrix}
    e \\ \jmath_L \\  \jmath_V
\end{pmatrix}  \\ 
& \quad+ 
\underbrace{\begin{pmatrix}
    -A_I & \phantom{-}0 \\ \phantom{-}0 & \phantom{-}0 \\ \phantom{-}0 & -I
\end{pmatrix}}_{\displaystyle B:=}
\begin{pmatrix}
    \imath(t)  \\ v(t)
\end{pmatrix}
\end{align*}
with output
\[
y = \begin{pmatrix} -A_I^\top & \phantom{-}0 & \phantom{-}0 \\
                    \phantom{-}0    & \phantom{-}0 &-I \end{pmatrix} x 
                    = \begin{pmatrix} - A_I^\top e \\ - \jmath_V \end{pmatrix}.
\]
This system has maximal index~2. 

In the following investigations, we use a projector $K_E$ onto $\ker(E)$, which reads here
\begin{align*}
 K_E=\text{blkdiag}(K_C, 0, I)
\end{align*}
with $K_C$ being a projector onto $\ker (A_C^\top)$. 
Now, for the $J$-$R$ decomposition, we have the two cases according to Assumption~\ref{ass:restriction-index1}:
\begin{enumerate}
\item[(a)] $K_E^\top R = 0$ and $K_E^\top B = 0$; and the matrix pencil $\{E, J\}$  is non-singular (i.e., algebraic constraints are described by $J$),\\
or
\item[(b)] $K_E^\top J = 0$; and the matrix pencil $\{E, R\}$ is non-singular (i.e., algebraic constraints are described by $R$ and $Bu$).
\end{enumerate}

Case (a) is of particular interest to us, as it requires the use of generalized Cayley transforms. In this case, the source term $Bu$ only contributes to the differential part. This excludes voltage sources, since $v(t)$ would appear in an algebraic equation. Since $R$ also only contributes to the differential part, we have $A_R^\top K_C=0$ (resistors in differential part) and $K_C^\top A_I =0$ (current sources in differential part).  This means that the terminals of each resistor and each current source are connected to a path of capacitors, therefore rank$(A_C)=\text{rank}(A_C, A_R, A_I)$ (the number of connected units is equal). To obtain a DAE instead of an ODE, it is necessary that $\text{rank}(A_C) < n$. This means that there are more than one connected unit in the subsystem $(A_C,A_R,A_I)$. Hence, the inductors must form a cutset. This immediately results in a pH-DAE of index~2. 

Regarding Assumption~\ref{ass:restriction-index1}, index-1 pH-DAEs in MNA modeling belong to case (b), cf.\ Example~\ref{ex:JR-DAE-R}.

%%%%%%%%%%%
\subsection{Loop-cutset modeling}

In the following, we consider an electric network with the same linear lumped elements as in Appendix~\ref{app:mna}, but now modeled by the branch-oriented loop-cutset equations. Assuming that the index-1 topological conditions are satisfied, the system of equations reads
 \begin{align} \nonumber
        D \dot x(t) + Jx(t) + My(t) &= r_x(t) \\\label{eq:loopcutset}
        -M^\top x(t) + Sy(t) &= r_y(t)  \\
        z(t) + K_xx(t) + K_yy(t) &= r_z(t) \nonumber
 \end{align}
with the variables $x^\top = (\jmath_L, v_C)$, 
$y^\top=(\jmath_R  ,  \,v_G)$ and $z^\top=(v_I ,  \,\jmath_V)$, the input
${u}^\top(t):=(\imath_s (t), v_s (t))$ (independent current and voltage sources) and the matrices 
\begin{equation*}
    D=\begin{pmatrix}
       L & 0 \\ 0 & C
    \end{pmatrix}, \quad J = \begin{pmatrix}
       0 & -Q_{CL}^\top \\ Q_{CL} & 0
    \end{pmatrix}, \quad M= \begin{pmatrix}
       0 & -Q_{GL}^\top \\ Q_{CR} & 0 
    \end{pmatrix}, \quad S=\begin{pmatrix}
       R & -Q_{GR}^\top \\ Q_{GR} & G 
    \end{pmatrix}
\end{equation*}
(with $L,\, C,\, R,\, G$ symmetric positive definite) and 
\begin{equation*}
    K_x = \begin{pmatrix}
       0 & -Q_{CI}^\top \\ Q_{VL} & 0 
    \end{pmatrix}, \quad 
    K_y = \begin{pmatrix}
       0 & -Q_{GI}^\top \\ Q_{VR} & 0
    \end{pmatrix}, \quad
    K_z = \begin{pmatrix}
       0 & Q_{VI}^\top \\ -Q_{VI} & 0
    \end{pmatrix}, 
        \end{equation*}
    \begin{equation*}
    r_x 
    = K_x^\top {u},  
    ,\quad 
    r_y  
     = K_y^\top {u}
    , \quad 
    r_z = 
      K_z {u}                                      .
\end{equation*}
The variables $v_X$ and $\jmath_X$ represent the branch voltages and currents through a $X$-type component, respectively. Note that we refer to the resistors whose currents are determined in terms of the voltages by the equation $\jmath_G = Gv_G$  as $G$-resistances and to the resistors whose voltages are determined in terms of the current by the equation $v_R = R \jmath_R$ as $R$-resistances. Moreover, denoting by $Q = [I,  Q_l] \in \mathbb{R}^{(n-1) \times b}$ the fundamental cutset matrix, we divide the columns of the submatrix $Q_l$ with reference to the different circuit elements we have on twigs (tree branches) and links. Each submatrix is labeled $Q_{XY}$ (cf.~\cite{OSM_DAE} for detailed notations).

The DAE system \eqref{eq:loopcutset} fits into the port-Hamiltonian framework. To avoid confusion, we use here the subscript $\Bar{\,\,\,}$ to refer to the variables, output, input and coefficient matrices of the pHS, i.e.,
\begin{equation*}
    \Bar{x}^\top:= (x^\top\!\!,\,\, y^\top), \quad 
    \Bar{y}:=     -z 
    \quad \text{and} \quad \Bar{u}:=u.
\end{equation*}
Set 
\begin{align*}
    S_1 = \begin{pmatrix}
        R & 0 \\ 0 & G
    \end{pmatrix}\!, 
     \quad
    S_2= \begin{pmatrix}
        0 & Q_{GR}^\top \\ -Q_{GR} & 0
    \end{pmatrix},
\end{align*}
then  \eqref{eq:loopcutset} can be rewritten as
\begin{align*}
\underbrace{\begin{pmatrix}
    D  & 0 \\
    0  & 0 
\end{pmatrix}}_{\displaystyle \Bar E:=}
\begin{pmatrix}
    \dot x \\  \dot y
\end{pmatrix}
 &=
\Biggl[
\underbrace{\begin{pmatrix}
    -J & - M \\
    M^ \top & S_2 
\end{pmatrix}}_{\displaystyle \Bar J:=} -
\underbrace{\begin{pmatrix}
     0 & 0\\
    0 & S_1
\end{pmatrix}}_{\displaystyle \Bar R:=}
\Biggr]
\begin{pmatrix}
    x \\ y
\end{pmatrix} + 
\underbrace{
\begin{pmatrix}
    K_x^\top  \\ K_y^\top
\end{pmatrix}
}_{\displaystyle \Bar{B}:=}
 \bar{u},\\
 \bar{y} &=  -z =   \underbrace{(K_x,K_y)}_{\displaystyle \Bar B^\top:=} \Bar x +\underbrace{(- K_z)}_{\displaystyle \Bar N:=} \bar u.
\end{align*}
with $\Bar{E}$, $\Bar{R}$ symmetric positive semi-definite and $\Bar{J}$ skew-symmetric matrices.

\begin{figure}[t]
    \centerline{\includegraphics[trim={0 1.28cm 0 0},clip,width=0.4\textwidth]{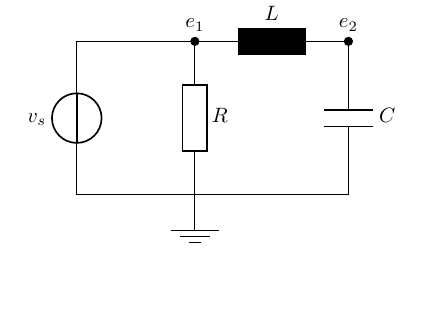}}
    \caption{Circuit with index-1 pH-DAE corresponding to case (b), cf.\ \eqref{eq:LCA_network}.}
    \label{fig:loopcutset}
\end{figure}

Since the resistive part belongs to the algebraic equation, we have to restrict $\Bar{J}$ to the differential equation for a valid $J$-$R$ decomposition. This corresponds to case (b) in Assumption~\ref{ass:restriction-index1}, just as with MNA modeling. 
Figure~\ref{fig:loopcutset} illustrates an example circuit whose loop-cutset equations are a respective index-1 pH-DAE corresponding to case (b), 
\begin{align}
   \label{eq:LCA_network}
    \begin{pmatrix}
        L & & \\
         & C & \\
         & & 0 
    \end{pmatrix} 
       \begin{pmatrix}
        \dot \jmath_L \\ \dot v_C \\ \dot \jmath_R
    \end{pmatrix} = \Biggl[
    \begin{pmatrix}
        0 & 1 & 0  \\
        -1 & 0 & 0 \\
        0 & 0 & 0
    \end{pmatrix} - \begin{pmatrix}
        0 & 0 & 0 \\
        0 & 0 & 0 \\ 
        0 & 0 & R
    \end{pmatrix} \Biggr]  \begin{pmatrix}
        \jmath_L \\ v_C \\\jmath_R
    \end{pmatrix}  + \begin{pmatrix}
        1 \\ 0 \\ 1
    \end{pmatrix} v_s(t).
\end{align}

%%%%%%%%%%%%%%%% REFERENCES %%%%%%%%%%%%%%%%%%%%%%%%%%%

\end{document}